\documentclass[3p, times]{elsarticle} 

\journal{}

\usepackage{graphicx}	
\usepackage{amsmath}	
\usepackage{amssymb}	
\usepackage{siunitx}
\usepackage{color, soul}

\newcommand{\mbf}[1]{\mathbf{#1}}			%
\newcommand{\x}{\mbf{x}}

\newcommand{\Q}{\mathbf{Q}}
\newcommand{\D}{\mathcal{D}}
\newcommand{\A}{\mathbf{A}}
\newcommand{\0}{\mathbf{0}}
\renewcommand{\u}{\mathbf{u}}

\newcommand{\q}{\mathbf{q}}
\newcommand{\F}{\mathbf{F}}

\newcommand{\Bn}{\mathcal{B}}
\newcommand{\Vn}{\mathcal{V}}
\newcommand{\f}{\mathbf{f}}
\newcommand{\g}{\mathbf{g}}

\newcommand{\B}{\mathbf{B}}
\newcommand{\de}[2]{\frac {\partial #1}{\partial#2}}

\renewcommand{\S}{\mathbf{S}}
\newcommand{\be}{\begin{equation} \begin{aligned} }
		\newcommand{\ee}{\end{aligned} \end{equation}}

\renewcommand{\epsilon}{\varepsilon}
\renewcommand{\phi}{\varphi}

\begin{document}

\begin{frontmatter}   

\title{A well balanced diffuse interface method for complex nonhydrostatic free surface flows}

\author[address1]{Elena Gaburro}
\ead{elena.gaburro@unitn.it}
\author[address2]{Manuel J. Castro}
\ead{castro@anamat.cie.uma.es}
\author[address1]{Michael Dumbser$^{*}$}
\ead{michael.dumbser@unitn.it}
\cortext[cor1]{Corresponding author}

\address[address1]{Department of Civil, Environmental and Mechanical Engineering, University of Trento, Via Mesiano, 77 - 38123 Trento, Italy.}
\address[address2]{Department of Mathematical Analysis, Statistics and Applied Mathematics, University of M\'alaga, Campus de Teatinos, 29071 M\'alaga, Spain.}

\begin{abstract}
In this paper we propose an efficient second order accurate well balanced finite volume method for modeling complex free surface flows
at the aid of a simple \textit{diffuse interface method}. 
The employed physical model is a two-phase model directly derived from the Baer-Nunziato system for compressible multi-phase flows. In particular, as proposed for the first time in \cite{dumbser2011simple}, the number of equations is reduced from seven to three by assuming that the relative pressure of the gas with respect to the atmospheric reference pressure is zero, and that the gas momentum is negligible  compared to the one of the liquid. 
The two-phase model does not make any of the classical assumptions of shallow water type systems, hence it does \textit{not} neglect  vertical accelerations and the free surface is \textit{not} constraint to be a single-valued function, so even \textit{complex} shapes 
as those of breaking waves can be properly captured. 

The resulting PDE system is solved by a novel \textit{well balanced} second order accurate path-conservative finite volume method on structured Cartesian grids, which is able to preserve \textit{exactly} the equilibrium states even in the presence of obstacles.  It furthermore automatically computes the location of the water-air interfaces, and assures low numerical dissipation at the free surface  thanks to a novel Osher-Romberg-type approximate Riemann solver. 
Finally, high computational performance is guaranteed by an efficient \textit{parallel} implementation on a GPU-based platform that 
reaches the efficiency of twenty million of volumes processed per seconds and makes it possible to employ even very fine meshes. 
The validation of our new well balanced scheme is carried out by comparing the obtained numerical results against existing 
analytical, numerical and experimental reference solutions for a large number of test cases, among which oscillating 
elliptical drops, dambreak problems, breaking waves, over topping weir flows, and wave impact problems. 
\end{abstract}


\end{frontmatter}



\section{Introduction}

Many important applications in environmental engineering and geophysics, such as the design of hydraulic structures for flood management (dams, harbors and levees) and even in mechanical, industrial and aerospace engineering, such as high speed water jet or high pressure fuel injection processes, are strictly related to the study of complex free surface hydrodynamics, often involving non-hydrostatic free surface 
flows, as well as flows where the free surface cannot be described by a single-valued function. 
The most commonly used physical models, in particular in the case of quasi incompressible flows, are depth-averaged hydrostatic shallow water type models.
They are derived from the incompressible Navier-Stokes equations by supposing that the horizontal length scale is much larger than the vertical one, hence any vertical acceleration can be neglected and a depth-averaging procedure is performed to reduce the number of involved space-dimensions by one. {Thus, the computational cost is reduced by one order of magnitude with respect to the corresponding two- and three-dimensional models, which is one of its principal advantages together with the explicit knowledge of the free surface position.} 
Their accurate and efficient numerical solution is a very active front of research in applied mathematics; it is out of the scope to  mention here all the relevant contributions to the field, we just mention a few of them like \cite{Bermudez1994, GarciaNavarro1, GarciaNavarro2,toro-book-swe,XingShu,Rhebergen2008} and \cite{Casulli1990,CasulliCheng1992,Casulli1999,CasulliWalters2000,Casulli2009,CasulliStelling2011,Giraldo2002,TumoloBonaventuraRestelli,GiraldoRestelli}. 

Moreover, we would like to underline that a very important topic in this context is the capability of the scheme in maintaining a certain class of \textit{equilibrium} solutions of the system \textit{exactly} up to machine precision. In the framework of Godunov-type schemes, 
this is not trivial and requires a special treatment even in the case of the simple steady state associated with water at rest. 
It is made possible by the introduction of the so-called \textit{well balanced} schemes that have been adapted to many different situations: multiple space dimensions, high order of accuracy, Cartesian, cylindrical and spherical geometries, and even to the Lagragian framework, see for example \cite{LeVequeWB, audusse2004fast,  Noelle1, pares2006numerical, Castro2006, Castro2008,  NCPRICE1D, NCPRICE2D, gaburro2016direct, castro_ortega2017well} and the reference therein for a more complete bibliography.

However, although the shallow water equations are properly applicable in many cases, their assumptions fail when unsteady phenomena with important vertical accelerations, such as dambreak problems and wave impact against solid structures have to be simulated. The problems become even more severe when interested in breaking waves or very general free surface shapes, which cannot be described with single-valued  functions any more, and finally when a compressible flow description has to be taken into account for higher Mach numbers.

For these reasons, in the present paper we have chosen a weakly compressible approach that can be directly derived from the more complete Baer-Nunziato model of compressible multiphase flows. In this model, the free surface is automatically captured by the volume fraction  function of the liquid, which is evolved in time by a scalar advection equation and distinguishes the volume occupied by the liquid phase from the one occupied by the gas phase.  
This model was first introduced in two space dimensions in \cite{dumbser2011simple} and then extended to the three dimensional case in~\cite{dumbser2013diffuse}. {The model we are considering here is  fully two-dimensional and it even allows to capture multivariate  free surface profiles. This makes its numerical simulation much more expensive with respect to the one-dimensional shallow water equations. To partially mitigate this problem and to keep computational times still acceptable, we propose an efficient CUDA implementation of our numerical method.} 
The use of a volume fraction function makes the approach similar to the well-known 
volume of fluid (VOF) method, see \cite{HirtNichols,Rieber,Kleefsman}, but in a compressible setting. 
In both cases the reliability of the proposed physical model has been shown by comparing the obtained
numerical results against available analytical, numerical and experimental reference solutions. The improved quality of a more complete model with respect 
to the simpler shallow water equations has been emphasized, in particular for the simulation of nonhydrostatic effects occurring  immediately after a dambreak, within breaking and overtopping waves and in recirculation phenomena. 
In the above references, the two-phase flow model has been solved with a third order ADER-WENO finite volume scheme based on the unstructured WENO reconstruction operator detailed in \cite{DumbserKaeser07}; the obtained results 
were in good agreement with the reference solutions, but the proposed numerical method was \textit{not} exactly well balanced. 

Thus, the aim of the present work is the design of new path-conservative algorithm that is  \textit{exactly well balanced} for the two  phase flow model introduced in \cite{dumbser2011simple,dumbser2013diffuse}, in such a way that the final scheme can be very accurate and  complete in the case of complex free surface flows and at the meantime can preserve steady equilibrium solutions up to machine precision.  
To preserve the equilibria in a system of equations with source terms, following \cite{Pares2006,castro2007well,MuellerWB}, we will rewrite it in terms of non-conservative products obtaining a system of the form
\be
\label{eq.generalform}
\de{\Q}{t} + \nabla \cdot \F(\Q)  + \B(\Q) \cdot \nabla \Q = \mathbf{S}(\Q), \quad \x \in \Omega \subset \mathbb{R}^2.
\ee
In this system $\x$ is the spatial position vector, $t$ represents the time, $\Omega$ is the considered computational domain, $\Q = (q_1,q_2, \dots, q_{\nu})$ is the vector of the conserved variables defined in the space of the admissible states $\Omega_{\Q} \subset \mathbb{R}^{\nu}$, $ \F(\Q) = (\,\f(\Q), \g(\Q)\,) $ is the non linear flux tensor, $\B(\Q) = (\, \B_1(\Q), \B_2(\Q) \,) $ is a matrix collecting the non-conservative terms, and $\mathbf{S}(\Q)$ represents a non linear algebraic source term. The system \eqref{eq.generalform} can also be written in the following quasi-linear form
\be 
\label{eq.quasilinear}
\de{\Q}{t}  + \A(\Q) \cdot \nabla \Q = \mathbf{S}(\Q), \quad \x \in \Omega(t) \subset \mathbb{R}^2, 
\ee 
with the system matrix $\A(\Q) = \partial \F / \partial \Q + \B(\Q)$. The system is \textit{hyperbolic} if for any normal direction $\mathbf{n} \neq 0$ the matrix $\A(\Q) \cdot \mathbf{n}$ has $\nu$ real 
eigenvalues and a full set of $\nu$ linearly independent eigenvectors for all $\Q \in \Omega_{\Q}$. 
The main difficulty of systems written in this form, both from the theoretical and the numerical points of view, comes from the presence of non-conservative
products that do not make sense in the standard framework of distributions when the solution $\Q$ develops discontinuities.
From the theoretical point of view, in this paper we assume the definition of non-conservative products as Borel measures given in \cite*{DalMaso1995}. 
This definition, which depends on the choice of a family of paths in the phase space $\Omega_{\Q}$, allows one to give a rigorous definition of weak solutions of (\ref{eq.generalform}).
We consider here the discretization of system (\ref{eq.generalform}) by means of numerical schemes which are \textit{path-conservative} in the sense introduced in \cite{Pares2006}. 
The concept of a path-conservative method, which is also based on a prescribed family of paths, provides a generalization of conservative schemes introduced by Lax for systems of conservation laws.
Then we will apply the recently introduced Osher-Romberg numerical flux \cite{gaburro2016direct, gaburro2018well} which is a general 
and very accurate well balanced technique. Provided the full set of eigenvalues and eigenvectors of $\A$ are known, it can be easily 
adapted to very different systems of equations and families of equilibrium solutions. 

As already mentioned above, in order to increase the computational efficiency of our numerical method, we have decided to use a massively parallel implementation, choosing in particular a GPU-based platform. 
We employ the NVIDIA CUDA (Compute Unified Device Architecture) framework, which is a hardware and software platform that allows to easily exploit the capacities of NVIDIA GPUs for numerical algorithms. We refer to the recent book \cite{mantas2016introduction} for a complete introduction to CUDA for Scientific Computing.

The rest of the paper is organized as follows. In Section~\ref{sec.PhysicalModel} we first introduce our mathematical model given by a special case of the more general Baer-Nunziato system for compressible multi-phase flows and the appropriate simplifications that allow us to reduce the number of equations from seven to three.
Then in Section~\ref{sec.C_GPU.NumericalMethod} and \ref{sec.C_GPU.CUDA} we present our second order accurate well balanced  path-conservative Osher-Romberg scheme in the particular case of a Cartesian mesh and of a parallel implementation in CUDA. 
Finally, in Section~\ref{sec.C_GPU.NumResults} we show equilibrium test cases together with a large variety of typical free surface 
flow problems and we exhibit the high efficiency of our code. The article is closed by some conclusive remarks and an outlook to 
future research.

\section{Physical model for complex free surface flow}
\label{sec.PhysicalModel} 

The simple and efficient two-phase interface-capturing algorithm proposed in this article is given by a special case of the more general Baer-Nunziato model for compressible multi-phase flows introduced for the first time by Baer and Nunziato in \cite{baer1986two}.
In this section, we present the physical model starting from the original Baer-Nunziato system and by introducing some appropriate simplifications that allow us to reduce the number of equations from seven to three, as was done for the first time in \cite{dumbser2011simple}. Then, we rewrite the sources via non-conservative products so that all the terms connected with the equilibrium could be treated together.

The original Baer-Nunziato system with gravity effects reads 
\be
\begin{cases}
	& \partial_t\left( \alpha_1 \rho_1\right)            + \nabla \cdot \left( \alpha_1 \rho_1 \mbf{u}_1 \right) = 0 \\
	& \partial_t\left( \alpha_1 \rho_1 \mbf{u}_1 \right) + \nabla \cdot \left( \alpha_1 \left( \rho_1 \mbf{u}_1 \otimes \mbf{u}_1 + p_1 \mbf{I} \right) \right) = p_I\nabla\alpha_1 + \alpha_1 \rho_1 \mbf{g} \\
	& \partial_t\left( \alpha_1 \rho_1 E_1 \right)       + \nabla \cdot \left( \alpha_1 \mbf{u}_1\left( \rho_1 E_1  + p_1  \right) \right) = - p_I \partial_t \alpha_1  \\
	& \partial_t\left( \alpha_2 \rho_2\right)            + \nabla \cdot \left( \alpha_2 \rho_2 \mbf{u}_2 \right) = 0 \\
	& \partial_t\left( \alpha_2 \rho_2 \mbf{u}_2 \right) + \nabla \cdot \left( \alpha_2 \left( \rho_2 \mbf{u}_2 \otimes \mbf{u}_2 + p_2 \mbf{I} \right) \right) = p_I\nabla\alpha_2 + \alpha_2 \rho_2 \mbf{g} \\
	& \partial_t\left( \alpha_2 \rho_2 E_2 \right)       + \nabla \cdot \left( \alpha_2 \mbf{u}_2\left( \rho_2 E_2  + p_2 \right) \right) = - p_I \partial_t \alpha_2  \\
	& \partial_t \alpha_1 + \mbf{u}_I \cdot \nabla \alpha_1  = 0, \\
\end{cases}
\ee 
where $\alpha_j$ is the volume fraction of phase number $j$ with $\alpha_1 +\alpha_2 = 1$, $\rho_j$ is the fluid mass density, $\mbf{u}_j =(u_j,v_j)$ the velocity vector, $p_j$ the pressure and $\rho_j E_j$ the total energy per mass unit of phase number $j$, respectively. Moreover, $\mbf{g}=(0,-g)$ is the vector of gravity acceleration, $g=9.81$.
The model must be closed by the equations of state (EOS) for each phase that link the pressure $p_j$ to the density and the internal energy, and furthermore the model requires a proper choice of the interface velocity $\mbf{u}_I$ and the interface pressure  $p_I$.
Baer and Nunziato proposed the following choice 
\begin{equation}
	\label{eq.interfaceLaw}
	p_I = p_2, \ \text{and } \ \mbf{u}_I=\mbf{u}_1,
\end{equation}
which we also use here.

As explained in detail in \cite{dumbser2011simple}, the three-equation~model we will study is based on the following \textit{simplifications}.
The first assumption is that all pressures are relative pressures with respect to the atmospheric reference pressure $p_0$, hence we can define $p_0 = 0$. 
Second, the gas surrounding the liquid is supposed to remain always at atmospheric reference conditions, i.e. the gas pressure is
\be
\label{eq.interfacePressure}
p_2 = p_0 = 0 = \textit{const},
\ee
which is a standard assumption for free surface flows in fluid mechanics, see for example \cite{truckenbrodt1980bd}.
It is based on the fact that for low Mach number flows the pressure fluctuations $p_0$ are
approximately proportional to $\rho\mbf{u}^2$ according to Bernoulli's law and
since the liquid density is several orders of magnitude larger than
the density of the gas ($\rho_1 \gg \rho_2$) the pressure fluctuations of the
gas phase are much smaller than the pressure fluctuations in the liquid
phase ($p'_2 \ll p'_1$). We therefore can neglect all evolution equations
related to the gas phase $j = 2$. 
Furthermore, according to \eqref{eq.interfaceLaw} and \eqref{eq.interfacePressure} the interface pressure automatically results $p_I = p_2 = p_0 = 0$. 
This is consistent with the usual standard assumption for free surface flows, where, at the
free surface of the liquid, atmospheric reference pressure boundary
conditions are imposed. Also the choice of the interface velocity
$ \mbf{u}_I = \mbf{u}_1$ according to \eqref{eq.interfaceLaw} is consistent, since the interface will obviously
propagate with the speed of the liquid phase. 
Third, the pressure of the liquid is computed by the Tait equation of state \cite{batchelor1967introduction}: the key idea is that according to the first
assumption it yields a relative pressure with respect to the atmospheric
reference pressure ($p_0 = 0$). We therefore have
\be
\label{eq.TaitLaw}
p_1 = k_0 \left( \left(\frac{\rho_1}{\rho_0}\right )^\gamma - 1 \right),
\ee 
where $k_0$ is a constant that governs the compressibility of the fluid
and hence the speed of sound, $\rho_1$ is the liquid density, $\rho_0$ is the liquid
reference density at atmospheric standard conditions and $\gamma$ is a
parameter that is used to fit the EOS with experimental data. Since
the EOS \eqref{eq.TaitLaw} does not depend explicitly on the internal energy, also
the liquid energy equation can be dropped. 

So, we obtain the following reduced \textit{three-equation model}
\be
\label{eq.ReducedBNmodel}
\begin{cases} 
	&\partial_t \left ( \alpha \rho \right ) + \nabla \cdot \left ( \alpha \rho \mbf{u} \right) = 0 \\
	&\partial_t \left ( \alpha \rho \mbf{u} \right ) + \nabla \cdot \left ( \alpha \left ( \rho \mbf{u}\otimes\mbf{u} + p \mbf{I} \right) \right) = \alpha \rho \mbf{g} \\
	&\partial_t \alpha + \mbf{u} \cdot \nabla \alpha = 0,
\end{cases}
\ee 
where indices have been dropped to simplify the notation. We work in a two-dimensional framework $\x= (x,y)$ where $y$ indicates the gravity direction.
We underline that the free surface is captured automatically by the volume fraction function $\alpha$ which is evolved in time by the last equation (advection equation) and allows us to distinguish between the portion of the domain occupied by the liquid $\Omega_\ell$, with $\alpha \sim 1$, and the one occupied by the gas $\Omega_a$, with $\alpha \sim 0$.

\subsection{Equilibria and non-conservative formulation}

An important \textit{family of equilibria} of \eqref{eq.ReducedBNmodel} consists in the water at rest solution given by 
\be
\label{eq.BNs_eq1}
\mbf{u}=\0, \quad \alpha= \text{ const},
\ee
and density and pressure obtained through the momentum equation in $y-$direction as follows
\be
\label{eq.equilibriumBNsimp}
\de{p}{y} = \de{}{y} k_0 \left( \left( \frac{\rho\left(y\right)}{\rho_0} \right ) - 1 \right ) = - \rho\left(y\right) g,
\ee 
which for the particular case $\gamma=1$ and $\rho(0) = \rho_0$ gives the simple solution
\be
\label{eq.BNs_eq3}
\rho\left(y\right) = \rho_0 \text{exp}\left( - \frac{g \rho_0}{k_0} \left(y-y_0\right) \right),
\ee 
and hence the pressure distribution
\be 
\label{eq.BNs_eq4}
p(y) = k_0 \left( \text{exp} \left( - \frac{g \rho_0}{k_0}\left(y-y_0\right) \right) - 1 \right),
\ee 
where $y_0$ is the free surface position.
Note that, since the constant $k_0$ is supposed to be very large and hence the argument of the exponential function is small, from the Taylor series expansion of $p$ 
\be
p(y) = k_0 \, \Biggl( 1&-\left(\frac{g \rho_0}{k_0}\left(y-y_0\right)\right) + \frac{1}{2} \left(\frac{g \rho_0}{k_0}\left(y-y_0\right)\right)^2  
-\frac{1}{6} \left(\frac{g \rho_0}{k_0}\left(y-y_0\right)\right)^3 + \dots -1 \Biggr),
\ee
we can deduce that in the limit $k_0 \rightarrow \infty$ 
\be 
\rho(y) = \rho_0, \quad p(y) \rightarrow -\rho_0 g \left(y-y_0\right),
\ee 
i.e. the pressure distribution tends to the hydrostatic one.

To be able to preserve the equilibrium \eqref{eq.equilibriumBNsimp}, that involves both terms in the flux and in the source, we need to rewrite all the involved terms in such a way that they could be treated together.
{For this purpose, the gravitational source term $\alpha \rho \mathbf{g}$ is converted into a differential term by rewriting it as the gradient of a gravitational
potential $Y(y) = y$, i.e. 
\be
\alpha \rho \mathbf{g} = \alpha \rho g \, \nabla {y}, 
\ee 
and introducing it into the advective operator on the left hand side of the system.
Note that now the system can be still written in terms of conserved variables through the introduction of $y$ as a new quantity, using the following trivial equation
\be
\partial_t{y} = 0, 
\ee \\
which is the typical strategy adopted also in \cite{greenberg1997analysis, gosse2000well, gosse2001well, castro2007well}. 
Indeed, once the steady state solution has been described through the balancing of fluxes and  non-conservative terms, the problem of its \textit{exact preservation} is transfered to the choice of the correct \textit{path} for writing the path-conservative method, and then one can apply some of the solutions available in the literature devoted to this topic.
Furthermore, we remark that the necessity of a path-conservative approach was however conveyed by the presence of the non-conservative advection equation for the volume fraction function $\alpha$; and finally, that the choice of considering the gravitational potential $Y$ (instead of, for example, the other space variable $x$) is guided by the aim of obtaining a balancing between forces acting in the gravity direction. }

{Our final non-conservative system therefore reads as follows} 
\be
\label{eq.ReducedBNmodel_nc}
\begin{cases} 
	&\partial_t \left ( \alpha \rho \right ) + \nabla \cdot \left ( \alpha \rho \mbf{u} \right) = 0 \\
	&\partial_t \left ( \alpha \rho u \right ) + \partial_x\left( \alpha  \rho u^2  + \alpha p\right) + \partial_y \left( \alpha  \rho uv \right)  = 0 \\
	&\partial_t \left ( \alpha \rho v \right ) + \partial_x\left( \alpha  \rho u v \right)  + \partial_y\left( \alpha  \rho v^2 \right)  + \partial_y\left(\alpha p\right)  + \alpha \rho g \, \partial_y{y} = 0 \\
	&\partial_t \alpha + \mbf{u} \cdot \nabla \alpha = 0,\\
	& \partial_t y = 0,
\end{cases}
\ee 
which can be cast in the general form \eqref{eq.generalform} with
\be
\label{eq.BNsimplified_genForm}
&\Q  \!=\! \left( \begin{array}{c} \alpha \rho\\ \alpha \rho u \\ \alpha \rho v  \\ \alpha \\ z              \end{array} \right)\!, 
\ \ \f \!=\! \left( \begin{array}{c}  \alpha \rho u \\  \alpha \rho u^2 + \alpha p\\   \alpha \rho u v  \\ 0 \\ 0    \end{array} \right)\!, 
\ \ \g \!=\! \left( \begin{array}{c} \alpha \rho v \\\alpha \rho u v \\ \alpha \rho v^2  \\ 0 \\ 0  \end{array} \right)\!, 
\ \ \B_1 \cdot \nabla \Q\!=\! \left( \begin{array}{c}  0 \\ 0 \\ 0 \\\!\! u\de{\alpha}{x}\!\! \\ 0 \end{array} \right)\!, \ \ \B_2 \cdot \nabla \Q \!=\! \left( \begin{array}{c}  0 \\ 0 \\ \!\!\! \de{\alpha p}{y} + \alpha \rho g \de{y}{y} \!\!\!\\ v\de{\alpha}{y} \\ 0 \end{array} \right)\!, \ \ \mathbf{S}\!=\!\0.  
\ee

\section{Well balanced path-conservative scheme}
\label{sec.C_GPU.NumericalMethod}

We solve the two-phase model \eqref{eq.ReducedBNmodel_nc} with a second order well balanced path-conservative scheme based on the Osher-Romberg numerical flux introduced for the first time for the shallow water equations in cylindrical coordinates in \cite{gaburro2016direct} and for the Euler equations with gravity in \cite{gaburro2018well}.

In these two works the scheme was presented in the case of complex unstructured nonconforming moving meshes; here instead, to increase the overall efficiency of the algorithm, we will work with fixed Cartesian grids.

{For the sake of completeness, and since it has never been done explicitly before, we will also briefly describe the Osher-Romberg scheme for standard system of conservation laws in Section \ref{ssec.OsherRomberg_X_conservativeSystem}.}

To discretize the computational domain, we consider a rectangular grid with $N$ elements in the $x-$direction and $M$ elements in the $y-$direction, i.e. a total number of $N \times M$ 
quadrilateral elements $T_{ij} = [x_{i-\frac{1}{2}}, x_{i+\frac{1}{2}}] \times [y_{j-\frac{1}{2}}, y_{j+\frac{1}{2}}]$, $i=1, \dots N$, $j=1\dots M$. We denote $|T_{ij}|$ the area of an element $T_{ij}$, $\Delta x$ the length of a horizontal edge, $\Delta y$ the length of a vertical edge and $\Delta t$ the size of the current time step.

The evolution of the conserved variables $\Q$ inside each element $T_{ij}$ from time $t^n$ to time $t^{n+1} = t^n + \Delta t$ is obtained, after having integrated \eqref{eq.ReducedBNmodel_nc} over $T_{ij}$, by a path-conservative scheme that reads  
\be
\label{eq.Scheme}
\Q_{ij}^{n+1} = \Q_{ij}^n 
&-\frac{\Delta t}{\Delta x} \left (\D_{i\!-\!\frac{1}{2},j}^+\left( \q_{i\!-\!\frac{1}{2},j}^{-}, \q_{i\!-\!\frac{1}{2},j}^{+}  \right) +  \D_{i\!+\!\frac{1}{2},j}^-\left( 
\q_{i\!+\!\frac{1}{2},j}^{-},  \q_{i\!+\!\frac{1}{2},j}^{+} \right) \right  ) 
-\frac{\Delta t}{\Delta y} \left (\D_{i,j\!-\!\frac{1}{2}}^+\left( \q_{i, j\!-\!\frac{1}{2}}^{-}, \q_{i, j\!-\!\frac{1}{2}}^{+}  \right) +  \D_{i, j\!+\!\frac{1}{2}}^-\left( 
\q_{i,j\!+\!\frac{1}{2}}^{-},  \q_{i,j\!+\!\frac{1}{2}}^{+} \right) \right  ) \\
& - \frac{\Delta t}{|T_{ij}|} \int_{T_{ij}}  \B(\q_{ij}^n) \cdot {\nabla} \q_{ij}^n
\ d\mathbf{x} dt \\
\ee 
where $\q_{ij}^n(\x,t)$ is the approximation of the conserved variables inside cell $T_{ij}$ at time $t^n$, computed via a well balanced reconstruction operator from the conserved variables $\Q_{i\pm 1, j\pm 1}$ in the neighborhood of $T_{ij}$, while $ \q_{i \pm\frac{1}{2}, j}^{\mp}=\q_i(x_{i\pm \frac{1}{2}}, y,t)$ and  $ \q_{i, j\pm \frac{1}{2}}^{\mp}=\q_i(x, y_{j\pm\frac{1}{2}},t)$ denote the evaluation of $\q_{ij}^n(\x,t)$ at the mid point of each edge of $T_{ij}$.

According to \cite{Pares2006} and \cite{Castro2012} $\D^{\pm}_{i+\frac{1}{2},j}$ and $\D^{\pm}_{i, j+\frac{1}{2}}$ can be defined as follows:
\be 
\label{eq.FIrstOrderDgeneric_x}
\D^{\pm}_{i\!+\!\frac{1}{2},j} & \left( \q_{i\!+\!\frac{1}{2},j}^{-}, \q_{i\!+\!\frac{1}{2},j}^{+} \right) = \frac{1}{2} \Biggl [ \ \f(\q_{i\!+\!\frac{1}{2},j}^{+}) - \f(\q_{i\!+\!\frac{1}{2},j}^{-})  + 
\Bn_{i\!+\!\frac{1}{2},j} \left ( \q_{i\!+\!\frac{1}{2},j}^{+}-\q_{i\!+\!\frac{1}{2},j}^{-} \right ) \pm \Vn_{i\!+\!\frac{1}{2},j} \left ( \q_{i\!+\!\frac{1}{2},j}^{+}-\q_{i\!+\!\frac{1}{2},j}^{-} \right ) \Biggr],
\ee 
\be 
\label{eq.FIrstOrderDgeneric_y}
\D^{\pm}_{i, j\!+\!\frac{1}{2}}& \left( \q_{i, j\!+\!\frac{1}{2}}^{-}, \q_{i, j\!+\!\frac{1}{2}}^{+} \right) = \frac{1}{2} \Biggl [ \ \g(\q_{i, j\!+\!\frac{1}{2}}^{+}) - \g(\q_{i, j\!+\!\frac{1}{2}}^{-})  + 
\Bn_{i, j\!+\!\frac{1}{2}} \left ( \q_{i, j\!+\!\frac{1}{2}}^{+}-\q_{i, j\!+\!\frac{1}{2}}^{-} \right ) \pm \Vn_{i, j\!+\!\frac{1}{2}} \left ( \q_{i, j\!+\!\frac{1}{2}}^{+}-\q_{i, j\!+\!\frac{1}{2}}^{-} \right ) \Biggr],
\ee 
where $\f$ and $\g$ are the physical fluxes, $\Bn_{i+\frac{1}{2},j}$ and $\Bn_{i, j+\frac{1}{2}}$ are the discretization of the non-conservative terms and  $\Vn_{i+\frac{1}{2},j}$ and $\Vn_{i,j+\frac{1}{2}}$ are the viscosity terms respectively in $x$- and $y$-directions. 
We would like to underline that the introduction of a path-conservative scheme for the discretization of \eqref{eq.ReducedBNmodel_nc} is motivated by two different reasons, i.e. the non-conservative advective equation for the evolution of $\alpha$, and the non-conservative terms that have been deliberately introduced to deal with equilibria \eqref{eq.BNs_eq1}, \eqref{eq.BNs_eq3} and \eqref{eq.BNs_eq4}. 
The presence of the gravity source term will lead to two different numerical fluxes in $x$- and $y$- direction, respectively. 
Indeed, in $x$-direction the only non-conservative term is 
\be 
b_4^x := u\,\partial_x{\alpha};
\ee  
which is trivially zero at the equilibrium, so that we can use the standard path-conservative Osher scheme described in \cite{OsherNC} (well balanced only for trivial equilibria).
Instead in $y$-direction, in addition to the simple term 
\be b_4^y := v \, \partial_y{\alpha},\ee
we have the crucial non-conservative part
\be b_3^y := \partial_y{\left( \alpha p \right)} + \alpha \rho g \, \partial_y{y}, \ee 
for which, to ensure that equilibria are preserved up to machine precision, we need to use a well balanced treatment, as the one offered by the Osher-Romberg scheme presented in \cite{gaburro2018well}.

In both the cases, $\Bn$ and $\Vn$ are defined in terms of a family of paths $\Phi(s, \q^{-},\q^{+})$, $s \in [0,1]$.  
In general,  according to  the theory of \cite{DalMaso1995}, the family of paths should be a Lipschitz continuous family of functions  $\Phi(s,\q^-, \q^+), \ s \in [0,1]$  satisfying some regularity and compatibility conditions, in particular,
\be 
\label{eq.Pathproperties}
\begin{aligned}
	\Phi(0,\q^-, \q^+) = \q^-, \quad \Phi(1; \q^-, \q^+) = \q^+, \quad \Phi(s, \q, \q) = \q.
\end{aligned}
\ee
Moreover, according to \cite{Pares2006}, $\D^{\pm}$  should satisfy the following properties:
\be 
\begin{aligned}
	\D^{\pm} (\q, \q) = \0 \quad \forall \q \in \Omega_{\q},
\end{aligned}
\ee
being $\Omega_{\q}$ the set of admissible states for the problem, and, for every $\q^-$, $\q^+$ $\in \Omega_{\q}$, 
\be  
\label{PC}
\begin{aligned}
	\D^-(\q^-,\q^+) + \D^+(\q^-,\q^+) = \int_0^1 \!\! \mbf{A}_k\left(\Phi\left(s;\q^-,\q^+\right)\right)\de{\Phi}{s}(s;\q^-, \q^+)ds, 	
\end{aligned}
\ee
where $\A_k$, $k=1,2$  is 
\be
\begin{aligned}
	&\mbf{A}_1 \!=\! \de{\f}{\q} + \mbf{B}_1 \!=\! \begin{pmatrix}
		0 & 1 & 0 & 0 & 0 \\ 
		c^2-u^2 & 2u & 0 & p-\rho c^2 & 0 \\ 
		-uv & v & u & 0 & 0 \\ 
		0 & 0 & 0 & u & 0 \\ 
		0 & 0 & 0 & 0 & 0
	\end{pmatrix}\!, \ \   
	\text{and} \ \ 
	\mbf{A}_2 \!=\! \de{\g}{\q} + \mbf{B}_2 \!=\! \begin{pmatrix}
		0 & 0 & 1 & 0 & 0 \\ 
		-uv & v & u & 0 & 0 \\ 
		c^2-v^2 & 0 & 2v & p-\rho c^2  & \alpha \rho g \\ 
		0 & 0 & 0 & v & 0 \\ 
		0 & 0 & 0 & 0 & 0
	\end{pmatrix}\!.
\end{aligned}
\ee 
with the speed of sound $c$ defined as
\be
c^2 = \de{p}{\rho} = \frac{\gamma k_0}{\rho_0}\left( \frac{\rho}{\rho_0} \right)^{\gamma-1}.
\ee 
The interested reader is referred to \cite{DalMaso1995} and \cite{Pares2006} for a rigorous and complete presentation of this theory.
{ 
The right eigenvectors and the eigenvalues of the matrices $\mathbf{A}_1$ and $\mathbf{A}_2$ of the augmented system including the gravity potential read  
\begin{equation}
  \mathbf{R}_1 = \left(  \begin{array}{ccccc} 
     1 & 0 & \rho c^2 - p & 0 & 1 \\ 
     u-c & 0 & u(\rho c^2 - p) & 0 & u+c \\ 
     v & 1 & 0 & 0 & v \\ 
     0 & 0 & c^2 & 0 & 0 \\ 
     0 & 0 & 0 & 1 & 0 
  \end{array} 
  \right), \qquad \Lambda_1 = \textnormal{diag}\left( u-c, u, u, 0, u+c \right),  
\end{equation}
and 
\begin{equation}
\mathbf{R}_2 = \left(  \begin{array}{ccccc} 
1 & 0 & \rho c^2 - p & \alpha \rho g & 1 \\ 
u & 1 & 0 & u \alpha \rho g & u \\ 
v-c & 0 & v(\rho c^2 - p) & 0 & v+c \\ 
0 & 0 & c^2 & 0 & 0 \\ 
0 & 0 & 0 & v^2-c^2 & 0 
\end{array} 
\right), \qquad \Lambda_2 = \textnormal{diag}\left( v-c, v, v, 0, v+c \right).  
\end{equation}
The inverse matrices of $\mathbf{R}_1$ and $\mathbf{R}_2$ are given by 
\begin{equation}  
    \mathbf{R}_1^{-1} = \left(  \begin{array}{ccccc} 
    \frac{1}{2} \frac{c+u}{c} & -\frac{1}{2c} & 0 & -\frac{1}{2}\frac{\rho c^2-p}{c^2} & 0 \\[2pt]  
    -v & 0 & 1 & v \frac{\rho c^2-p}{c^2} & 0 \\[2pt]
    0 & 0 & 0 & \frac{1}{c^2} & 0 \\[2pt] 
    0 & 0 & 0 & 0 & 1 \\[2pt] 
  \frac{1}{2} \frac{c-u}{c} & +\frac{1}{2c} & 0 & -\frac{1}{2}\frac{\rho c^2-p}{c^2} & 0 \\[2pt]      
    \end{array} 
    \right), \qquad \textnormal{and} \qquad 
     \mathbf{R}_2^{-1} = \left(  \begin{array}{ccccc} 
  \frac{1}{2} \frac{c+v}{c} & 0 & -\frac{1}{2c} & -\frac{1}{2}\frac{\rho c^2-p}{c^2} & \frac{1}{2c} \frac{\alpha \rho g }{c-v} \\[2pt]  
  -u & 1 & 0 & u \frac{\rho c^2-p}{c^2} & 0 \\[2pt]
  0 & 0 & 0 & \frac{1}{c^2} & 0 \\[2pt] 
  0 & 0 & 0 & 0 & \frac{1}{v^2-c^2} \\[2pt] 
  \frac{1}{2} \frac{c-v}{c} & 0 & +\frac{1}{2c} & -\frac{1}{2}\frac{\rho c^2-p}{c^2} & 
  \frac{1}{2c} \frac{\alpha \rho g }{c+v}  \\[2pt]      
  \end{array} 
  \right).  
\end{equation} 
The full eigenstructure of the PDE system and will be needed later for the construction of the  numerical scheme.  
}
The rest of this section is organized as follows.  First, 
{we recall the Osher scheme and we subsequently introduce the Osher-Romberg scheme for conservative systems.} Then, we give the details on the first order path-conservative scheme used in the $x$-direction (the standard Osher scheme) and in the $y$-direction (the well-balanced Osher-Romberg scheme). Next, we propose a second order scheme constructed using the previous first order schemes in combination with a second order well balanced reconstruction operator.
We close the section by explaining how also an equilibrium solution that may \textit{not} be known \textit{a priori} can 
be preserved, by automatically deducing it from the cell averages at the beginning of each time step. 

\subsection{{Osher and Osher-Romberg schemes for conservative systems}}
\label{ssec.OsherRomberg_X_conservativeSystem}

{Let us consider a system of conservation law given by \eqref{eq.generalform} with $\B=\S=\0$, for which the numerical flux at the cell interface can be written as 
\begin{equation}
\f\left( \q^{-}, \q^{+} \right) = 
\frac{1}{2} \left( \ \f(\q^{+}) + \f(\q^{-})  \ \right)
- \frac{1}{2}  \Vn \left ( \q^{+}-\q^{-}  \right ).  
\end{equation}
It contains two standard parts:
the first one given by the arithmetic average of the fluxes computed at the left and right of each interface, 
and a second one $\Vn(\q^+ - \q^-)$, called \textit{numerical viscosity}, that stabilizes the method and that depends on the eigenstructure of the Jacobian matrix and on the jump of the conserved variables at the interface. 
Depending on the way in which the term $\Vn(\q^+ - \q^-)$ is computed one obtains different numerical schemes.		
}

{
The \textit{Osher--Solomon} flux goes back to \cite{osherandsolomon} and has been recently extended to general conservative and non--conservative systems within the framework of path--conservative schemes   \cite{OsherUniversal, OsherNC,ApproxOsher}. Its numerical viscosity reads    
\be
\label{eq.OsherMatrixFormal_2b}
\Vn \left ( \q^{+}-\q^{-}  \right )= \int_0^1  \!\left | \mbf{A} \left ( \Phi(s)  \right) \cdot \mathbf{n}    \right | \partial_s \Phi(s) ds, \ \   0 \le s \le 1,
\ee 
where $\A = [\A_1, \A_2]$ is the Jacobian of the system and $\mathbf{n}$ the outward pointing normal at the interface between two states $\q^+$ and $\q^-$.
The absolute value of $\A$ is evaluated as usual as 
\begin{equation}
|\mathbf{A}| = \mathbf{R} |\boldsymbol{\Lambda}| \mathbf{R}^{-1},  \qquad |\boldsymbol{\Lambda}| = \textnormal{diag}\left( |\lambda_1|, |\lambda_2|, ..., |\lambda_\nu| \right),  
\end{equation}
where $\mathbf{R}$ and $\mathbf{R}^{-1}$ denote, respectively, the right eigenvector matrix of $\A \cdot \mathbf{n}$ and its inverse.
}

{
Finally, in order to compute the integral in \eqref{eq.OsherMatrixFormal_2b} whose primitive is in general not available, it has been proposed in \cite{ApproxOsher} to connect $\q^+$ and $\q^-$ using a simple straight--line segment path  
\begin{equation}
\boldsymbol{\Psi}(s) = \q^- + s \left( \q^+ - \q^- \right), \qquad 0 \leq s \leq 1.
\label{eqn.path} 
\end{equation} 
By choosing a quadrature formula with nodes $s_k$ and weights $\omega_k$ one obtains
\be
\label{eq.Osher_discreto}
\Vn \left ( \q^{+}-\q^{-}  \right )= \sum_{k=1}^{\ell} \omega_k \left | \mbf{A} \left ( \Phi(s_k)  \right) \cdot \mathbf{n}   \right |  \left(\q^+ - \q^- \right).
\ee 
}

{
Next, with the scope of \textit{modifying} \eqref{eq.Osher_discreto} in such a way that it does not depend on the jump between the conserved variables $\left(\q^+ -\q^-\right)$, one notices that 
\eqref{eq.OsherMatrixFormal_2b} can be written as
\be
\Vn(\q^+ - \q^-) = \! \int_0^1  \!\!\text{sign} \left( \mbf{A}  \left ( \Phi(s)  \right) \cdot \mathbf{n}  \right )\mbf{A}\left ( \Phi(s)  \right) \cdot \mathbf{n} \, \partial_s \Phi(s) ds,  
\ee 
with $\text{sign}( \A ) = \mathbf{R} \, \text{sign} ( \boldsymbol{\Lambda} ) \mathbf{R}^{-1}$ and $\text{sign} ( \boldsymbol{\Lambda} )$ the diagonal matrix containing the signs of all eigenvalues of 
$\A \cdot \mathbf{n}$.
Now, by exploiting the fact that $\A(\Q) = \partial \F / \partial \Q$ and by using a quadrature formula we can approximate the previous expression with 
\be 
\label{eq.OR_cons}
\Vn(\q^+-\q^-) = \sum_{k=1}^l \omega_k \text{sign} \left(\mbf{A}(\Phi(s_k) \cdot \mathbf{n} \right) 
\left( \frac{ \F(\Phi(s_k+\epsilon_k)) - \F(\Phi(s_k-\epsilon_k)) }{2\epsilon} \right),
\ee
where, in particular, the use of the Romberg method with $l=3$ and
\be 
\label{eq.RombergWeights}
s_1=1/4, \,s_2=3/4, \,s_3=1/2, \quad 
\omega_1=2/3, \,\omega_2=2/3, \,\omega_3=-1/3, \quad 
\epsilon_1=1/4, \,\epsilon_2=1/4, \,\epsilon_3=1/2,
\ee
reduces the computational effort allowing to re-use values of $\F$ already computed in the first part of the flux.
}

{
The obtained result \eqref{eq.OR_cons} has been named Osher-Romberg scheme in \cite{gaburro2016direct} and \cite{gaburro2018well}, because it employs the Romberg method as quadrature formula  and it has been derived by slightly modifying the Osher--type scheme of \cite{OsherNC,OsherUniversal}. Since also the Osher-Romberg method makes use of the full eigenstructure of the PDE system, it is a complete Riemann solver that is little dissipative. 
We emphasize that the principal feature of this Riemann solver is that it does \textit{not} explicitly depend on the jump of the conserved variables at the interface any more, which has been substituted by the difference of the fluxes. This fact will play a fundamental role when also non-conservative terms are present and greatly simplifies the derivation of a well balanced scheme.
}

\subsection{First order scheme in $x$-direction}

Let us first remark that the scheme \eqref{eq.Scheme} across a vertical edge reduces to 
\be 
\label{eq.SemidiscreteScheme_1}
\begin{aligned} 
	\Q_{ij}^{n+1} \!= \Q_{ij}^n 
	&\!-\!\frac{\Delta t}{\Delta x} \left (\!\D_{i\!-\!\frac{1}{2},j}^+\left( \q_{i\!-\!\frac{1}{2},j}^{-}, \q_{i\!-\!\frac{1}{2},j}^{+}  \! \right) \!+  \D_{i\!+\!\frac{1}{2},j}^-\left( 
	\q_{i\!+\!\frac{1}{2},j}^{-},  \q_{i\!+\!\frac{1}{2},j}^{+}\!  \right) \! \right  )\!, 
\end{aligned} 
\ee 
where $\D^{\pm}_{i\!+\!\frac{1}{2},j}$ is given by \eqref{eq.FIrstOrderDgeneric_x}, 
if $\q_{ij}^n(\x,t)$ is constant within each cell for every time $t \in [t^n, t^{n+1}]$ and coincides with the cell average $\Q_{ij}^n$. Thus, the resulting  scheme is first order accurate in space and time. Moreover,  $\q_{i+\frac{1}{2},j}^{-}=\q_{ij}=\Q_{ij}$ and $\q_{i+\frac{1}{2},j}^{+}=\q_{i+1,j}=\Q_{i+1,j}$.

For what concerns the discretization of the non-conservative terms $ \Bn_{i\!+\!\frac{1}{2},j} \left ( \q_{i\!+\!\frac{1}{2},j}^{+}-\q_{i\!+\!\frac{1}{2},j}^{-} \right )$, there is only one non zero component that can be easily written as 
\be
\label{eq.b_4^x}
\left(b_{i\!+\!\frac{1}{2},j}\right)_4^{x}  = \overline{u} \Delta \alpha,
\ee 
where $\Delta \alpha = \alpha_{i+1,j} - \alpha_{ij}$, and we propose to use the Roe average velocity to compute $\overline{u}$ 
\be
\overline{u} = \frac{u_{ij}\sqrt{\rho_{ij}} + u_{i+1,j}\sqrt{\rho_{i+1,j}}}{\sqrt{\rho_{ij}}+\sqrt{\rho_{i+1,j}}}.
\ee 
The above expression is so simple because of the equilibrium constraint \eqref{eq.BNs_eq1} that ensures $b^{x}_4=0$ when $\q_{ij}^n$ and  $\q_{i+1,j}^n$ lie on the same stationary solution.

Then, for what concerns the viscosity term, {adapting \eqref{eq.OsherMatrixFormal_2b} to the context}, we can write
\be
\label{eq.OsherMatrixFormal}
\Vn_{i\!+\!\frac{1}{2},j} \left ( \q_{i\!+\!\frac{1}{2},j}^{+}-\q_{i\!+\!\frac{1}{2},j}^{-} \right )=\int_0^1  \!\left | \mbf{A}_1 \left ( \Phi(s)  \right)  \right | \partial_s \Phi(s) ds, \ \   0 \le s \le 1.
\ee 
In particular, to approximate the integral term in \eqref{eq.OsherMatrixFormal} we use a Gauss-Legendre quadrature rule with $k=3$ quadrature points $s_k = (\frac{1}{2}-\frac{\sqrt{15}}{10},\frac{1}{2},\frac{1}{2}+\frac{\sqrt{15}}{10} ) $ in the unit interval $[0,1]$ and associated weights $\omega_k = ( \frac{5}{18}, \frac{8}{18}, \frac{5}{18} ) $
\be
\label{eq.OsherMatrixFormal_approx}
\Vn_{i\!+\!\frac{1}{2},j} \left ( \q_{i\!+\!\frac{1}{2},j}^{+}-\q_{i\!+\!\frac{1}{2},j}^{-} \right ) = \left( \sum_{k=1}^3   \omega_k \left | \mbf{A}_1  \left ( \Phi(s_k)  \right)  \right |  \right)  \left ( \q_{i\!+\!\frac{1}{2},j}^{+}-\q_{i\!+\!\frac{1}{2},j}^{-} \right ) \! .
\ee

As pointed out in \cite{Pares2006}, a sufficient condition for a first order path-conservative scheme to be well balanced is that $\D^{\pm}_{i+\frac{1}{2},j} \left(\q_{ij}^n, \q_{i+1,j}^n \right) =\0$, if $\q_{ij}^n$ and $\q_{i+1,j}^n$ lie on the same stationary solution. 
Therefore in the $x$-direction the condition is trivially fulfilled because at the equilibrium  $u=0$, see \eqref{eq.BNs_eq1}, \eqref{eq.BNs_eq3} and \eqref{eq.BNs_eq4}, so the term \eqref{eq.b_4^x} vanishes, and $\q_{ij}^n = \q_{i+1,j}^n$, hence also \eqref{eq.OsherMatrixFormal_approx} is null.

\subsection{First order scheme in $y$-direction}

Across a horizontal edge a first order path-conservative scheme reads as follows  
\be 
\label{eq.SemidiscreteScheme_2}
\begin{aligned} 
	\Q_{ij}^{n+1} = \Q_{ij}^n 
	&\!-\!\frac{\Delta t}{\Delta y} \left (\D_{i,j\!-\!\frac{1}{2}}^+\left( \q_{i, j\!-\!\frac{1}{2}}^{-}, \q_{i, j\!-\!\frac{1}{2}}^{+}   \right) +  \D_{i, j\!+\!\frac{1}{2}}^-\left( 
	\q_{i,j\!+\!\frac{1}{2}}^{-},  \q_{i,j\!+\!\frac{1}{2}}^{+}  \right) \right  )\!,  \\
\end{aligned} 
\ee 
where in particular $\D^{\pm}_{i, j\!+\!\frac{1}{2}}$ is given by \eqref{eq.FIrstOrderDgeneric_y} and $\q_{i, j+\frac{1}{2}}^{-}\!=\q_{ij}=\Q_{ij}$ and $\q_{i, j+\frac{1}{2}}^{+}=\q_{i,j+1}=\Q_{i,j+1}$.

The discretization of the non-conservative term $b^y_4$ can be given again in an easy way by choosing 
\be
\label{eq.b_4^y}
\left(b_{i, j\!+\!\frac{1}{2}}\right)^{y}_4 = \overline{v} \Delta \alpha, 
\ee 
where $\Delta \alpha = \alpha_{i,j+1} - \alpha_{ij}$ and 
\be
\overline{v} = \frac{v_{ij}\sqrt{\rho_{ij}} + 
	v_{i,j+1}\sqrt{\rho_{i,j+1}}}{\sqrt{\rho_{ij}}+\sqrt{\rho_{i,j+1}}},
\ee 
which is zero when $\q_{ij}^n$ and $\q_{i,j+1}^n$ lie on the same stationary solution, thanks to the equilibrium constraint \eqref{eq.BNs_eq1}.
We remark that this expression results from choosing a segment path 
\be
\label{eq.SegmentPath}
\Phi(s; \q_{ij}, \q_{i, j+1}) = \q_{ij} + s(\q_{i, j+1} - \q_{ij}),
\ee
to connect the left states $v_{ij}$ and $\alpha_{ij}$ to the right ones $v_{i,j+1}$ and $\alpha_{i,j+1}$.

If the same choice had been made for density and pressure in order to discretize $b_{3}^y$ the final scheme would not have been well balanced (as it happened in \cite{dumbser2011simple}).
So, in this work, we propose the following family of paths so that stationary solutions given by \eqref{eq.BNs_eq1}, \eqref{eq.BNs_eq3} and \eqref{eq.BNs_eq4} are preserved.
Let $\Phi^E(s,\Q^E_{ij}, \Q^E_{i, j+1})$ be a reparametrization of a stationary solution given by \eqref{eq.BNs_eq1}, \eqref{eq.BNs_eq3} and \eqref{eq.BNs_eq4} that connects the state $\Q^E_{ij}$ with $\Q^E_{i, j+1}$, where $\Q^E_{ij}$ is the cell average of the given stationary solution in the cell $T_{ij}$. Note that in the case of first and second order schemes $\Q^E_{ij}$ could be approximated by the evaluation of the stationary solution at the center of the cell.  Then we define $\Phi(s,\q_{ij}, \q_{i,j+1})$ as follows
\be
\label{eq:path}
\Phi(s,\q_{ij},\q_{i,j+1})=\Phi^E(s,\Q^E_{ij}, \Q^E_{i,j+1})+\Phi^f(s,\q^f_{ij}, \q^f_{i, j+1}), 
\ee
where $\q^f_{ij}=\q_{ij}- \Q^E_{ij}$ and $\q^f_{i,j+1}=\q_{i,j+1}- \Q^E_{i, j+1}$ and
\be 
\Phi^f(s,\q^f_{ij}, \q^f_{i,j+1})=  \q^f_{ij} + s(\q^f_{i,j+1} - \q^f_{ij}).
\ee
That is, $\Phi^f$ is a segment path on the {\it fluctuations} with respect to a given stationary solution. (Note that trivially the segment path used to obtain \eqref{eq.b_4^y} can be cast in form \eqref{eq:path} simply by choosing $\Phi^E(s, v_{ij}/\alpha_{ij}, v_{i, j+1}/\alpha_{i,j+1})=0$).
For the sake of simplicity, in the following we will use the notation $\Phi(s)$ instead of $\Phi(s; \q_{ij}, \q_{i, j+1})$ when there is no confusion.

With this choice, it is clear that if $\q_{ij}$ and $\q_{i,j+1}$ lie on the same stationary solution satisfying \eqref{eq.BNs_eq1}, \eqref{eq.BNs_eq3} and \eqref{eq.BNs_eq4}, then $\q^f_{ij}=\q^f_{i,j+1}=\0$ and $\Phi$ reduces to $\Phi^E$. 
In such situations we have that $\g(\q_{i,j+1})=\g(\q_i)=\0$ and
\be
\Bn_{i, j+\frac{1}{2}} & (\q_{i, j+1}-\q_{ij})= \int_0^1 \!\!\! \B_2(\Phi^E(s)) \de{\Phi^E}{s}(s)ds =\0.
\ee
Therefore
\be 
\g(\q_{i,j+1})-\g(\q_{ij})+\Bn_{i, j+\frac{1}{2}}(\q_{i, j+1}-\q_{ij})=\0. 
\ee

Let us now define $\left(b_{i, j\!+\!\frac{1}{2}}\right)^y_3$ in the general case, where $\q_{i,j+1}$ and $\q_{ij}$ do not lie on a stationary solution. In this case we have that

\be
\left(b_{i, j\!+\!\frac{1}{2}}\right)^y_3 = \int_0^1 \de{\Phi_\alpha(s)}{s}\de{\Phi_p(s)}{s} + \Phi_{\alpha}(s)\Phi_{\rho}(s)g  \de{\Phi_y(s)}{s}ds,
\ee
where $\Phi_y(s)=\Phi_y(s; y_{ij},y_{i,j+1})=y_{ij} + s(y_{i,j+1}-y_{ij})$, $\Phi_{\alpha, \rho, p}(s)=\Phi^E_{\alpha, \rho, p}(s)+\Phi^f_{\alpha, \rho, p}(s)$ respectively.
Taking into account that 
\be
\int_0^1  \de{\Phi_\alpha^E(s)}{s}\de{\Phi_p^E(s)}{s} + \Phi_{\alpha}^E(s)\Phi_{\rho}^E(s)g  \de{\Phi_y(s)}{s}ds = 0
\ee
$\left(b_{i, j\!+\!\frac{1}{2}}\right)^y_3$ can be rewritten as follows
\be
& \!\! \int_{0}^{1} \! \left (  
\de{\Phi_\alpha^E\!\left(s\right)}{s}\de{\Phi_p^f\!\left(s\right)\!}{s} \!+\! \de{\Phi_\alpha^f\!\left(s\right)}{s}\de{\Phi_p^E\!\left(s\right)\!}{s} \!+\! \de{\Phi_\alpha^f\!\left(s\right)}{s}\de{\Phi_p^f\!\left(s\right)\!}{s} \right)ds 
  \! + \!\! \int_0^1  \!\! \left( \Phi_{\alpha}^E(s)\Phi_{\rho}^f(s) + \Phi_{\alpha}^f(s)\Phi_{\rho}^E(s) + \Phi_{\alpha}^f(s)\Phi_{\rho}^f(s) \right) \!g  \de{\Phi_y(s)}{s}ds.  \\ 
\ee
In general, the integral term could be difficult to compute, therefore we propose to use a numerical quadrature formula. Here the mid-point rule is used. In this case, we define $\left(b_{i, j\!+\!\frac{1}{2}}\right)^y_3$ as follows
\be
\label{eq.b^y_3_final}
\Delta \alpha^E\!\Delta p^f \!\!+ \Delta \alpha^f\!\Delta p^E \!\!+ \Delta \alpha^f\!\Delta p^f  \!\!+ \left(\!\overline{\alpha^E}\overline{\rho^f} \!\!+ \overline{\alpha^f}\overline{\rho^E} \!+ \overline{\alpha^f}\overline{\rho^f} \right) \!g \Delta y, 
\ee
where
\be
&\Delta p^{E,f}\!\!= p^{E,f}_{i,j+1}\!-\!p^{E,f}_{ij}\!, \quad  \Delta \alpha^{E,f}\!\!= \alpha^{E,f}_{i,j+1}\!-\!\alpha^{E,f}_{ij}\!,  
\quad  \Delta y\!= y_{i,j+1}\!-\!y_{ij}, \\
& \overline{\alpha^{E,f}}= \frac{\alpha^{E,f}_{ij} + \alpha^{E,f}_{i,j+1}}{2}, \ \text{ and } \quad  \overline{\rho^{E,f}}= \frac{\rho^{E,f}_{ij} + \rho^{E,f}_{i,j+1}}{2}.
\ee
It is clear from the definition that $\left(b_{i, j\!+\!\frac{1}{2}}\right)^y_3=0$ if $\q_{ij}^n$ and $\q_{i,j+1}^n$ lie on the same stationary solution because all the fluctuation terms are zero.

\bigskip

\noindent \textbf{{Well balanced} Osher-Romberg viscosity.} 

\noindent Finally, for the viscosity term $\Vn_{i, j\!+\!\frac{1}{2}} \left ( \q_{i, j\!+\!\frac{1}{2}}^{+}-\q_{i, j\!+\!\frac{1}{2}}^{-} \right )$ we use the well balanced Osher-Romberg scheme described in detail in \cite{gaburro2018well}, which guarantees automatically a well balanced viscosity term $\Vn$ whenever $\Bn$ is discretized in a well balanced way (even in the case of complex equilibria).

We re-propose here the derivation of the Osher-Romberg scheme starting  
from the standard Osher scheme 
\be
\label{eq.OsherMatrixFormal_2}
\Vn_{i, j\!+\!\frac{1}{2}} \left ( \q_{i,j\!+\!\frac{1}{2}}^{+}-\q_{i, j\!+\!\frac{1}{2}}^{-}  \right )= \int_0^1  \!\left | \mbf{A}_2  \left ( \Phi(s)  \right)  \right | \partial_s \Phi(s) ds, \ \   0 \le s \le 1,
\ee 
{in the case $\B \ne \0$, i.e. when the presence of non-conservative terms requires the use of a well balanced path-conservative scheme.}
First we notice that the previous expression can be written as
\be
\Vn_{i,j\!+\!\frac{1}{2}}(\q_{i,j+1} - \q_{ij}) = \! \int_0^1  \!\!\text{sign} \left( \mbf{A}_2 \left ( \Phi(s)  \right) \right )\mbf{A}_2\left ( \Phi(s)  \right)\partial_s \Phi(s) ds.  
\ee 
Then, we approximate the previous expression by a quadrature formula as follows
\be 
\Vn_{i, j\!+\!\frac{1}{2}} (\q_{i, j+1}\!-\!\q_{ij}) = \sum_{k=1}^l \omega_k \text{sign} \left(\mbf{A}_2(\Phi(s_k)\right) \mbf{A}_2(\Phi(s_k))\partial_s \Phi(s_k).
\ee
Now, we approximate $\mbf{A}_2(\Phi(s_k))\partial_s \Phi(s_k)$ by the following expression
\be 
\mbf{A}_2(\Phi(s_k)) \partial_s \Phi(s_k) \approx \frac{\mbf{A}_{\Phi_k}}{2\epsilon_j}\left(\Phi(s_k+\epsilon_k)-\Phi(s_k-\epsilon_k)\right),
\ee
where $\mbf{A}_{\Phi_k}=\A_2(\Phi(s_k-\epsilon_k),\Phi(s_k+\epsilon_k))$ is a Roe-matrix associated to the system (see  \cite{Pares2006} for details), i.e. is a matrix satisfying
\be 
\mbf{A}_{\Phi_k}&\left(\Phi(s_k+\epsilon_k)- \Phi(s_k-\epsilon_k)\right)= \g(\Phi(s_k+\epsilon_k))-\g(\Phi(s_k-\epsilon_k)) 
 +\Bn_{\Phi_k}\left(\Phi(s_k+\epsilon_k)-\Phi(s_k-\epsilon_k)\right),  
\ee 
where $\Bn_{\Phi_k}\left(\Phi(s_k+\epsilon_k)-\Phi(s_k-\epsilon_k)\right)$ is defined as in \eqref{eq.b_4^y}, \eqref{eq.b^y_3_final} using the states $\Phi(s_k-\epsilon)$ and $\Phi(s_k+\epsilon)$.
Therefore, the viscosity term reads as follows
\be
\label{viscosity_n_osher}
\Vn_{i, j\!+\!\frac{1}{2}}(\q_{i, j+1}-\q_{ij}) =\sum_{k=1}^l \omega_k \text{sign}\left( \mbf{A}_2(\Phi(s_k)\right) \frac{\mathcal{R}_k}{2\epsilon_k},
\ee
where 
\be
& \mathcal{R}_k=\g(\Phi(s_k+\epsilon_k))-\g(\Phi(s_k-\epsilon_k))+ \Bn_{\Phi_k}\left(\Phi(s_k+\epsilon_k)-\Phi(s_k-\epsilon_k)\right)\!. 
\ee

Note that if $\q_{ij}^n$ and $\q_{i,j+1}^n$ lie on the same stationary solution we have $\Phi(s)=\Phi^E(s)$ and $\mathcal{R}_k=\0$, $k=1, \dots, l$ and $\Vn$ vanishes. Therefore, the numerical scheme \eqref{eq.SemidiscreteScheme_2} with \eqref{eq.FIrstOrderDgeneric_y}, where $\Bn$ is defined with \eqref{eq.b_4^y}, \eqref{eq.b^y_3_final} and $\Vn$ being defined by \eqref{viscosity_n_osher} is exactly well balanced for stationary solutions given by \eqref{eq.BNs_eq1}, \eqref{eq.BNs_eq3} and \eqref{eq.BNs_eq4}.

Then, following \cite{gaburro2018well}, we use the Romberg method with $l=3$ (see \eqref{eq.RombergWeights}). 
With this choice, the  viscosity term of the Osher-Romberg method reads as follows:
\be
\label{eq.ViscosityOR}
\Vn_{i, j+\frac{1}{2}}(\q_{i, j+1}-\q_{ij}) = 
\ &\frac{4}{3} \text{sign}\left(\mbf{A}_2(\Phi(\frac{1}{4}))\right) \left( \g(\Phi(\frac{1}{2}))-\g\left(\q_{ij}\right)+\Bn_{i,j+\frac{1}{4}}\left(\Phi(\frac{1}{2})-\q_{ij}\right)\right) + \\
&\frac{4}{3}\text{sign}\left(\mbf{A}_2(\Phi(\frac{3}{4}))\right)\left(\g\left(\q_{i, j\!+1}\right)-\g(\Phi(\frac{1}{2}))+\Bn_{i, j+\frac{3}{4}}\left( \q_{i,j\!+1}-\Phi(\frac{1}{2})\right)\right) \\
-&\frac{1}{3}\text{sign}\left(\mbf{A}_2(\Phi(\frac{1}{2}))\right)\left(\g\left(\q_{i, j\!+1}\right)-\g\left(\q_{ij}\right)+\Bn_{i, j+\frac{1}{2}} \left(\q_{i,j\!+1}-\q_{ij}\right)\right)\!.
\ee

\subsection{2nd order well balanced reconstruction}
\label{ssec.2d}

According to \cite{Pares2006} and \cite{Castro2006}, scheme \eqref{eq.Scheme} is well balanced if both, the underlying first order scheme and the reconstruction operator are well balanced, and all the integrals that appear in \eqref{eq.Scheme} are computed exactly.   
Therefore, in order to define a second order scheme, a second order well balanced reconstruction operator should be defined. 	 

The most popular way to define a second order reconstruction  operator  is based on the MUSCL method introduced by van Leer in \cite{leer5} that we use joint with the Barth \& Jespersen limiter \cite{BarthJespersen}.
It is clear that the standard MUSCL method is only well balanced for linear stationary solutions, which is not the case here. In this paper we therefore follow the strategy proposed in \cite{Castro2008},  where the reconstruction operator is defined as a combination of a smooth stationary solution together with a standard reconstruction operator to reconstruct the {\it fluctuations} with respect to the  given stationary solution, that is
\be
\label{eq.FinalReconstruction2}
\q_{ij}^n(\x,t) = \Q^E_{ij}(\x,t) + \mathcal{P}^{f,n}_{ij}(\x,t), \quad (\x,t) \in T_{ij}^n \times [t^n, t^{n+1}],
\ee 
where, $\mathcal{P}^{f,n}_{ij}(\x,t)$ is the standard MUSCL method applied in order to reconstruct the fluctuations with respect to the given stationary solution computed for the neighbors $T_{i\pm 1,j}^n$,  $T_{i,j\pm 1}^n$,  $T_{i\pm 1,j\pm 1}^n$ of $T_{ij}^n$ as 
\be
\label{eq.FluctuationForRec2}
\Q^{f\!,n}_{i\pm 1,j\pm 1} \!= \Q_{i\pm 1,j\pm 1}^{n} \!- \Q^{E\!,n}_{i\pm 1,j\pm 1}.
\ee

Then, the expression of the reconstruction operator for the fluctuations is
\be
\label{reconstruction_fluctuation2}
\mathcal{P}^{f,n}_{ij} (\x,t) = \Q_{ij}^{f,n} + \Psi_{ij} \nabla \Q_{ij}^{f,n} (\x-\x_{ij}) + \partial_t\Q_{ij}^n (t -t^n),
\ee 
where $\x_{ij}$ is the barycenter of cell $T_{ij}^n$.

To compute $\nabla \Q_i^{f\!,n}$ we use the standard MUSCL method (see \cite{leer5})  together with the Barth and Jespersen limiter to recover $\Psi_{ij}$ (see \cite{BarthJespersen}). We would like to remark that the employed methods are standard, the novelty is in the fact that both are applied only to the fluctuations in order to obtain an exactly well-balanced scheme.


Next, the term $\partial_t \Q_{ij}^n$ indicates the time derivative of $\Q$ and  it can be computed using a discrete version of the governing equation
\be
\partial_t \Q_{ij}^n = \left( \mbf{J}_{\f} + \B_1 \right )\!|_{\x_{ij}} \partial_x \Q (\x_{ij}) + \left( \mbf{J}_\g  + \B_2 \right)|_{\x_{ij}} \, \partial_y \Q(\x_{ij}),
\ee
evaluated at the barycenter $\x_{ij}$ of $T_{ij}$. 
In particular the gradient of the conserved variables must be expressed as the gradient of the equilibrium plus the previously computed gradient of the fluctuation, i.e. 
\be
\label{eq.gradientEq+fluct}
\nabla \Q = \nabla \Q^E_{ij} + \Psi_{ij}  \nabla \Q_{ij}^{f,n}, 
\ee 
in order to preserve the well balancing.

Finally, we have to compute the integral term of \eqref{eq.Scheme} that can be expanded as follows
\be
\int_{T_{ij}}\!{\B}\left(\q_{ij}^n\right)\cdot{\nabla} \q_{ij}^n\, d\mathbf{x} dt = 
\int_{T_{ij}}\!\left(0, 0, \de{\alpha p}{y}+\alpha \rho g \de{y}{y}, u\de{\alpha}{x}+v\de{\alpha}{y}, 0\right)\,  d\mathbf{x} dt.
\ee
For the fourth component one can simply use the mid point quadrature rule without affecting the well balancing (because $\u=\0$  at the equilibrium). Instead for the third component a more careful approach should be employed.
Formally, one has to write $\q_{ij}^n(\x,t)$ as  $\Q^E_{ij}(\x) + \mathcal{P}^{f,n}_{ij}(\x,t)$ and using the fact that 
\be 
\int_{T_{ij}} \B( \Q^E_{ij} )\cdot{\nabla} \Q^E_{ij} \, d\mathbf{x} dt =\0,
\ee
the only remaining terms will depend only on fluctuations and so can be approximated using any quadrature formula, in particular we use the the mid point rule which is second order accurate and we obtain
\be
\label{eq.b^y_3_final_integral }
|T_{ij}| \, \Biggl( \Delta \alpha^E\Delta p^f & + \Delta \alpha^f \Delta p^E + \Delta \alpha^f\Delta p^f  
 + \left({\alpha^E_{ij}}{\rho^f_{ij}} +{\alpha^f_{ij}}{\rho^E_{ij}} + {\alpha^f_{ij}}{\rho^f_{ij}} \right) g \Delta y \Biggr ), 
\ee
where the symbol $\Delta q^E$ refers to the exact gradient in the $y$-direction of the equilibrium profile and the symbol $\Delta q^f$ refers to the reconstructed gradient always in the $y$-direction of the fluctuations. 

To conclude, we point out that in dambreak problems in order to avoid wasting computational time we do not use the second order version of the code when approximating cells with vanishing liquid volume fraction (cells occupied by the surrounding zero pressure gas).
To be more specific, when simulating water-air problems we  always choose the initial conditions for the volume fraction $\alpha$ as
\be
\alpha(\x,0) = \begin{cases} 1-\epsilon, \quad &\text{if } \  \x \in \Omega_l \\\epsilon, \quad &\text{if } \ \x \in \Omega_a \end{cases} 
\ee
where $\epsilon > 0$ is a small parameter.
Then, we use the second order scheme  only  wherever $\alpha >10 \epsilon$, otherwise we switch to the first order method.

\subsection{Automatic detection of the equilibrium profile}

In the previous section we have presented a novel well balanced algorithm that is able to preserve up to machine precision equilibrium profiles of the form \eqref{eq.BNs_eq1}, \eqref{eq.BNs_eq3} and \eqref{eq.BNs_eq4}, where $k_0$ and $\rho_0$ are \textit{a priori} known constants that characterize the problem, and $y_0$ indicates the free-surface position.
If also $y_0$ is fixed at the beginning of the simulation, then the application of the scheme is straightforward. 

However, in general $y_0$ is locally constant for water at rest solutions, but its value is not known \textit{a priori}.  
Moreover, it can even be a function of $x$, i.e. $y_0(x)$, for example in the presence of obstacles that divide water 
at rest at different levels from each other. 
In these cases, we need a procedure able to recover an approximate value of $y_0(x)$ for any element $T_{ij}^n$; in particular, for a vertical strip $i$ of elements $T_{ij}^n$, $i$ fixed and $j=1,\dots, M$ covering a domain $[x_L, x_R]\times[y_L, y_R]$, we want to fix an unique value $y_0(x_{ij}).$ 
So at the beginning of each time step we define
\be
\label{eq.Free_surf_detector}
y_0(x_{i},t^n) = \int_{y_L}^{y_R} \!\! \alpha\left(x_{i}, y\right)\, dy \simeq y_L + \sum_{j=1}^{M} \alpha_{ij}^n, \quad \text{for any} \ i \text{ fixed in } \ 1,\dots, N,
\ee 
and the equilibrium profile in each vertical strip $\Q^E_{ij}$, $i$ fixed and $j=1,\dots, M$ will be defined accordingly.
This procedure returns the exact value of the free surface in the case of the water at rest test case (even in presence of obstacles), and approximates its value in the other situations, allowing the application of our method also far from the equilibrium. 
We underline that the values of $y_0$ obtained through \eqref{eq.Free_surf_detector} are always the same in adjacent cells of type $T_{ij}^n\backslash T_{i,j+1}^n$, which is necessary for the definition of the first order scheme in the $y$-direction, but can be different in two cells of type $T_{ij}^n\backslash T_{i+1,j}^n$; however this fact does not affect the scheme in the $x$-direction which does not depend on the equilibrium profile.

This automatic detection of the free surface is particularly useful in the case of small perturbations around an equilibrium profile but also for dambreaks problem. 
In this last situation, in particular, the possibility of computing dynamically the free surface position allows first to maintain perfectly the two initial states of a dambreak problem until they are not reached by the dam waves and then to restore the equilibrium when the water level comes back to a stable configuration at rest.

\section{CUDA-based parallel implementation on GPU}
\label{sec.C_GPU.CUDA}

The simplified two-phase flow model discretized through our new second order well balanced method is particularly well suited to  simulate breaking waves and dambreak problems: in particular, the free surface is perfectly described both in equilibrium situations as well as in the case where the free surface is a multi-valued function, e.g. in the case of splashing, 
breaking waves and overtopping flows. 
However, to capture all the details of a complex free surface flow we need also a very fine mesh.
Thus, in order to increase the computational efficiency of our method, we have decided for a parallel implementation, choosing in particular a GPU-based platform. 

The NVIDIA CUDA (Compute Unified Device Architecture) framework is a hardware and software platform that allows to easily exploit the NVIDIA GPUs and parallelize many costly algorithms. We refer to the recent book \cite{mantas2016introduction} for a complete introduction to CUDA for Scientific Computing. The images of this section are indeed courtesy of its authors.
Here we limit ourselves to a small introduction and to some notes on our specific implementation.

\begin{figure}
	\begin{center}
		\includegraphics[width = 0.45\linewidth]{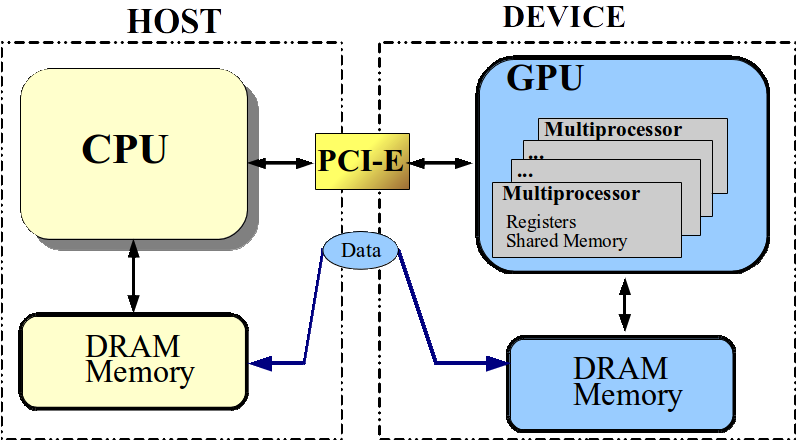} 
	\end{center}
	\caption{GPU device: coprocessor for the host system.}
	\label{fig.GPUdevice}
\end{figure}

According to the CUDA framework, a GPU is viewed as part of a computing
device, with their own Dynamic Random Access (DRAM) memory, which works
as a coprocessor for a host which includes the CPU and its (DRAM) memory, see
Figure \ref{fig.GPUdevice}. It is possible to copy data from the host DRAM memory (linked to CPU)
to the DRAM memory of the device (linked to GPU) and viceversa.

\begin{figure}
	\begin{center}
		\includegraphics[width = 0.45\linewidth]{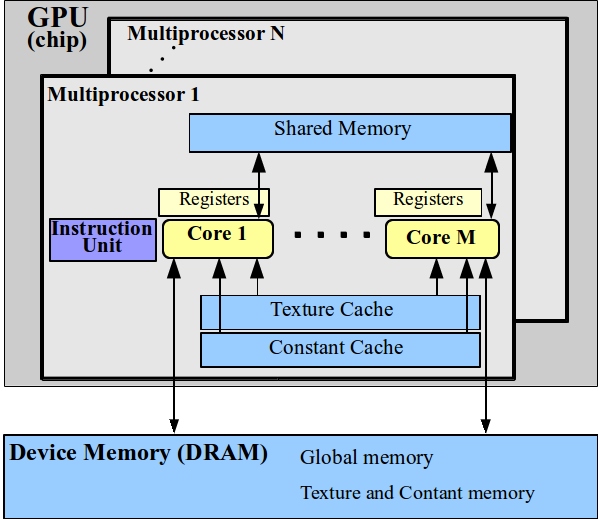} 		
	\end{center}
	\caption{CUDA hardware model.}
	\label{fig.GPUhardware}
\end{figure}

The CUDA framework assumes a unified architectural view of the GPU where a
GPU is formed by many multiprocessors, each one having multiple cores (see Figure \ref{fig.GPUhardware}).
At any clock cycle, each core of the multiprocessor executes the same instruction, but operates on different
data (SIMD paradigm, i.e. single instruction, multiple data). 
The data which are processed in each multiprocessor can be stored in the \textit{global memory}, that can be shared by all execution units but is slower than the other type of memory i.e. \textit{registers} and \textit{shared memory} which belong respectively to each core and each multiprocessor.

\begin{figure}
	\begin{center}
		\includegraphics[width = 0.5\linewidth]{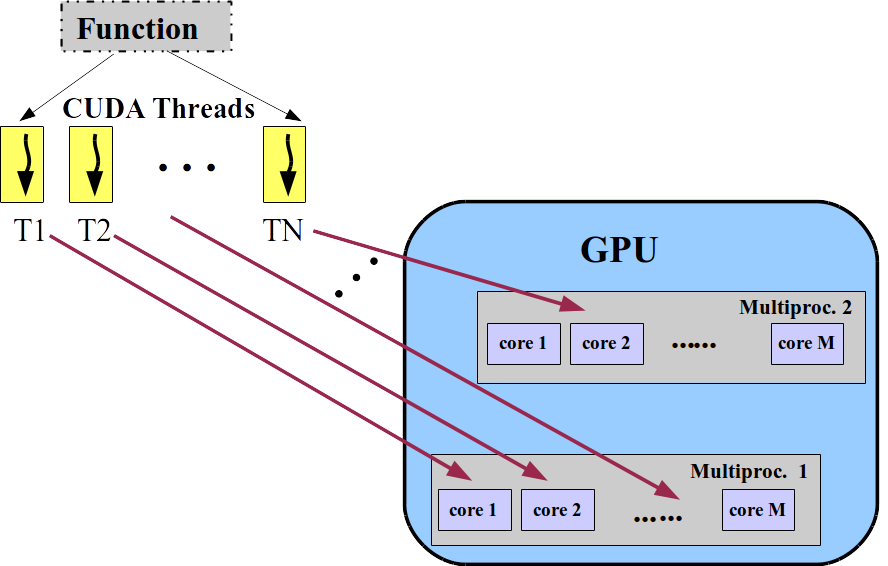} 
	\end{center}
	\caption{Execution of multiple CUDA threads associated with the same kernel.}
	\label{fig.GPUthreads}
\end{figure}

\begin{figure}
	\begin{center}	
		\includegraphics[width = 0.35\linewidth]{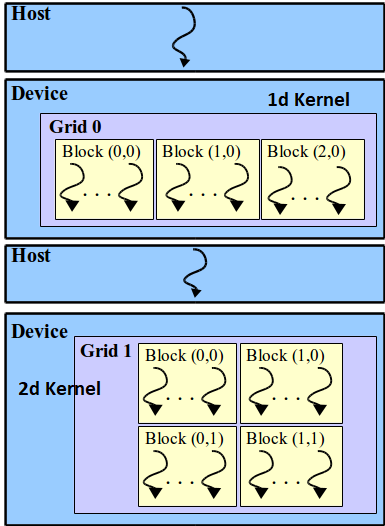} 
	\end{center}
	\caption{ Structure of a CUDA Program. A kernel is a function to be executed by each thread (represented by the curved arrow). According to the data structure the threads are organized in a $1D$, $2D$ or $3D$ grid.}
	\label{fig.GPUgrids}
\end{figure}

The multiprocessors of a GPU are specialized in the parallel execution of a huge
number of CUDA threads. A CUDA thread represents a sequential computational
task which executes an instance of a function and is executed logically in parallel
with respect to other CUDA threads (associated to the same function but operating
on different data) on the cores of a GPU multiprocessor (see Figure \ref{fig.GPUthreads}).
To specify the function to be executed by each thread on the GPU, the programmer
must define a special C function, called CUDA \textit{kernel}. A CUDA kernel is called
from the CPU and is executed on the GPU by each CUDA thread.
When a kernel is launched, every thread executes the same function but the processed data
depend on the value of several built-in variables which identify
the position of each particular thread, from which one can deduce for example the element or the edge on which the thread is operating.
According to the data structure the threads are organized in a $1D$, $2D$ or $3D$ grids, see Figure \ref{fig.GPUgrids}. In this way data referring to neighbor cells in the physical problem will be grouped together even during the computation contributing in accelerating the algorithm, since the interactions between neighbors are frequent and the proximity of the data reduces the  memory accesses.

\subsection{Notes on our implementation}

Our algorithm is composed by three main kernel functions.
The first one is used to reconstruct the polynomial representing the conserved variables inside each cell. Data are organized in a $2D$ grid and each thread works on one element of the mesh. 

The second one is used to solve the Riemann problem at each interface, i.e. each edge of the mesh. 
We compute the numerical flux between two elements across an edge: the obtained quantity should be added (with different signs) to each of the two elements.
As in any parallel implementation, it is important to guarantee that there are no superpositions, i.e. that we do not modify at the same time the same element. 
So we cannot simply launch the kernel function indiscriminately on all the edges at the same time, because each element has 4 edges and so many possibilities of superpositions.
To overcome this problem the idea consists simply in subdividing the edges in four groups. 
Since the mesh is Cartesian we can distinguish the edges between vertical and horizontal ones and we can number each group incrementally form left to right and from the bottom up. The four groups are then obtained by considering: the vertical edges with odd numbering, the vertical edges with even numbering, the horizontal edges with odd numbering, the horizontal edges with even numbering.
So the computation of the flux at the interfaces is subdivided into four kernel functions each one operating on a non intersecting set of data.

The third kernel function is organized as the first one and is used to update the solution in each cell from time $t^n$ to time $t^{n+1}$ summing up all the flux contributions.

\section{Numerical Results}
\label{sec.C_GPU.NumResults}

The numerical results presented in this section show that our second order well balanced method preserves the equilibrium solutions up to machine precision and at the same time can also solve complex free surface flows, both with and without gravity, always with an efficiency of ten/twenty millions of volumes treated per seconds.
In particular we will take into account three groups of test cases.
\begin{enumerate}
	\item
	First, in Section \ref{sec.EquilibriaResults} we will show the ability of our second order well balanced scheme in maintaining the equilibrium profiles of water at rest. 
	\item
	Then, we will analyze problems far away from the equilibrium, with a velocity different from zero and a free surface position that cannot be described by a single-valued function already at time $t=0$. This kind of test cases, presented in Sections \ref{sec.Drop}-\ref{sec.JetImpinging}, are studied without gravity, i.e. by taking $g=0$ and hence with a zero source term. It is clear that in these cases the use of a well balanced scheme would not be required, but however this does not prevent the  proposed well balanced scheme from working with robustness and efficiency. Moreover, it is useful to emphasize one main property of this kind of schemes: indeed, the fact of preserving equilibria up to machine precision and of being more accurate than standard  schemes \textit{close} to steady states does \textit{not} worsen its efficiency and accuracy with respect to standard schemes  of the same order, when being \textit{far away} from the equilibria.  
	\item
	Last, in Sections \ref{sec.FlowOver}-\ref{sec.Impacttest} we will concentrate on dambreak type problems. In this group of test the gravity is switched on with value $g=9.81$ and the initial density and pressure profiles are taken in equilibrium at time $t=0$ inside the liquid domain by using \eqref{eq.rho_in_liquid}; in the gas phase instead we employ the constant values $p=0$ and $\rho=\rho_0$. Thanks to our well balanced scheme a liquid cell remains stationary until the effects of the dam removal propagate up to that cell. 
\end{enumerate}

The great variety of the presented tests is intended to show both the wide range of applicability of the proposed well balanced scheme and the completeness of the studied two-phase model. 
In particular, we emphasize the capability of this model in capturing complex free surface flows that cannot be described by a  single-valued function (see for example \ref{sec.SpinningSquare} and \ref{sec.Impacttest}), flows where vertical accelerations  
are rather important, as in the first stages after a dambreak, and recirculation phenomena that can be noticed for example in  Figure \ref{fig.Dambreak_withStep_wet_zoom_2}. 
Moreover, we would like to underline that these applications are of real interest and have been longly studied analytically, numerically and carrying out numerous experiments. For this reason we can compare the obtained numerical results against analytical or numerical solutions of the shallow water system, reference solutions obtained with the smooth particle hydrodynamics (SPH) scheme or through experiments, and the numerical results obtained in \cite{dumbser2011simple} for the same two-phase system 
studied with a third order accurate but not well-balanced ADER-WENO finite volume scheme.

In all the case our numerical results are in perfect agreement with the experiments and the SPH results. Moreover, it is worthy to notice that our second order well balanced scheme shows an excellent correspondence with the results obtained in \cite{dumbser2011simple} with a third order method. Finally, for long time simulations we recover also the shallow water free surface profile, which instead it is rather different at early times, as expected .

The majority of the proposed tests have been taken from \cite{dumbser2011simple} and in particular for the choice of the constant $k_0$ we refer to it. The other parameters of the model are chosen as $\epsilon=10^{-3}$, $\rho_0=1000$, $\gamma=1$ and the initial condition inside the air domain $\Omega_a$ is always set to $\u=\0$, $\rho=\rho_0$.
If not otherwise specified, the reflective wall boundary conditions are implemented by assigning a state at the wall boundary which solves the inverse Riemann problem at the element interface such that the normal velocity $\u\cdot \mbf{n} = 0$.

\subsection{Equilibria} 
\label{sec.EquilibriaResults}

In this section, we show the ability of our second order well balanced scheme in maintaining the water at rest equilibrium profile up to machine precision for very long times, both with a flat bottom topography and in the presence of obstacles.
We consider an equilibrium profile characterized by \eqref{eq.BNs_eq1}, \eqref{eq.BNs_eq3} and \eqref{eq.BNs_eq4}, with $k_0=2.78\cdot 10^5$, $\rho_0=1000$, $\alpha=1$ and $\gamma=1$ on a domain $\Omega=[-0.5,0.5]\times[0,1]$.
In the first test case we consider a flat bottom topography so $\Omega_\ell = \Omega$. In the second case instead we model the presence of an obstacle by removing a portion of the domain, i.e. 
$\Omega=\Omega_\ell=([-0.5,0.5]\times[0,1])\backslash([-0.25,0.25]\times[0,0.33])$ whose boundary are treated  with reflective boundary conditions.

\begin{figure}
	\begin{center}
		\includegraphics[width = 0.35\linewidth]{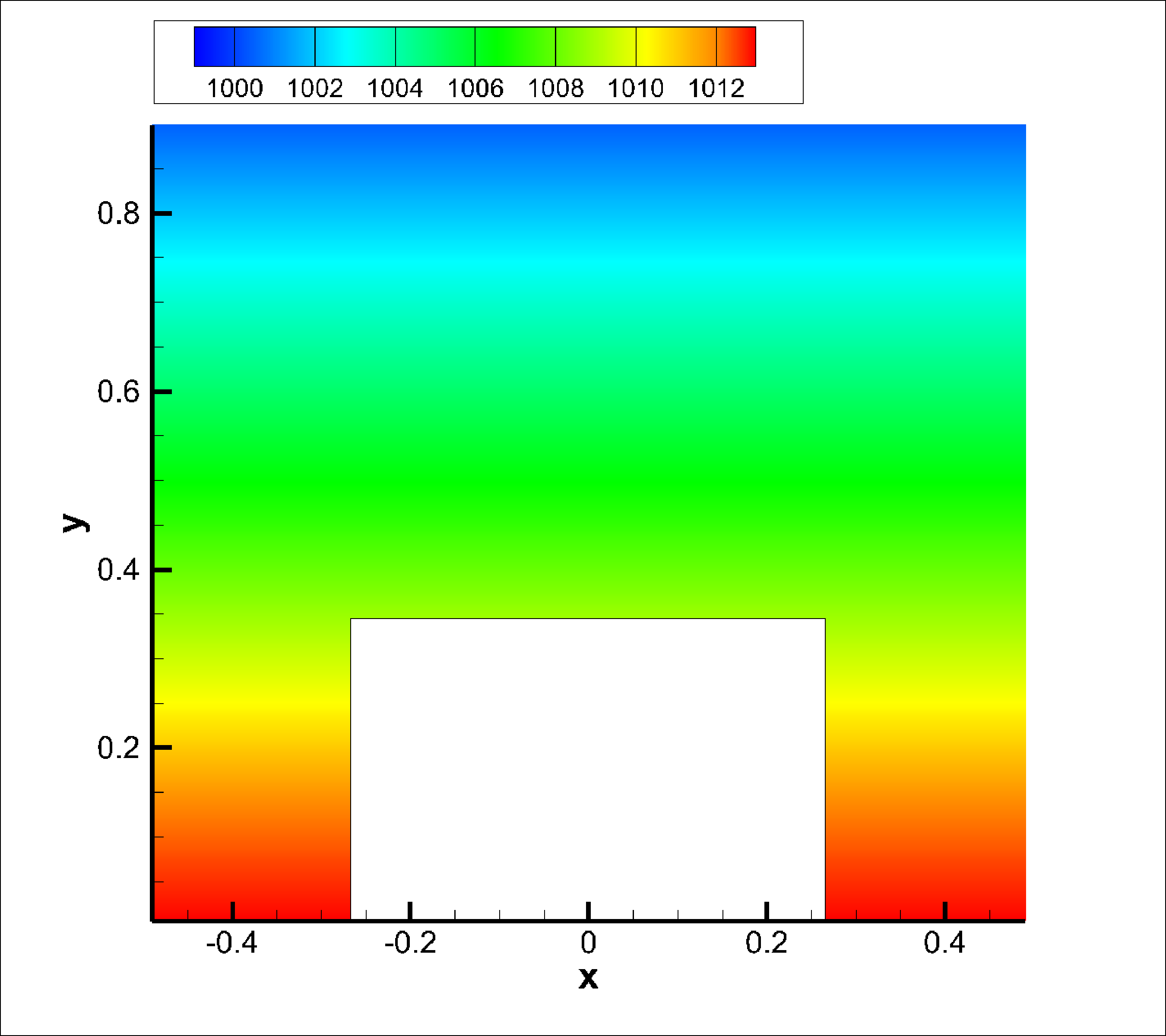} \qquad
		\includegraphics[width = 0.35\linewidth]{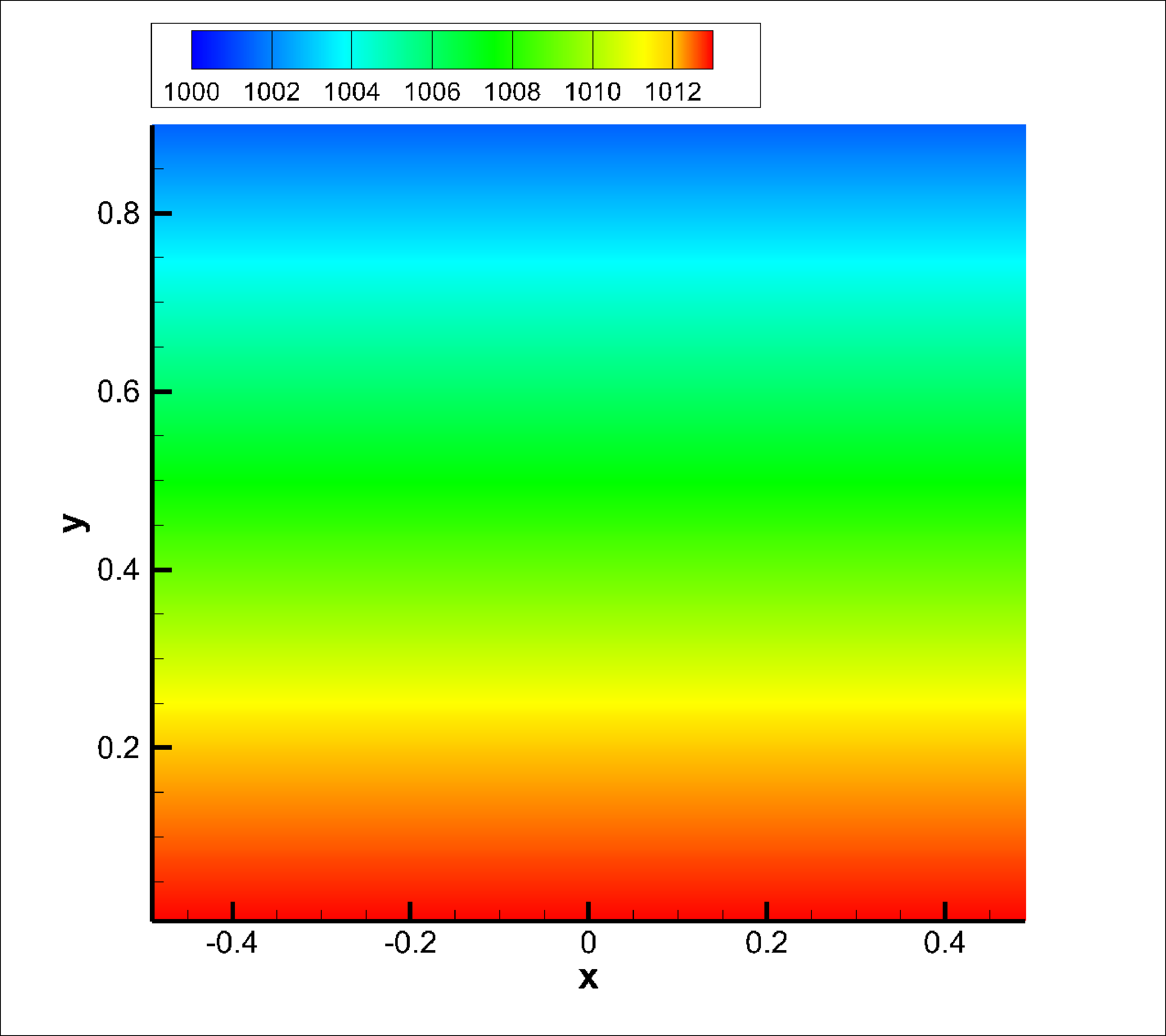}
		\caption{Contour colors of the quantity $\alpha \rho$ at the equilibrium with (left) and without (right) obstacle.}
		\label{fig.BNsimp.Equilibria}
	\end{center} 
\end{figure}

The considered equilibrium profiles are shown in Figure \ref{fig.BNsimp.Equilibria}. Our scheme preserves the equilibria  up to machine precision: in Table \ref{tab.BNsimp_equilibria} we report the maximum error between our numerical results, and the exact solution for $rho$, $u$, $v$, and $P$ in the two studied cases. We remark that the values of density and pressure depend on $k_0$ which is of the order of $10^5$ and in double precision the machine precision referred to $10^5$ is of the order of $10^{-11}$.

\begin{table} 
	\begin{center}	
		\caption{Maximum error between the numerical and the exact solution for $\rho$, $u$, $v$ and $P$ with (left) and without (right) obstacle. The numerical results have been obtained with our second order well balanced Osher-Romberg on a mesh of $100 \times 100$ elements. Equilibria are maintained up to machine precision for very long times. }
		\label{tab.BNsimp_equilibria}	
		\begin{tabular}{c|cccc} 		
			 time  & $E_{\rho}$ & $E_{u}$ & $E_{v}$ & $E_{P}$  \\
			\hline
			1    & 4.55E-13 &  1.00E-13 & 8.99 E-13  & 2.91E-10\\
			10   & 5.57E-12 &  1.43E-12 &  2.55E-12  & 3.63E-09 \\
			50  & 7.04E-12 &  4.46E-12 &  4.84E-12  & 5.50E-09\\
			100 & 8.41E-11 &  5.38E-12 &  5.80E-12  & 7.74E-09   \\		
			\hline 
		\end{tabular}	\qquad  
		\begin{tabular}{c|cccc} 		
			time  & $E_{\rho}$ & $E_{u}$ & $E_{v}$ & $E_{P}$  \\
			\hline
			1    & 5.38E-12 &  0 &  5.71E-13  & 3.61E-10\\
			10   & 5.50E-12 &  0 &  6.83E-13  & 3.62E-10\\
			50  & 5.80E-12 &  0 &  7.42E-13  & 3.70E-10 \\
			100 & 6.41E-12 &  0 & 1.93E-13  &  4.12E-10 \\		
			\hline 
		\end{tabular}	
	\end{center}
\end{table}

\subsection{Elliptical drop}
\label{sec.Drop}

\begin{figure*}
	\begin{center}
		\includegraphics[width = 0.30\linewidth]{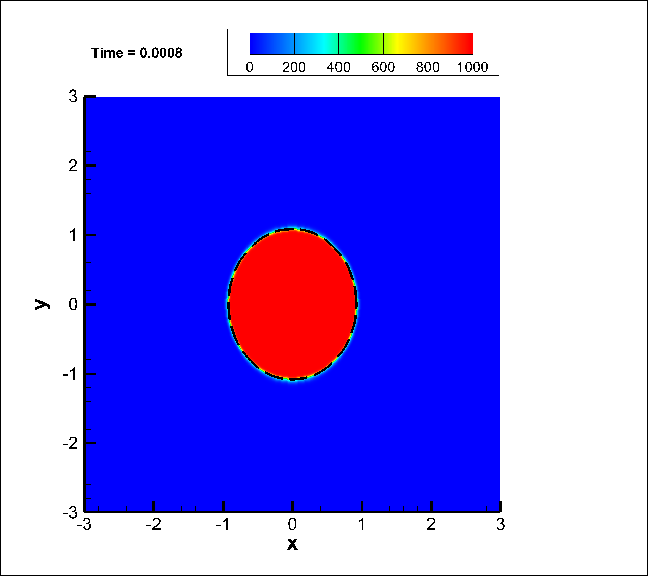}  \ 
		\includegraphics[width = 0.30\linewidth]{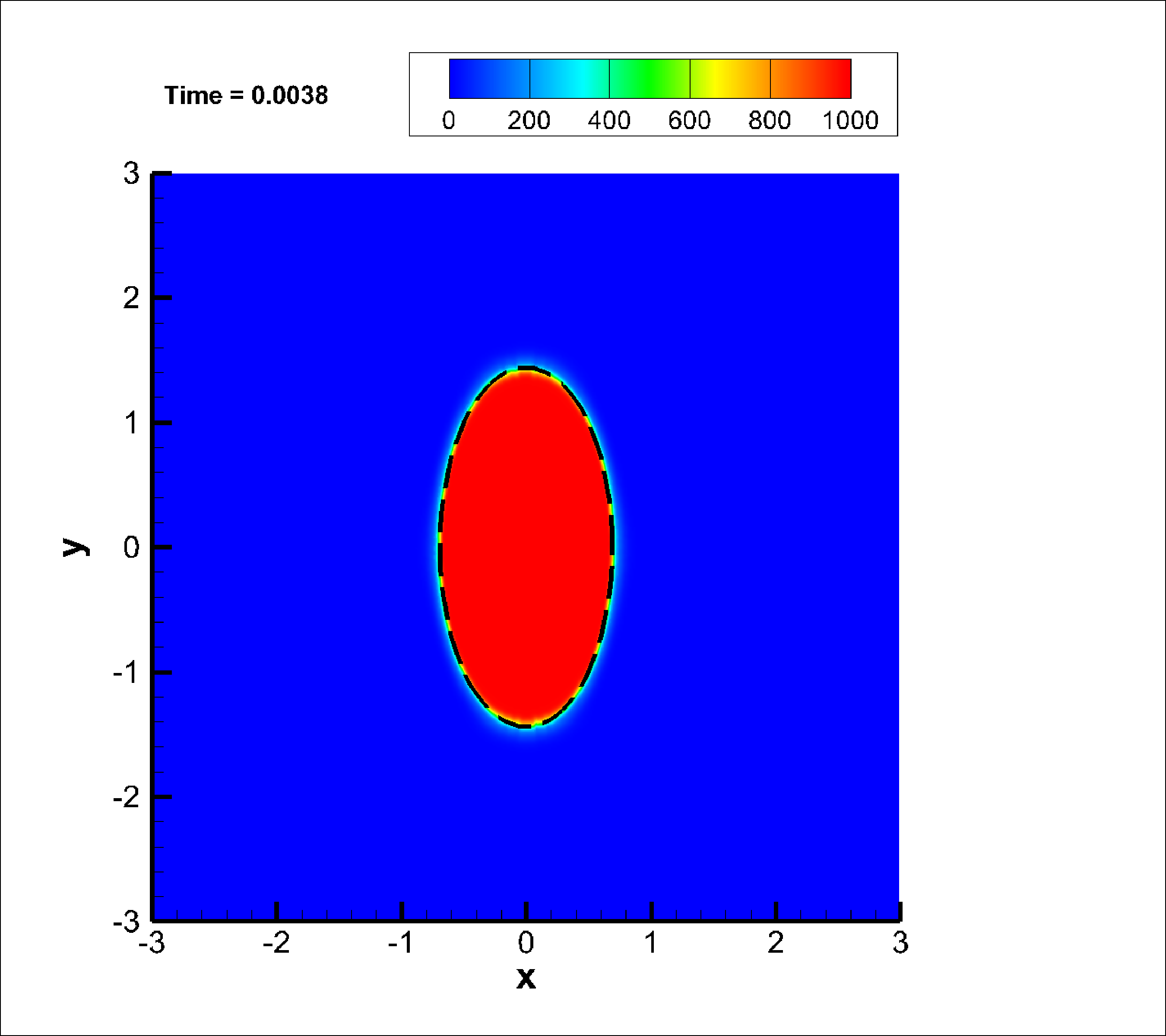}  \ 
		\includegraphics[width = 0.30\linewidth]{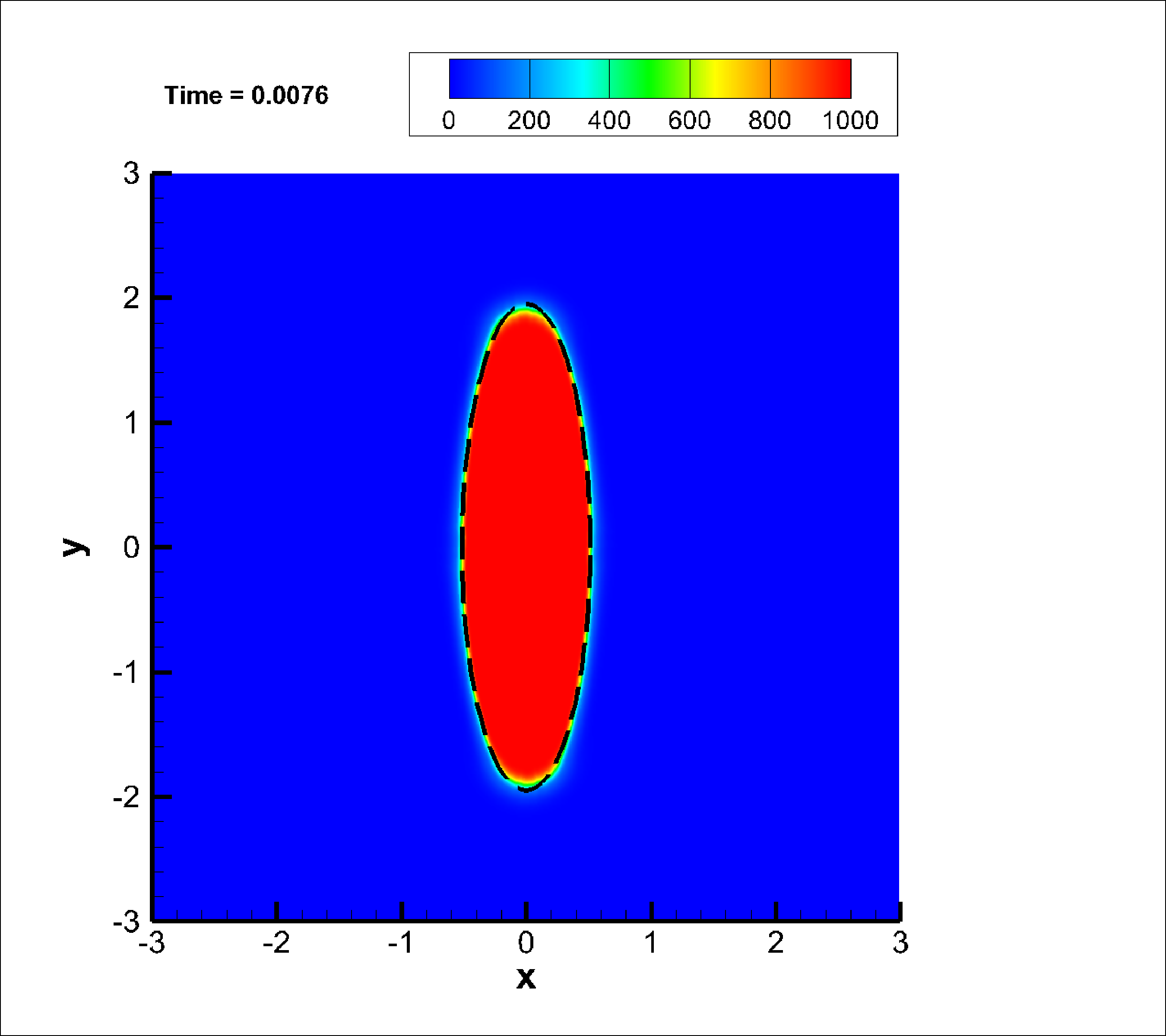}  \!\!\!\!\!
		\caption{Elliptical drop test case at time $t=0.0008$ (left), $t=0.0038$ (middle) and $t=0.0076$ (right). The density contour represents the quantity $\alpha\rho$ and the dashed line represents the exact analytical solution.}
		\label{fig.Drop}
	\end{center} 
\end{figure*}

Here, we validate our well balanced Osher-Romberg scheme for the 2D diffuse interface model in the case of $g=0$ and a liquid domain completely surrounded by the gas, where the free-surface cannot be described with a single-valued function. 
As computational domain we take a square $\Omega=[-3,3]\times[-3,3]$; the liquid domain $\Omega_\ell$ consists of a drop which at time $t=0$ has a circular shape of radius $r=1$ centered at the origin.
Inside the liquid domain we choose $\rho = \rho_0$, $\alpha = 1 - \epsilon$ with $\epsilon = 10^{-3}$ and a velocity field $\u=(-100x, 100y)$; the sound speed is chosen to be close to the real sound speed of water $(c=1500)$ by taking $k_0=2.25\cdot10^9$.
We solve the problem with our second order well balanced scheme by discretizing $\Omega$ with a Cartesian mesh of $400\times400$ elements. We would like to underline that in this test the initial condition is not in equilibrium and also that a well balanced scheme is not strictly necessary since with $g=0$ there are no source terms; however the scheme works and shows an excellent agreement with the exact analytical solutions, as shown in Figure \ref{fig.Drop}. 
The exact solution of this problem is given in \cite{Monaghan1994}, where the test was proposed for the first time to validate the smooth particle hydrodinamics (SPH) scheme. The exact solution consists in an elliptical drop with constant area $ab=1$, with $a$ and $b$ respectively the minor and the major half-axis. In particular, according to \cite{Monaghan1994} $b=1.083$ at time $t=0.0008$, $b=1.44$ at time $t=0.0038$ and $b=1.95$ at time $t=0.0076$.

To reach the final time $t=0.0076$ the scheme requires $4110$ time steps, so the total number of processed volumes is $\texttt{nVol} = 400 \times 400 \times 4110 = 6.576\cdot 10^8$; by employing an Nvidia GeForce Titan Black GPU, 
the algorithm is executed in $45$ seconds, with an average of $1.5\cdot10^7$ volumes processed per second.  


\subsection{Spinning square}
\label{sec.SpinningSquare}

\begin{figure*}
	\begin{center}
		\includegraphics[width = 0.24\linewidth]{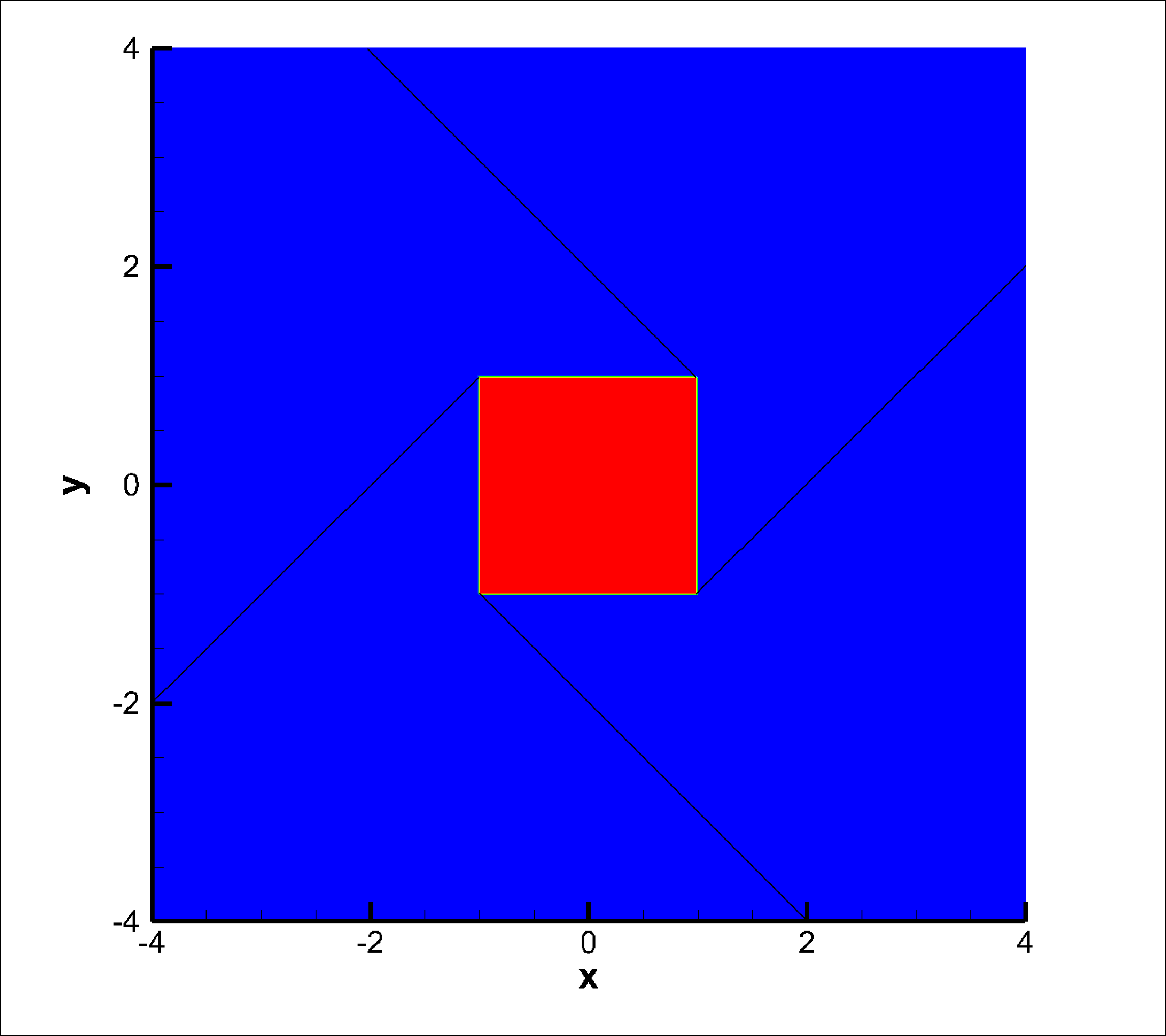} 
		\includegraphics[width = 0.24\linewidth]{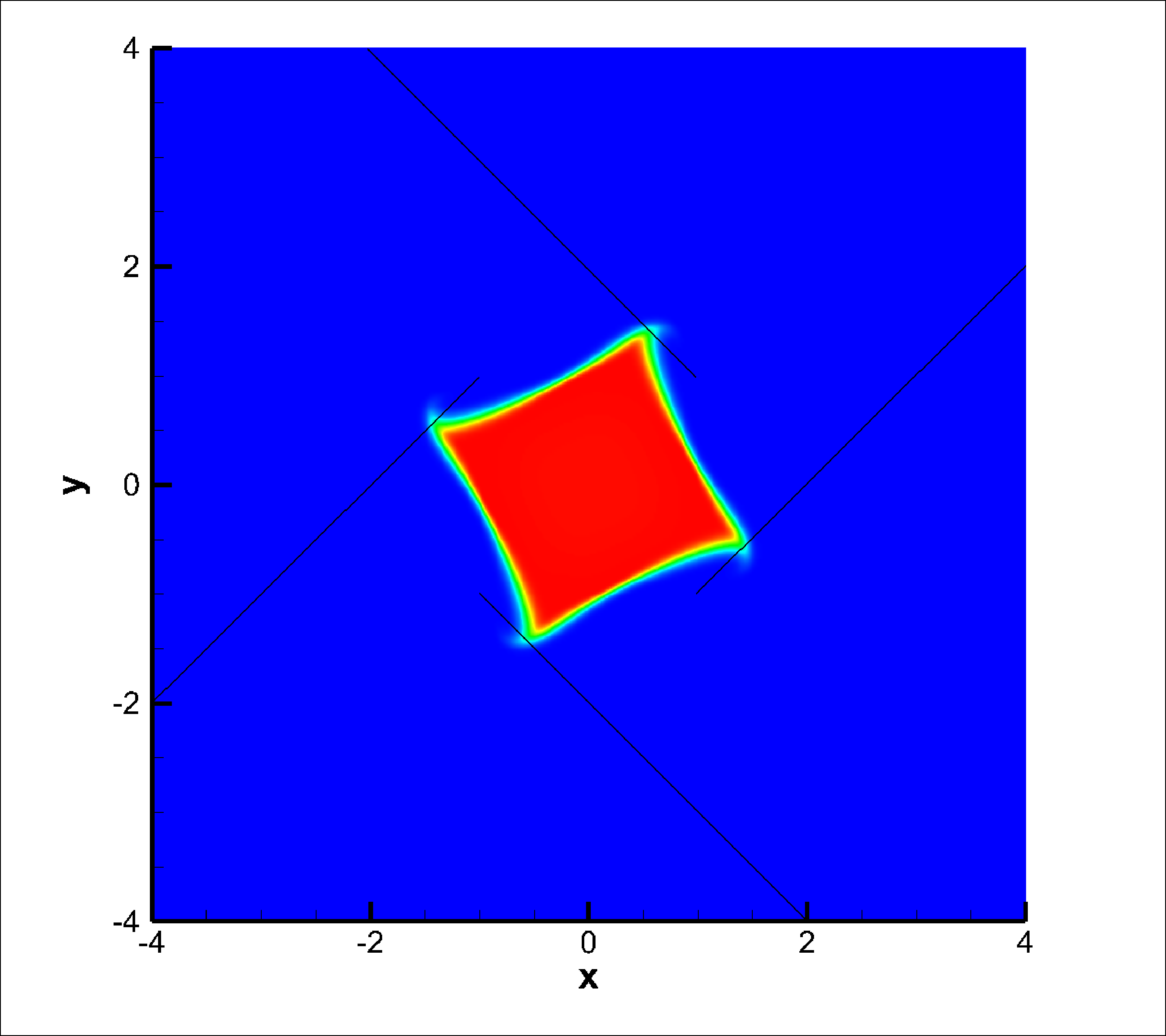} 
		\includegraphics[width = 0.24\linewidth]{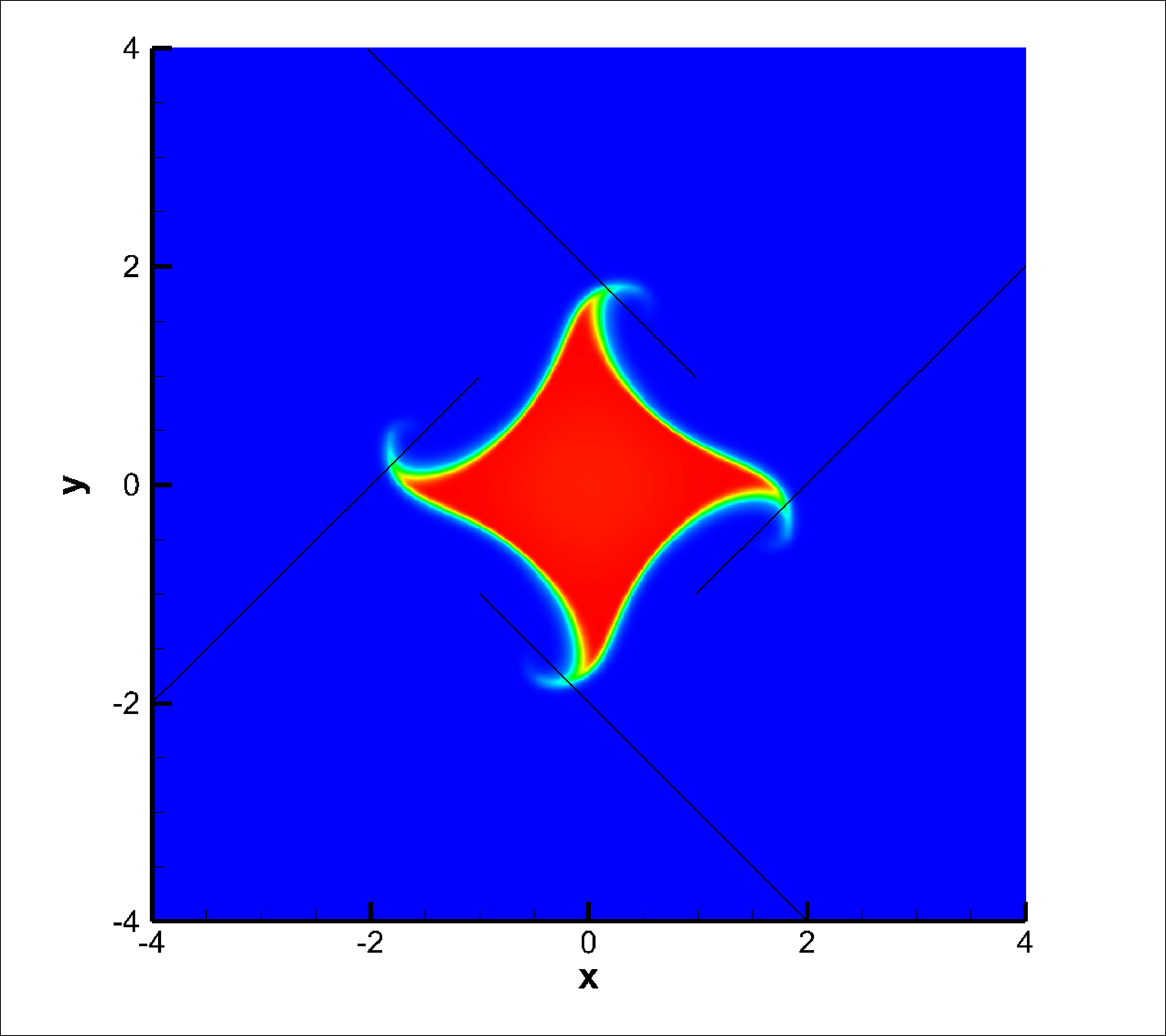} 
		\includegraphics[width = 0.24\linewidth]{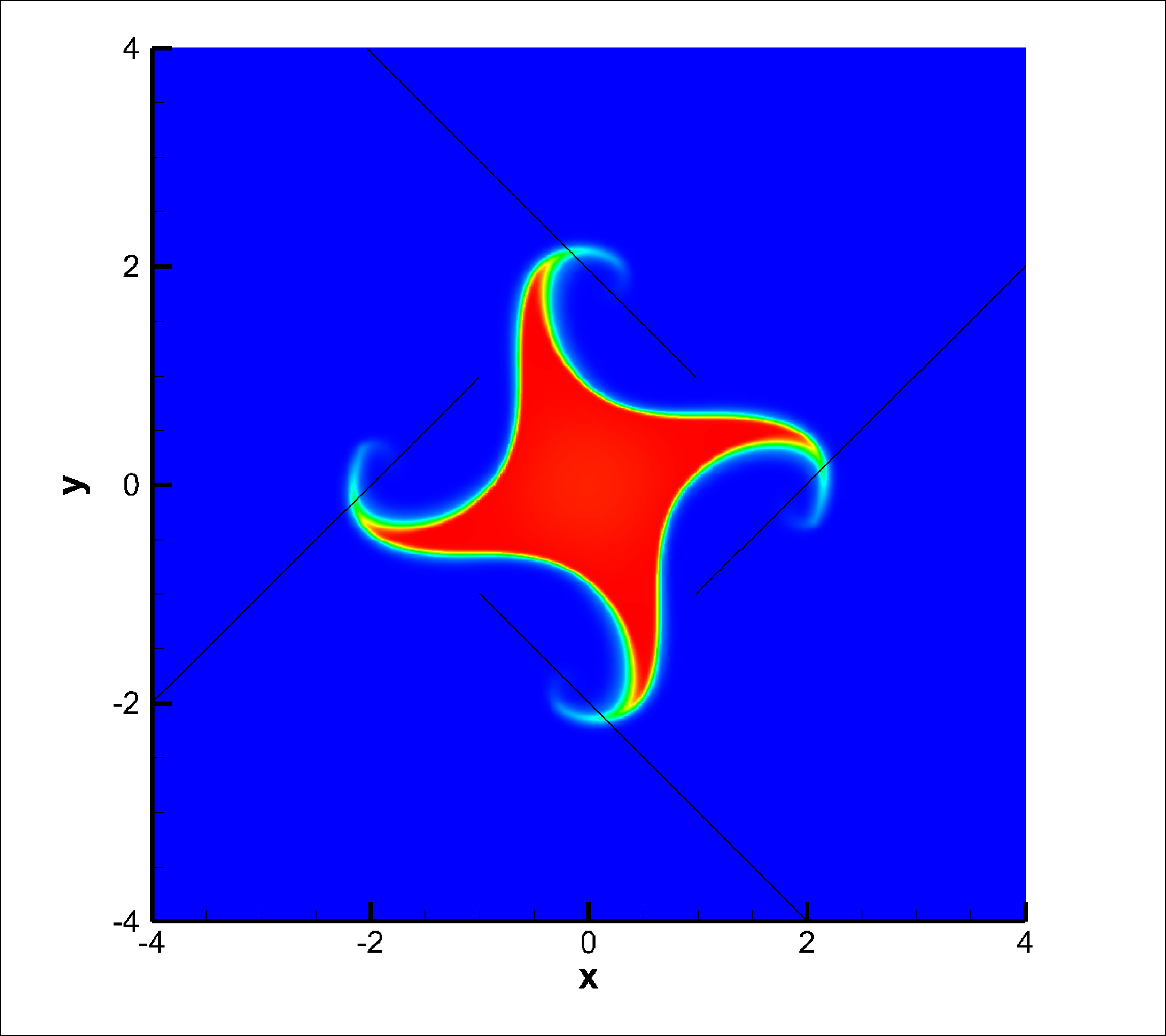}\\[1pt] 
		\includegraphics[width = 0.24\linewidth]{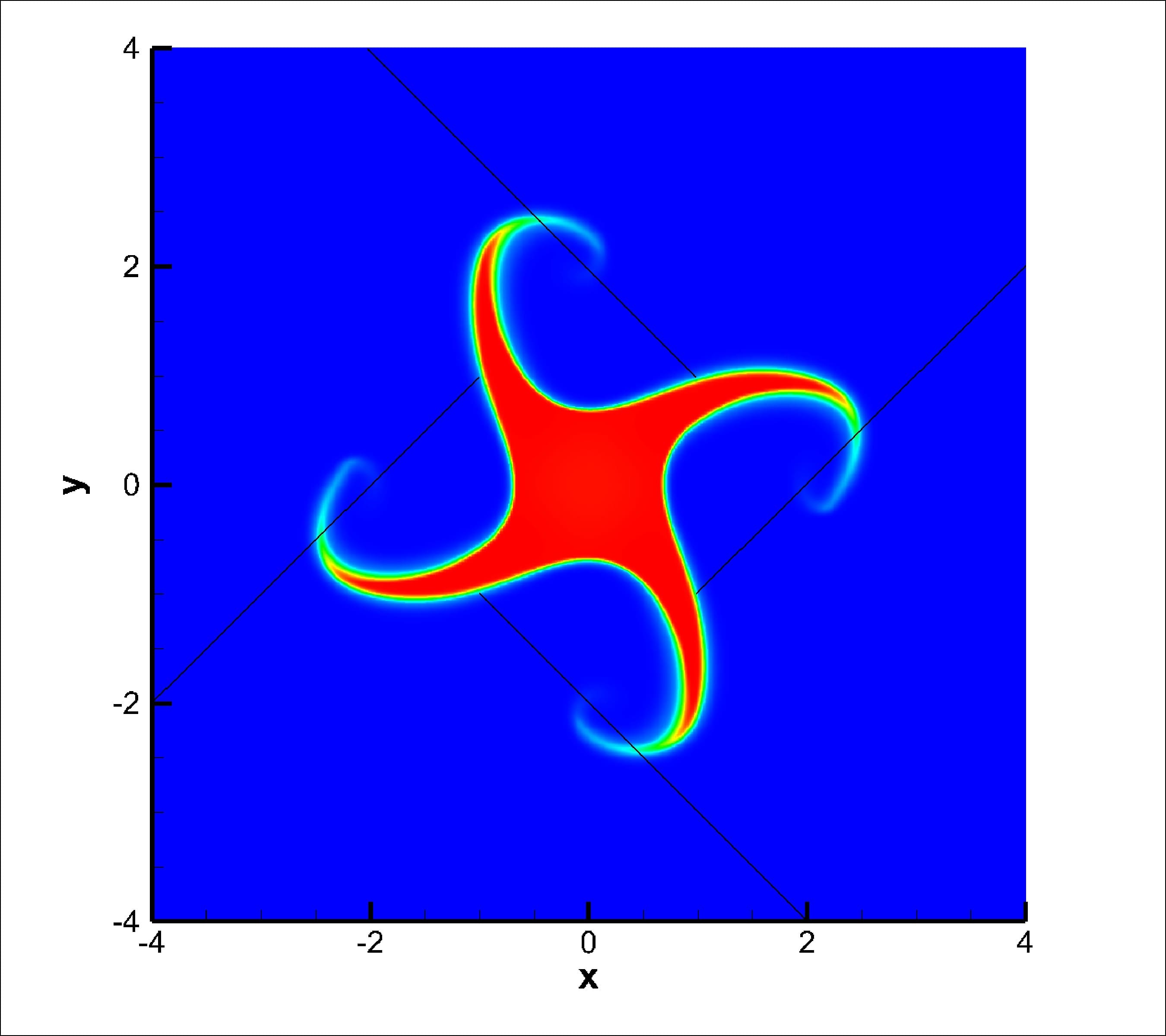} 
		\includegraphics[width = 0.24\linewidth]{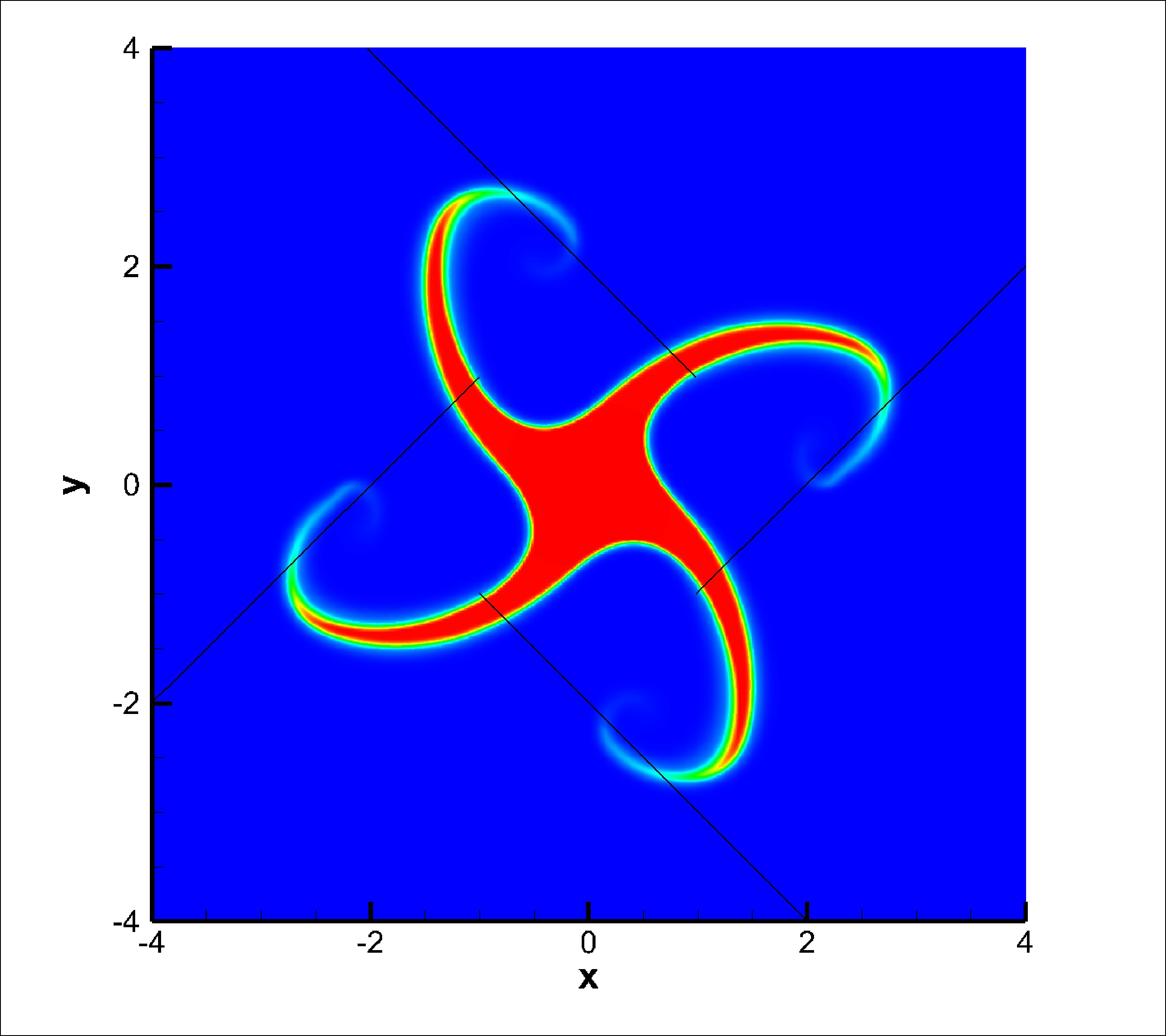}
		\includegraphics[width = 0.24\linewidth]{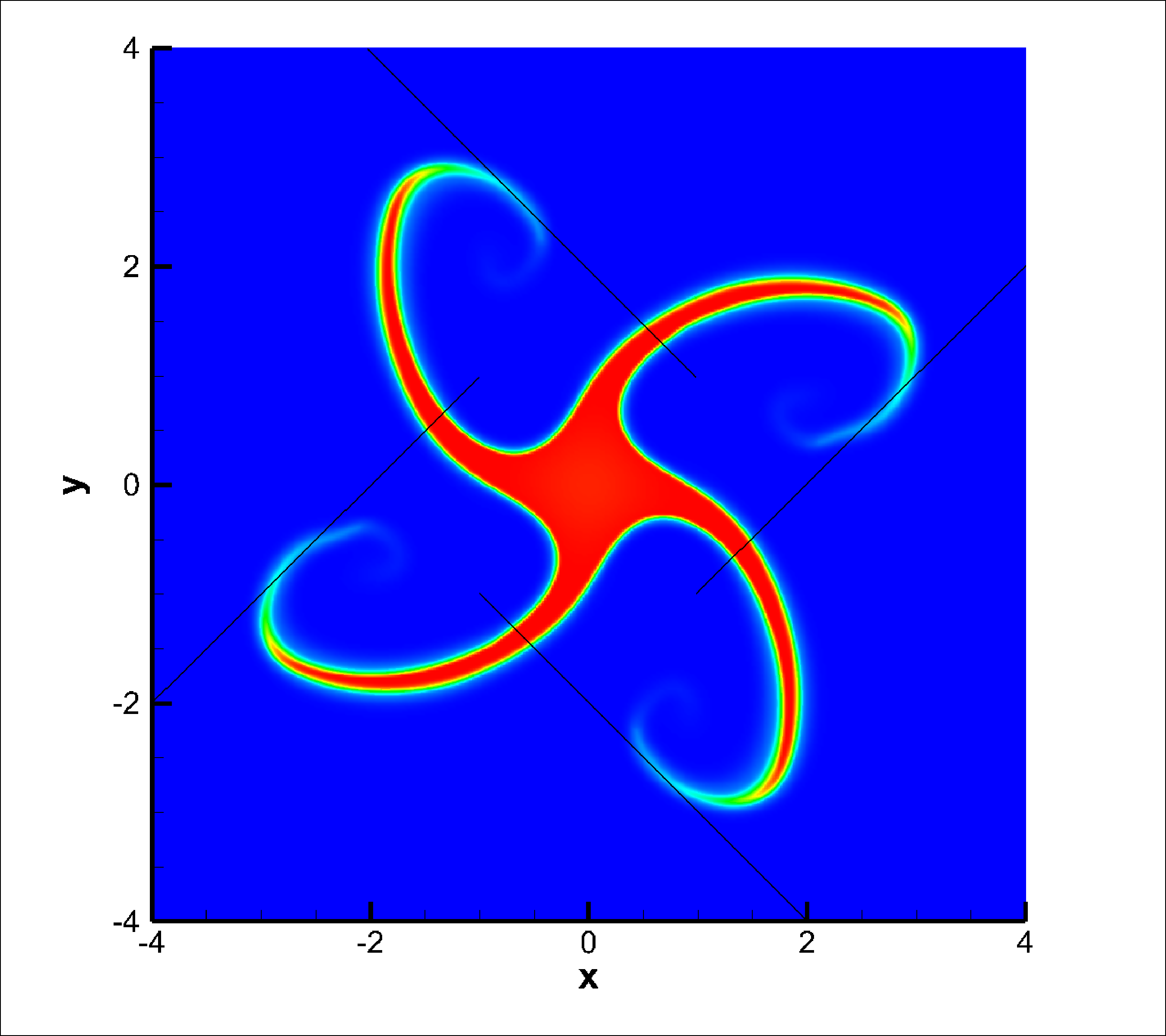} 
		\includegraphics[width = 0.24\linewidth]{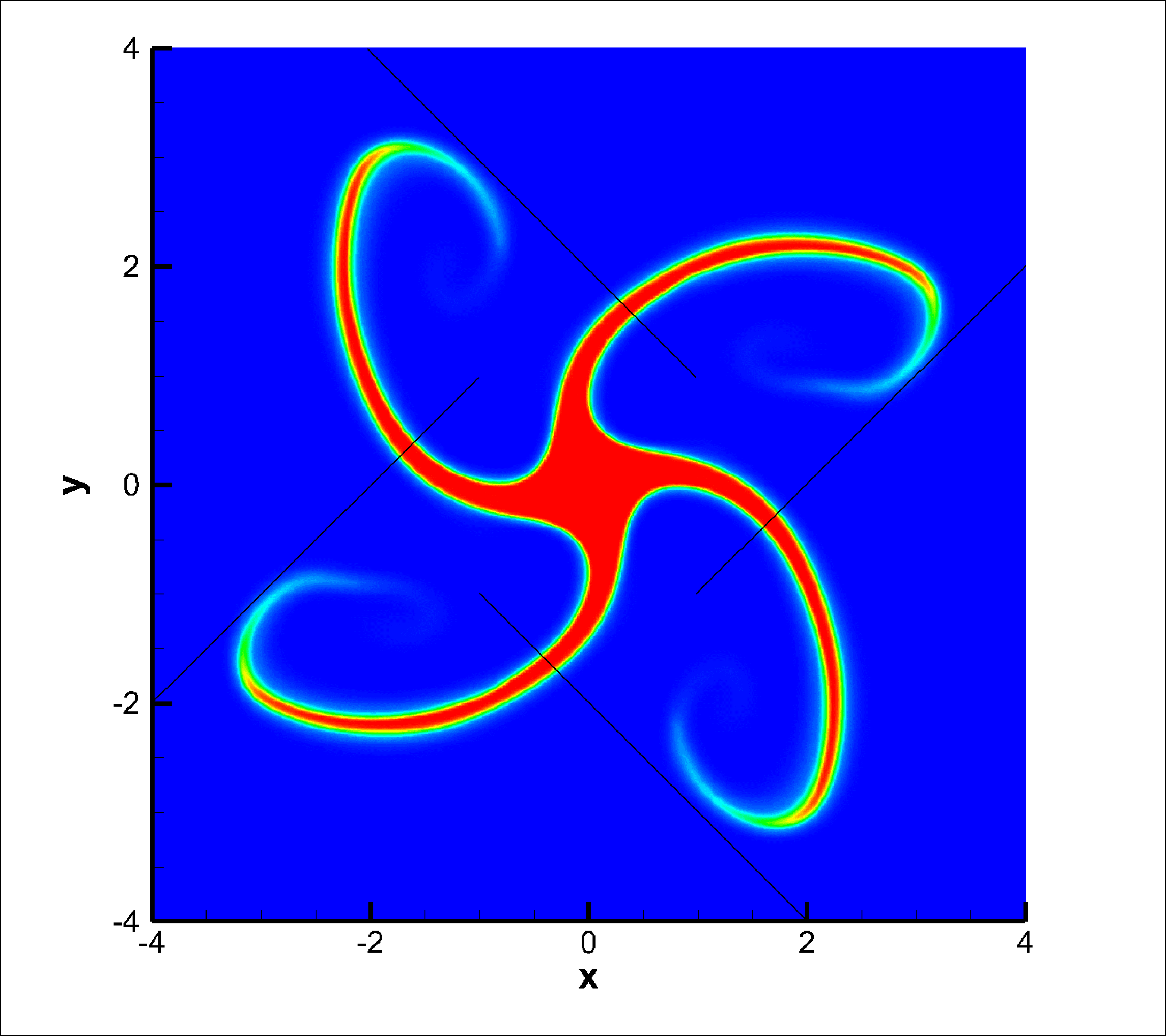}\\
		\caption{Spinning square test case. In the figure we depict the density contour of $\alpha\rho$ at times $t=0, 0.5/\omega, 1/\omega, 1.5/\omega, 2/\omega, 2.5/\omega, 3/\omega, 3.5/\omega$ with $\omega =2\pi$. The black straight lines represent the theoretical trajectories expected for the evolution of the square corners.}
		\label{fig.SpinningSquare}
	\end{center} 
\end{figure*}

Another test with characteristics similar to the previous one is the spinning square. We consider a computational domain $\Omega = [-5,5] \times [-5,5] $ covered with a Cartesian mesh of $850 \times 850$ elements. The liquid domain $\Omega_\ell$ at time $t=0$ is confined into the square $[-1,1] \times [-1,1]$, where we impose a density $\rho = \rho_0$, $\alpha = 1 - \epsilon$ with $\epsilon = 10^{-3}$ and a velocity field $\u=(2\pi y, -2\pi x)$ corresponding to a rigid-body rotation at angular frequency $\omega=2\pi$; then we choose $k_0= 8.78 \cdot 10^5$ to obtain a Mach number of $M=0.3$ based on the maximum velocity in the corners of the square. 
In Figure \ref{fig.SpinningSquare} we report the obtained results at successive instances until $t=3.5/\omega$. The rather complex free-surface given by our two phase model is in agreement with the theoretical trajectories describing the evolution of the corner position (depicted with black straight lines) and does not show any problem due to tensile instabilities which instead is a common issue with traditional SPH schemes. One can refer to \cite{ColagrossiPhD} and \cite{Oger} for a complete presentation of this test and to \cite{GingoldMonaghan} for a possible cure of instabilities in the SPH context. 

To reach the final time $t=3.5/\omega$ the scheme requires $9813$ time steps, so the total number of processed volumes is $\texttt{nVol} = 850 \times 850 \times 9813 = 7.09\cdot 10^9$; on an Nvidia GeForce Titan Black GPU the algorithm is executed 
in $447$ seconds, with an average of $1.6\cdot10^7$ volumes processed per second.

\subsection{A water jet impinging obliquely on a plate}
\label{sec.JetImpinging}

\begin{figure}
	\begin{center}
		\includegraphics[width = 0.7\linewidth]{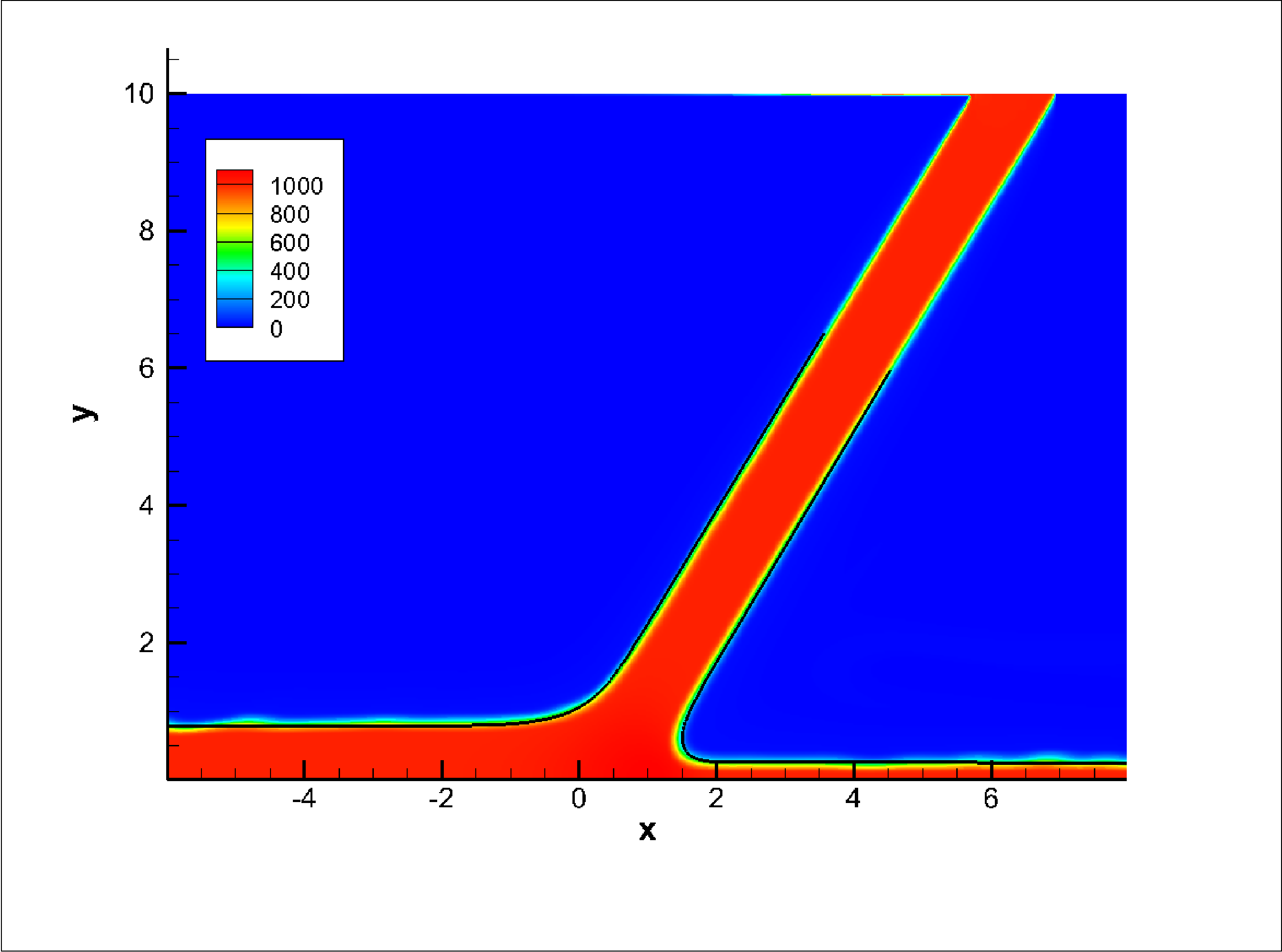} 
		\caption{Water jet impinging  obliquely on a plate. We depict the density contour of the quantity $\alpha\rho$ obtained at time $t=6.5$ and the exact analytical solution for steady state \eqref{eq.Jet_freesurface_exact_A} and \eqref{eq.Jet_freesurface_exact_B} (black line).}
		\label{fig.JetFreeSurface}
	\end{center} 
\end{figure}

\begin{figure}
	\begin{center}
		\includegraphics[width = 0.5\linewidth]{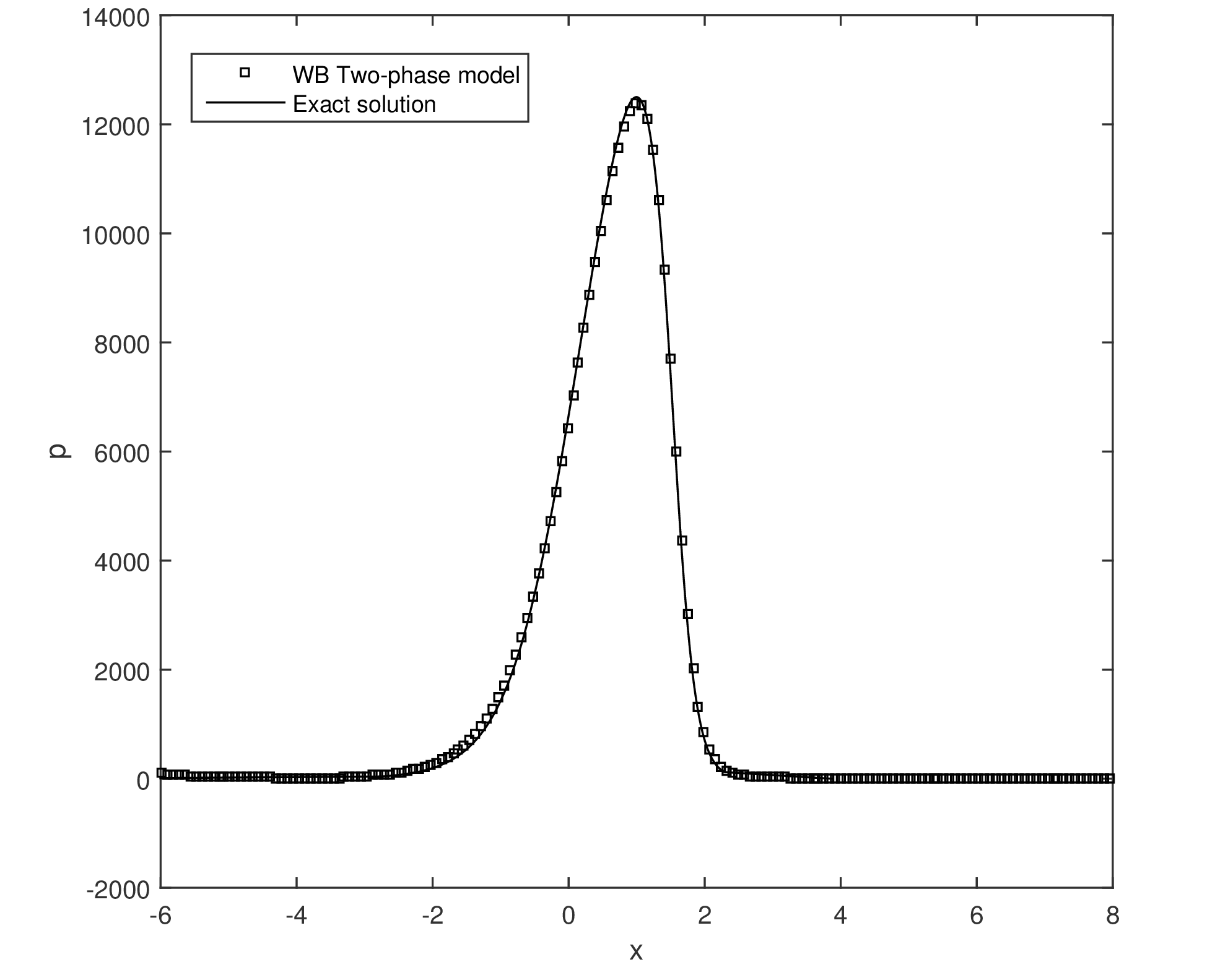} 
		\caption{Water jet impinging impinging obliquely on a plate. We depict the pressure profile along the flat plate at time $t=6.5$ (squares) compared against the exact solution }
		\label{fig.JetPressure}
	\end{center} 
\end{figure}

In this section we present another test problem without gravity, i.e. with $g=0$, that consists of a water jet impinging on a flat plane with an angle of $\theta=\ang{60}$.
The computational domain is $\Omega=[-6,8]\times[0,10]$ and is covered with a Cartesian mesh of $500\times350$ elements. At time $t=0$ we suppose that the liquid domain $\Omega_\ell$ has just impacted the flat plane and occupies the strip $\{ \x \in \Omega \,|\,  y-\sqrt{3} x \le 0 \wedge y-\sqrt{3}x \ge -2\}$, hence the jet has a thickness of $H=1$.
Inside the liquid domain we impose $\rho=\rho_0$, $\alpha = 1 - \epsilon$ with $\epsilon = 10^{-3}$ and a velocity field $\u=5 \left(-\cos({\theta}), -\sin({\theta})\right)$ and we use $k_0=2.78\cdot 10^5$ in order to obtain a Mach number of $M=0.3$.
We use reflective boundary conditions on any side of the domain except for $\{ 10/\sqrt{3} \le x \le 12/\sqrt{3} \wedge y=10\}$ where we assume an inflow boundary condition.
The jet impinges on the wall and gets reflected asymmetrically due to the incident angle $\theta$;  after $t=6.5$ a steady state is reached for the free surface and the pressure on the plate.
According to \cite{MilneThomson}, the exact solution for the free surface is given by the following two parts
\be
\label{eq.Jet_freesurface_exact_A}
\begin{cases}
	x =\!&\frac{1}{\pi} \biggl( \left(\theta-\pi\right)\sin\left(\theta\right)+\log\left(\tan\left(\frac{1}{2}\beta\right)\right)+\cos\left(\theta\right)\Bigl(   \log\left(\frac{1}{2}\sin\left(\beta\right)\right)  
	-\log\left(\sin\left(\frac{1}{2}\left(\theta+\beta\right)\right)\sin\left(\left(\theta-\beta\right)\frac{1}{2}\right)\right) \Bigr)    \biggr ),  
	\\[2pt]     
	y  =&\frac{1}{\pi}\biggl(  \frac{1}{2}\pi\left(1+\cos\left(\theta\right)\right)+ \sin\left(\theta\right)\Bigl(\log\left(\sin\left(\frac{1}{2}\left(\theta+\beta\right)\right)\right) 
	-\log\left(\sin\left(\frac{1}{2}\left(\theta-\beta\right)\right)\right)\Bigr)     \biggr), \qquad \quad \ang{0} < \beta < \theta, 
\end{cases}
\ee 
\be
\label{eq.Jet_freesurface_exact_B}
\begin{cases}
	x =&\frac{1}{\pi} \biggl( \theta \sin\left(\theta\right)+\log\left(\tan\left(\frac{1}{2}\beta\right)\right)+\cos\left(\theta\right)\Bigl(   \log\left(\frac{1}{2}\sin\left(\beta\right)\right) 
	-\log\left(\sin\left(\frac{1}{2}\left(\theta+\beta\right)\right)\sin\left(\left(\beta-\theta\right)\frac{1}{2}\right)\right) \Bigr)    \biggr ), 
	\\[2pt]
	y  =&\frac{1}{\pi}\biggl(  \frac{1}{2}\pi\left(1-\cos\left(\theta\right)\right)+ \sin\left(\theta\right)\Bigl(\log\left(\sin\left(\frac{1}{2}\left(\theta+\beta\right)\right)\right) 
	-\log\left(\sin\left(\frac{1}{2}\left(\beta-\theta\right)\right)\right)\Bigr)     \biggr), \qquad \quad \theta < \beta < \ang{90}.
\end{cases}
\ee 
In Figure \ref{fig.JetFreeSurface} the results obtained with our two phase model at time $t=6.5$ are compared with the above analytical solution and an excellent agreement can be noticed.
Moreover, in \cite{Taylor} the pressure along the plate is given according to the Bernoulli's law by the relations
\be
\label{eq.jet_pressure_exact}
\begin{cases} 
	x = &\frac{1}{2\pi} \left(\left(1+\cos \theta \right) \log(1+q) - \left(1-\cos\theta\right)\log\left(1-q\right)\right) + 
	\frac{1}{\pi} \sin \theta \sin^{-1}\!q + \text{const}, 
	\\[3pt]
	p = &\frac{1}{2}\rho_0|\u|^2(1-s^2) \\
\end{cases}
\ee 
with 
\be
s = \frac{-(1-q\cos \theta) + \sqrt{1-q^2}\sin \theta}{q-\cos\theta} \quad \text{and} \quad 0 < q < 1.
\ee   
In Figure \ref{fig.JetPressure} the pressure obtained with our method along the plate at time $t=6.5$ is compared against the exact solution. 

Finally, to reach the final time $t=6.5$ the scheme requires $27722$ time steps, so the total number of processed volumes is $\texttt{nVol} = 500 \times 350 \times 27722 = 4.85\cdot 10^9$; by employing the same GPU as in the previous tests, 
the algorithm is executed in $424$ seconds, with an average of $1.1\cdot10^7$ volumes processed per second. 

\subsection{Flow over a sharp-crested weir}
\label{sec.FlowOver}

In this section we study the flow over a sharped-crested weir with a setting similar to the one proposed in \cite{ferrari2010sph}. In particular from now gravity is switched on using $g=9.81$.
We take a rectangular computational domain $\Omega=[-7.5,7.5]\times [0,2.1]$ and we model the sharp-crested weir through a strip  situated at $x=0$, of height $h=0.7$ and a thickness of one cell, whose boundary is treated as reflective slip wall boundary condition.
At time $t=0$ the liquid domain is $\Omega_\ell=[-7.5,0]\times[0,1.5]$ and inside $\rho$ is chosen in order to be in equilibrium
\be
\label{eq.rho_in_liquid}
\rho(x,y) = \rho_0 \exp \left(-\frac{g\rho_0}{k_0} \left(y-y_0\left(x\right)\right) \right),
\ee 
where $y_0(x)$ is the initial vertical position of the free surface as a function of $x$, and the initial velocity is $\u=\0$ everywhere. 
Based on the reference velocity $|\u|_{\max} = 2 \sqrt{gH} = 2 \sqrt{1.5g} = 7.67$, to obtain a Mach number of $M=0.3$ the constant $k_0$ is chosen to be $k_0=6.54\cdot 10^5$.

We have covered the computational domain with a Cartesian mesh of $800\times250$ elements and to reach the final time $t=1.0$ the scheme requires $16306$ time steps, so the total number of processed volumes is $\texttt{nVol} = 800 \times 250 \times 16306 = 3.26\cdot 10^9$; on the GPU the algorithm is executed in $233$ seconds, with an average of $1.4\cdot10^7$ volumes processed per second.

\begin{figure}
	\begin{center}
		\includegraphics[width = 0.9\linewidth]{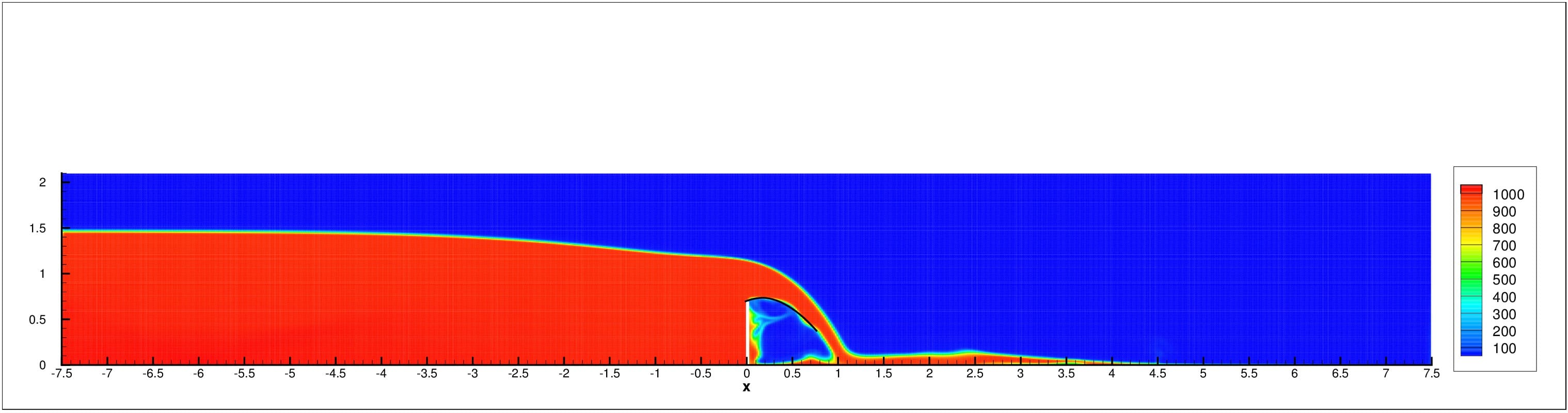}  
		\caption{Flow over a sharp-crested weir at time $t=1$. We report the density profile of $\alpha\rho$ and the experimental reference solution for the over topping profile \eqref{eq.OverToppingProfile}  (black line).}
		\label{fig.FlowOver_big}
	\end{center} 
\end{figure}

\begin{figure}
	\begin{center}
		\includegraphics[width =0.6\linewidth]{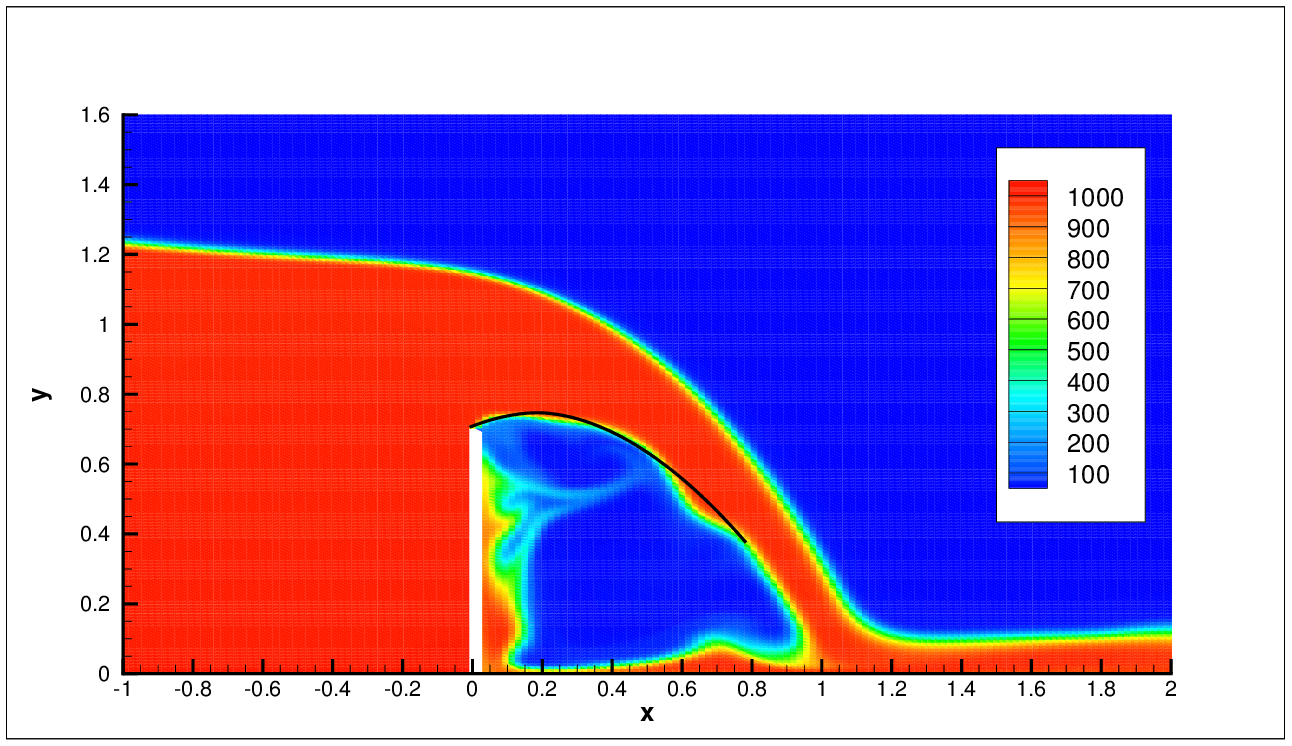}  
		\caption{Zoom into the weir flow at time $t=1$. Density contour of $\alpha\rho$ and experimental reference solution (black line).}
		\label{fig.FlowOver_small}
	\end{center} 
\end{figure}

The final time $t=1$ has been chosen because a rather stationary over topping profile is obtained that can be compared with the experimental reference solution of the lower streamline of the profile found in \cite{Scimemi} and reported in \cite{ferrari2010sph} as
\be
\label{eq.OverToppingProfile}
y(x) = y_m - 0.47\overline{h}_0\left( \frac{x-x_m}{\overline{h}_0} \right)^{1.85},
\ee 
where $x_m$ and $y_m$ is the location of the maximum of the curve and $\overline{h}_0$ is the vertical distance between the crest of the weir and the free surface at the location of the weir.
The obtained numerical results are depicted in Figure \ref{fig.FlowOver_big} and a zoom on the over topping profile compared against the reference solution \eqref{eq.OverToppingProfile} is reported in Figure \ref{fig.FlowOver_small}. In particular we observe $x_m=0.18$, $y_m=0.78$ and $\overline{h}_0=0.46$ and we note a good agreement between the numerical results obtained with our second order well balanced two-phase model and the experiments.

\subsection{Dambreak}
\label{sec.dambreak}

A very typical application for shallow water-type models is the so-called dambreak. It consists of the sudden collapse/removal of a vertical wall that separates two different piecewise constant states of water from each other. Since in the initial stages of dambreak flows, the classical shallow water assumption of small vertical velocities and accelerations does not hold, it is of interest to apply our more complete model to this well-studied phenomenon.

In this section, we will compare our second order well balanced scheme for the three-equation two-phase model
with the third order P$_0$P$_2$ finite volume scheme presented in \cite{dumbser2011simple} for the same system of equations.
Moreover we will analyze the differences between our scheme and exact or numerical solutions of the shallow water equations.
In particular, we will discuss the different behavior at small times (shortly after the dambreak) 
confirmed also by the results obtained with the new SPH method of Ferrari proposed in \cite{SPH3D} \cite{FerrariPhD} and whose results for the following four tests are shown in \cite{dumbser2011simple} obtained using 224,282 SPH particles with a characteristic particle distance of h = 0.03. 
In all the four test problems shown below, we can note an excellent agreement between
the new two-phase flow model and the 2D SPH simulations, whereas there are significant discrepancies between the full 2D models and the 1D shallow water model at short times. At large times, however, both weakly compressible 2D models agree very well with the shallow water theory.
Finally, our well balanced second order scheme on Cartesian meshes shows an excellent agreement with the results obtained in \cite{dumbser2011simple} through a third order method on unstructured meshes using  more or less the same characteristic mesh spacing. 


\subsubsection{Dry bed}
\label{ssec.drybedWithoutStep}

\begin{figure*}
	\centering
	\includegraphics[width=0.49\linewidth]{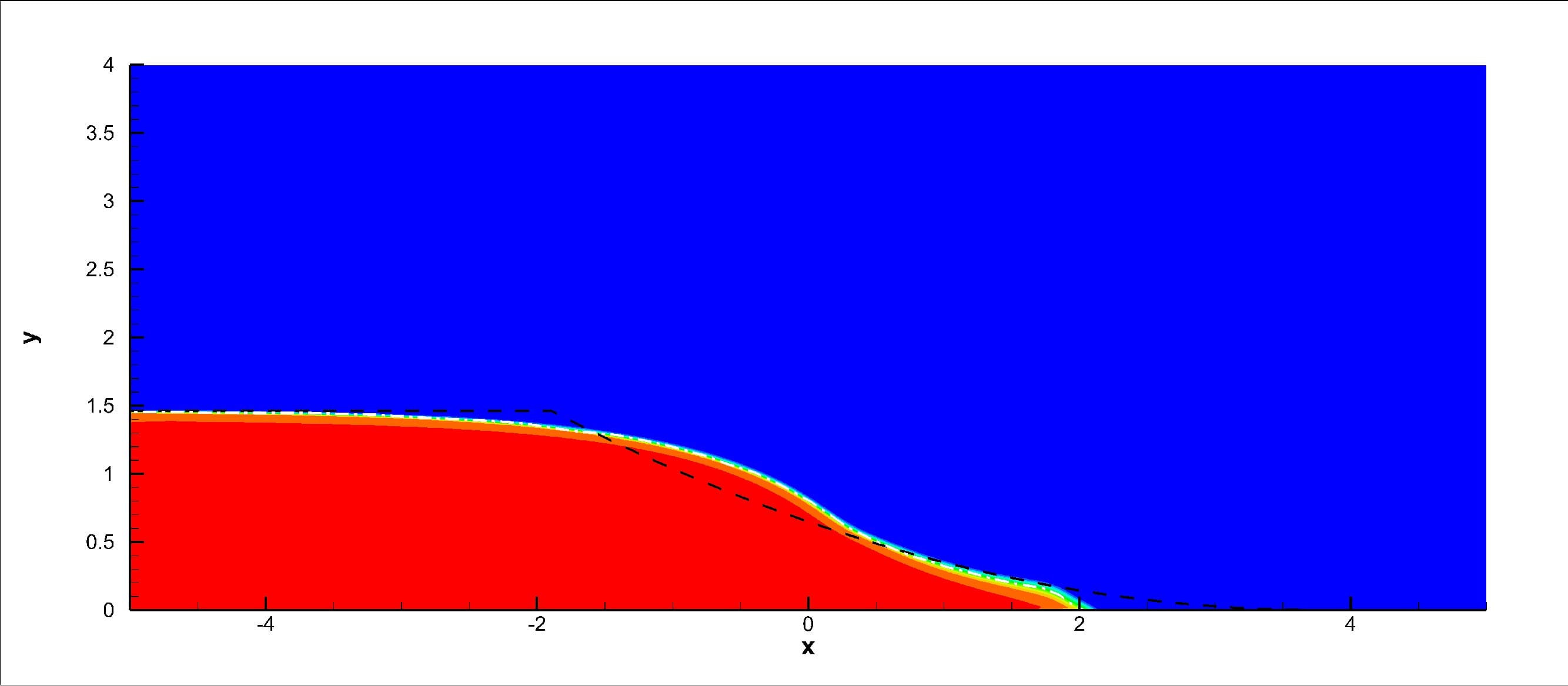} 
	\includegraphics[width=0.49\linewidth]{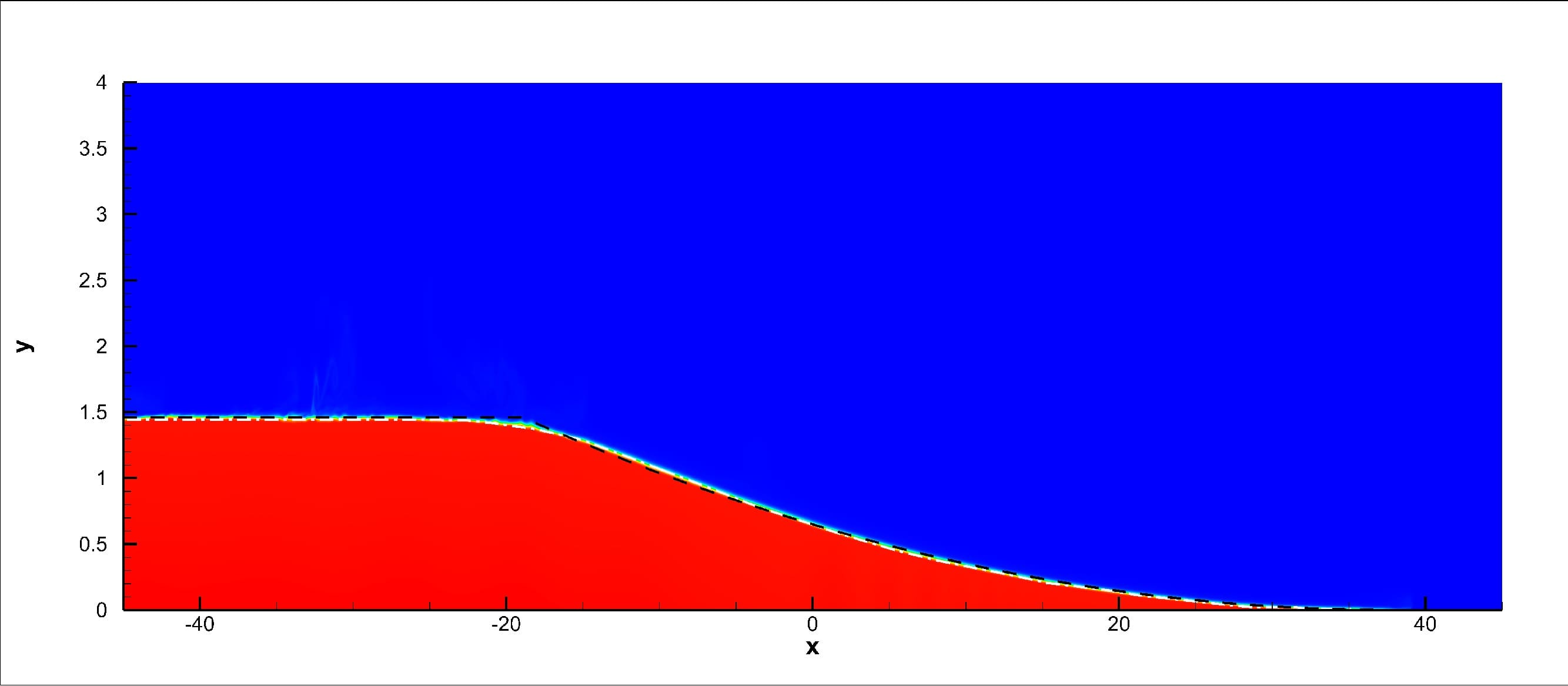}\\ \,
	\includegraphics[width=0.49\linewidth]{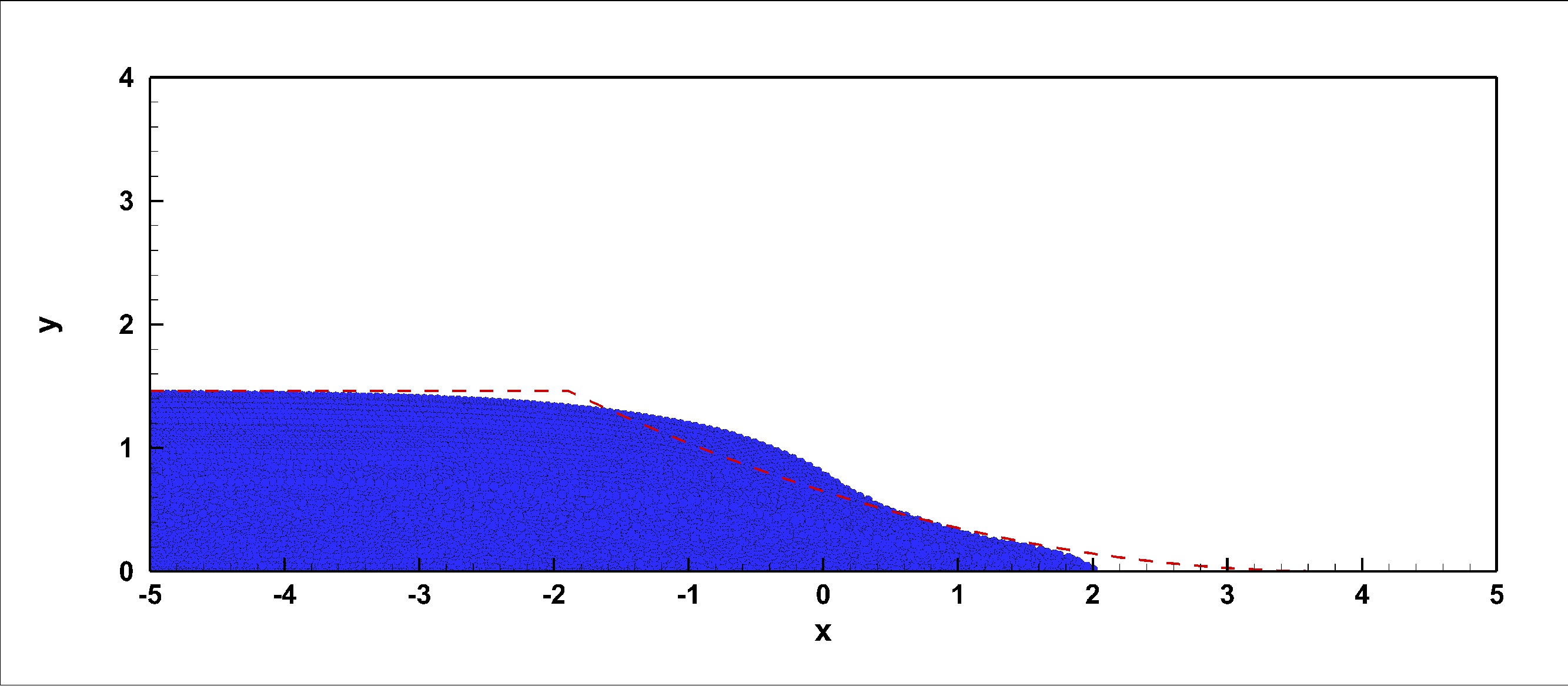}
	\includegraphics[width=0.49\linewidth]{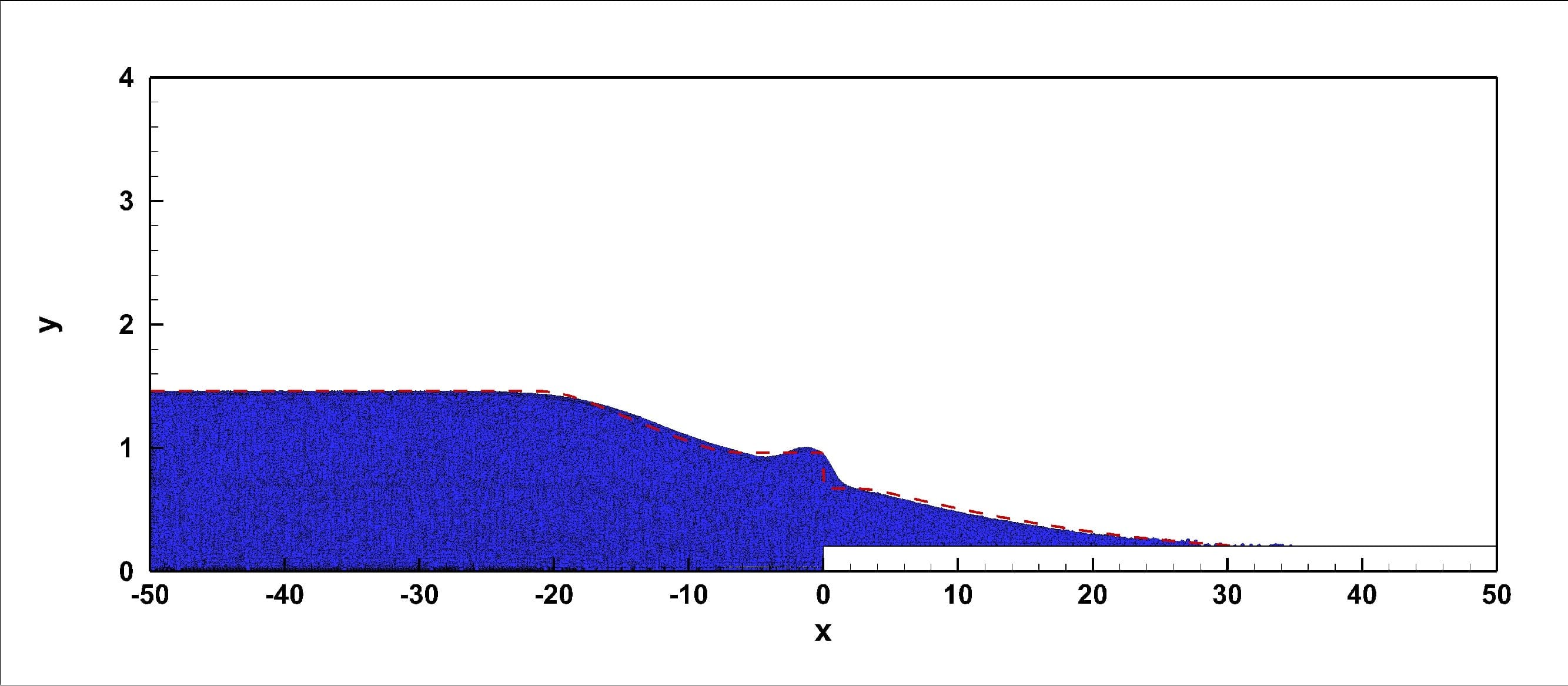} \\ \,
	\caption[Dambreak with dry bed.]{Dambreak with dry bed at time $t=0.5$ (left) and $t=5.0$ (right). Top row: we compare our numerical results (density contours) with the interface water-air taken from \cite{dumbser2011simple} (white dashed line), and the solution obtained using the shallow water model (black dashed line).  Bottom row: we report the results obtained with the SPH method of \cite{SPH3D} (blue circles) and the shallow water model (red dashed line). }
	\label{fig.Dambreak_withoutStep_dry}
\end{figure*}

First, we consider the classical dambreak over a dry bed. 
We take a computational domain $\Omega=[-50,50]\times[0,4]$ covered with a Cartesian mesh of $4000\times400$ elements.
At the initial time $t=0$ the liquid is contained in $\Omega_\ell = [-50,0]\times[0,1.4618]$.
The constants that characterize the problem are chosen to be $k_0~=~6.37\cdot10^5$, $\rho_0=1000$ and $\gamma=1$ (in this way the Mach number is $M=0.3$ with a maximum expected velocity $|\mbf{u}|=7.57$).
The initial velocity is $\mbf{u}=\0$ everywhere, and the initial density distributions in $\Omega_\ell$ is given by \eqref{eq.rho_in_liquid}.
In the air, as in all the other tests, we impose $p(\x,0) = 0$, and $\rho(\x,0)=\rho_0=1000$.
The boundary conditions are reflective wall on the left, bottom and right border of the computational domain and transmissive boundaries on the top.

The results presented in Figure \ref{fig.Dambreak_withoutStep_dry} clearly show the superiority of the two-phase model, which is more complete with respect to the classical shallow water assumptions in particular for modeling the first instances after the dambreak. Moreover the results obtained with our second order well balanced scheme are in perfect agreement with the reference  solutions taken from \cite{dumbser2011simple} and \cite{SPH3D}.
After longer times all the presented methods and models perform equally well. These results show the reliability of our scheme  both for short and long times after the dambreak. 

To reach the final simulation time $t=5$ the scheme requires $84016$ time steps, hence the total number of processed volumes is  $\texttt{nVol} = 4000 \times 400 \times 84016 = 1.34\cdot 10^{11}$; by employing the NVIDIA GeForce Titan Black GPU and 
rectangular blocks of $32\times2$ threads, the algorithm is executed in $7808$ seconds, with an average of $1.7\cdot10^7$ 
volumes processed per second. 

\subsubsection{Wet bed}
\label{ssec.wetBedWithoutStep}

\begin{figure*}
	\centering
	\includegraphics[width=0.49\linewidth]{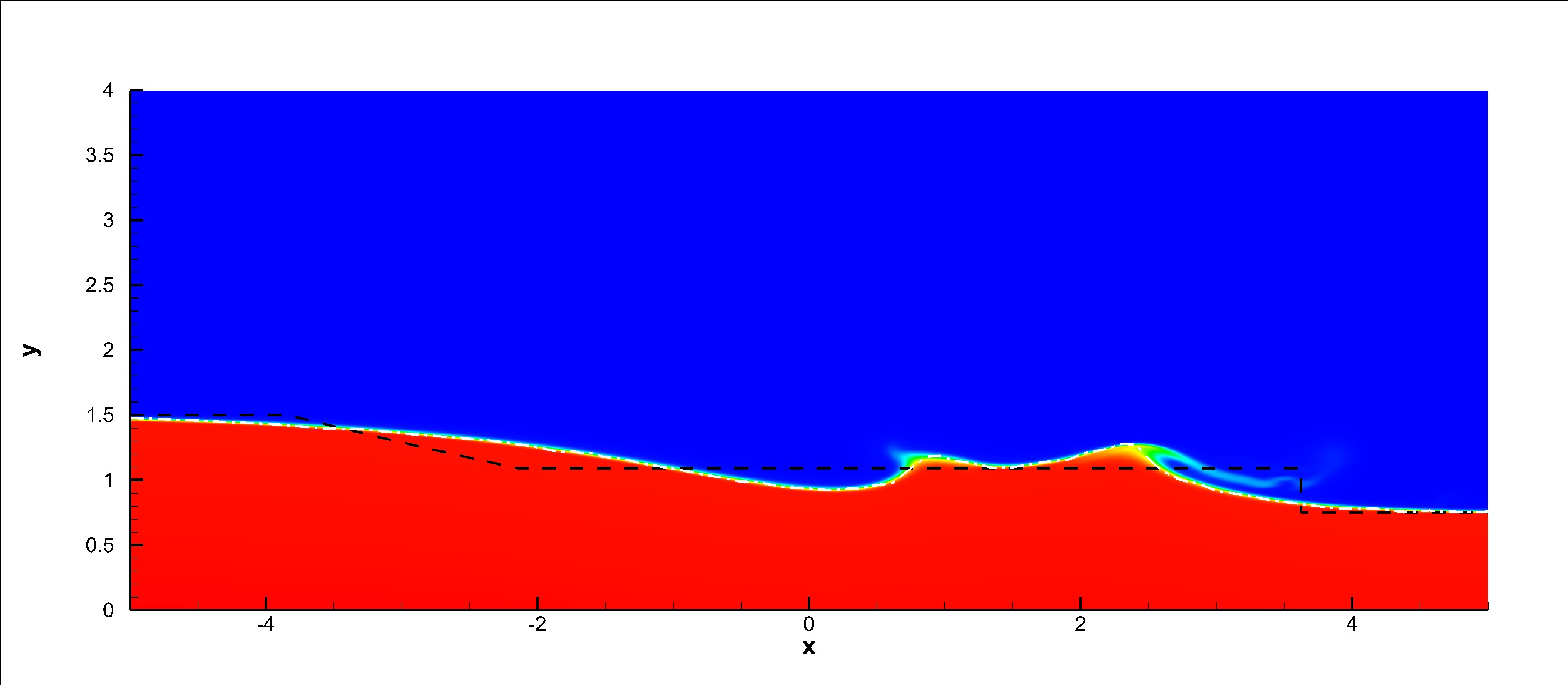}
	\includegraphics[width=0.49\linewidth]{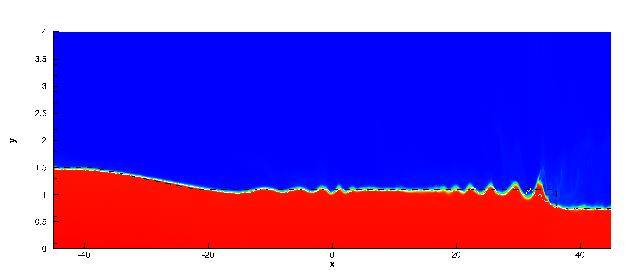} \\ \,
	\includegraphics[width=0.49\linewidth]{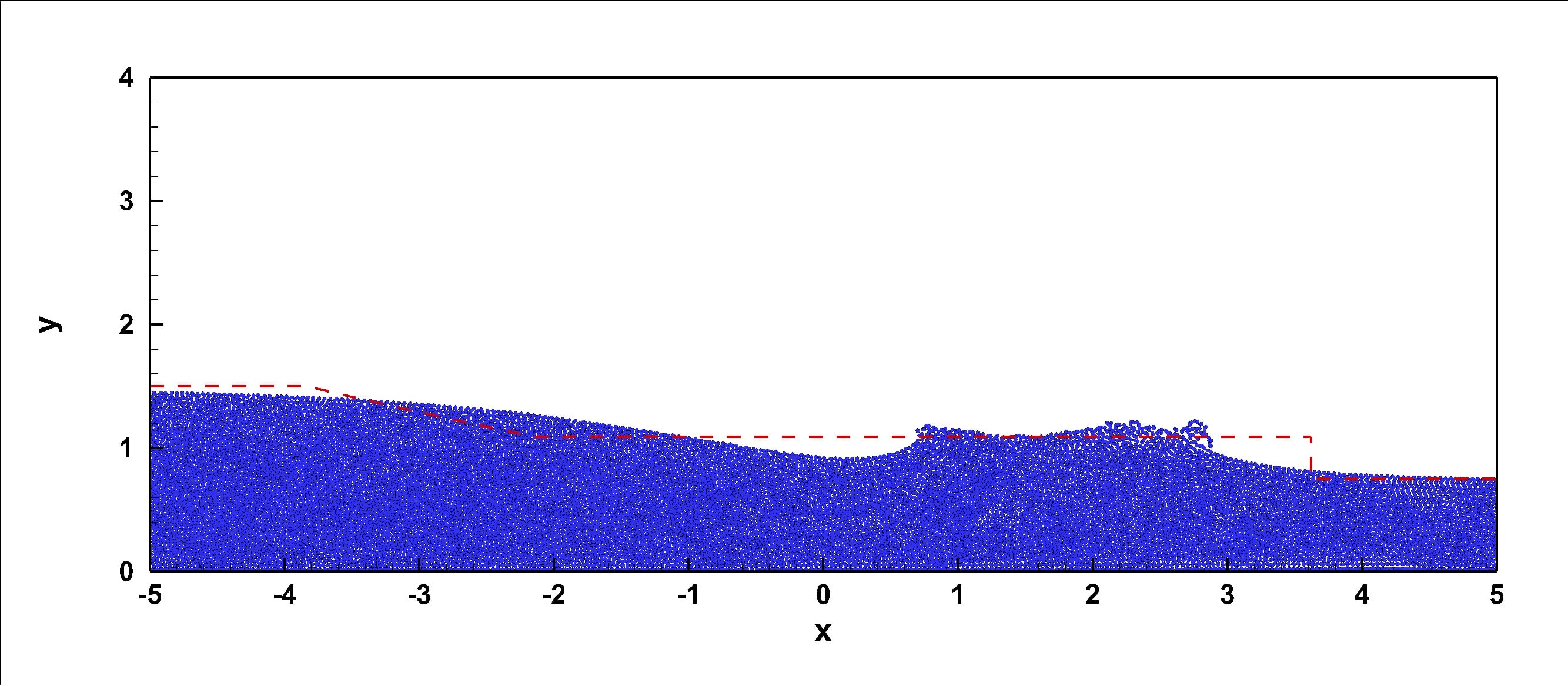}
	\includegraphics[width=0.49\linewidth]{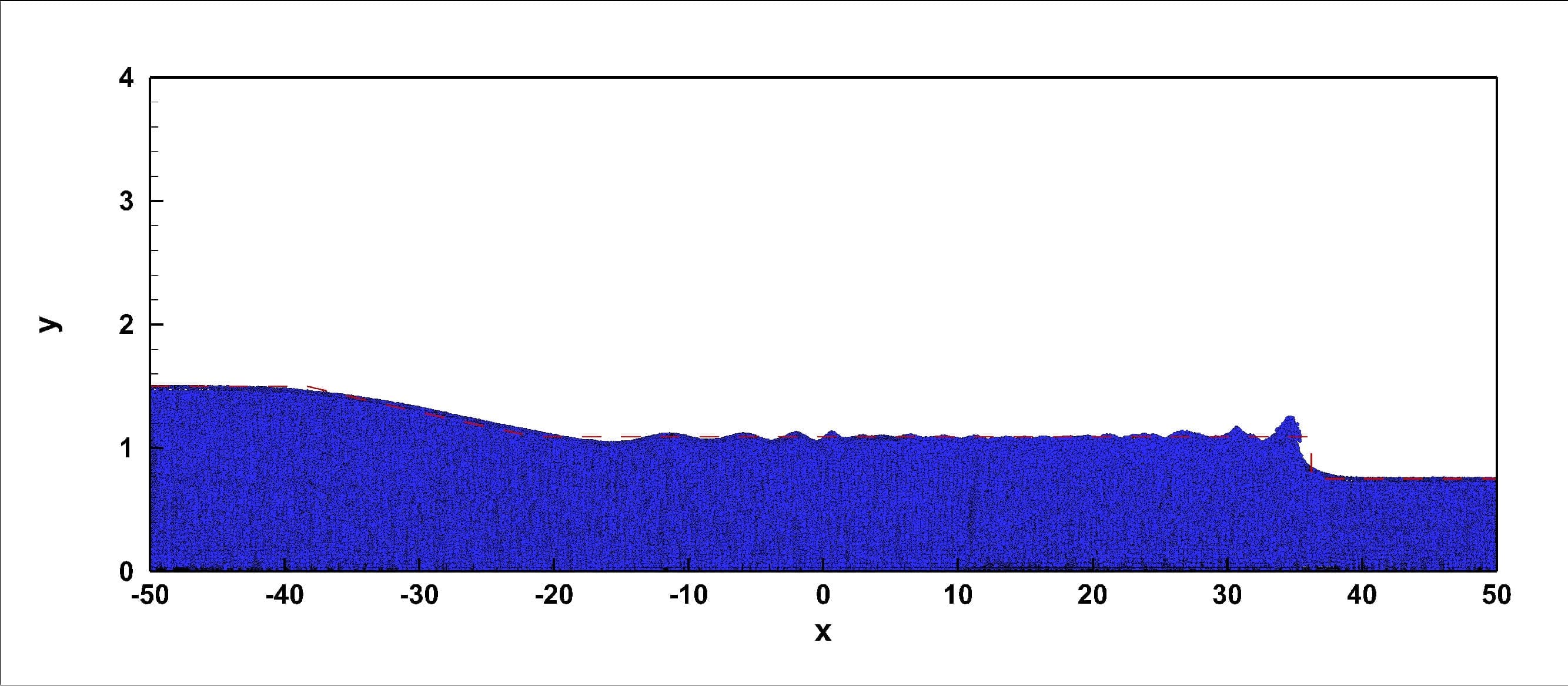}\\ \,
	\caption[Dambreak with wet bed.]{Dambreak with wet bed at time $t=1.0$ (left) and $t=10.0$ (right). Top row: we compare our numerical results (density contours) with the interface water-air taken from \cite{dumbser2011simple} (white dashed line), and the solution obtained using the shallow water model (black dashed line).  Bottom row: we report the results obtained with the SPH method of \cite{SPH3D} (blue circles) and the shallow water model (red dashed line). }
	\label{fig.Dambreak_withoutStep_wet}
\end{figure*}

We now consider the dambreak into a wet bed. 
We take a computational domain $\Omega=[-50,50]\times[0,4]$ covered with a Cartesian mesh of $4000\times400$ elements.
At the initial time $t=0$ the liquid is contained in $\Omega_\ell = [-50,0]\times[0,1.5] \cup [0,50] \times [0, 0.75]$.
As in the previous case we choose a Mach number of $M=0.3$ based on the maximum expected velocity of $|\mbf{u}|=7.67$, leading to $k_0~=~6.54\cdot10^5$, $\rho_0=1000$ and $\gamma=1$. All the other conditions are taken as in the previous test.

Experimental observations \cite{Janosi} show that breaking of waves can occur in the wet bed case. However, the wave breaking is the smaller the higher the water level on the right side of the dam. For the present test case, the level of the right water layer is rather high so that only very limited wave breakings occur during the simulation, as already observed in \cite{dumbser2011simple}.
For both the short time $t = 1$ and the final time $t = 10$ we do not observe any wave breaking phenomena in Figure \ref{fig.Dambreak_withoutStep_wet} but only the formation of small-scale free surface waves is visible in the flat region between the bore and the rarefaction wave. For comparison, we also plot the exact solution of the shallow water equations, which is in very good agreement with the results of the two-phase model for large times, in particular the water depth behind the bore (post-shock value). The disagreement for short times ($t = 1.0)$ is due to the fact that the shallow water model neglects the vertical accelerations that are very important in the early stages of the dambreak flow.

\subsubsection{Dry bed with step}

\begin{figure*}
	\centering
	\includegraphics[width=0.49\linewidth]{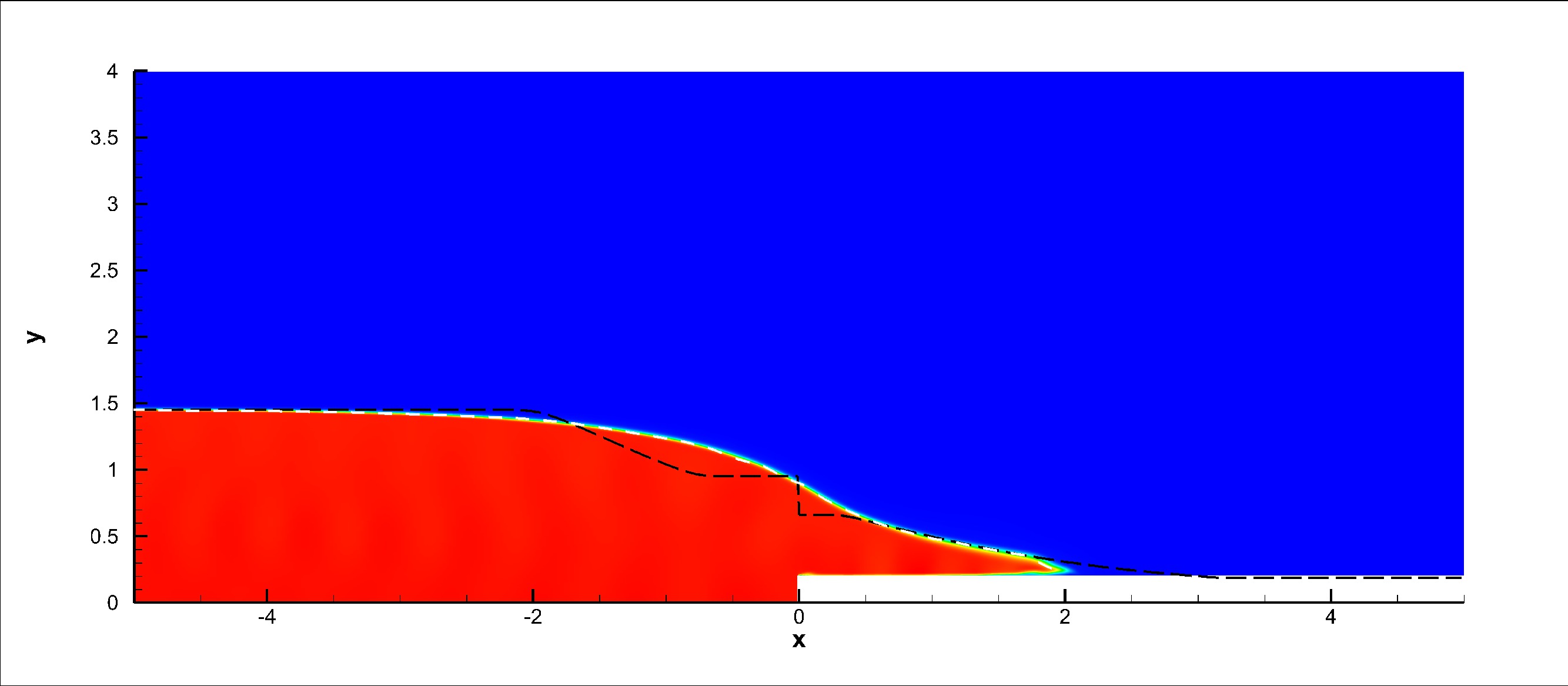}
	\includegraphics[width=0.49\linewidth]{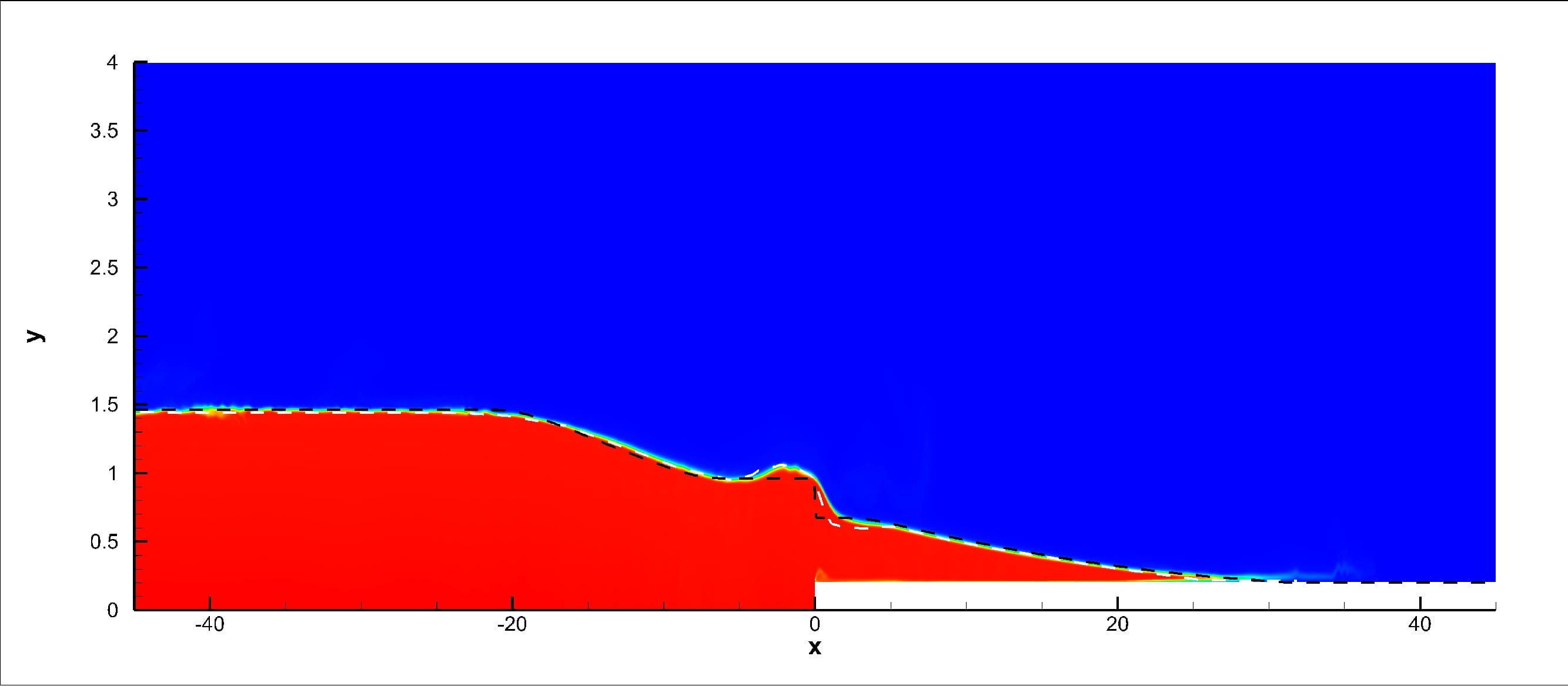} \\ \,
	\includegraphics[width=0.49\linewidth]{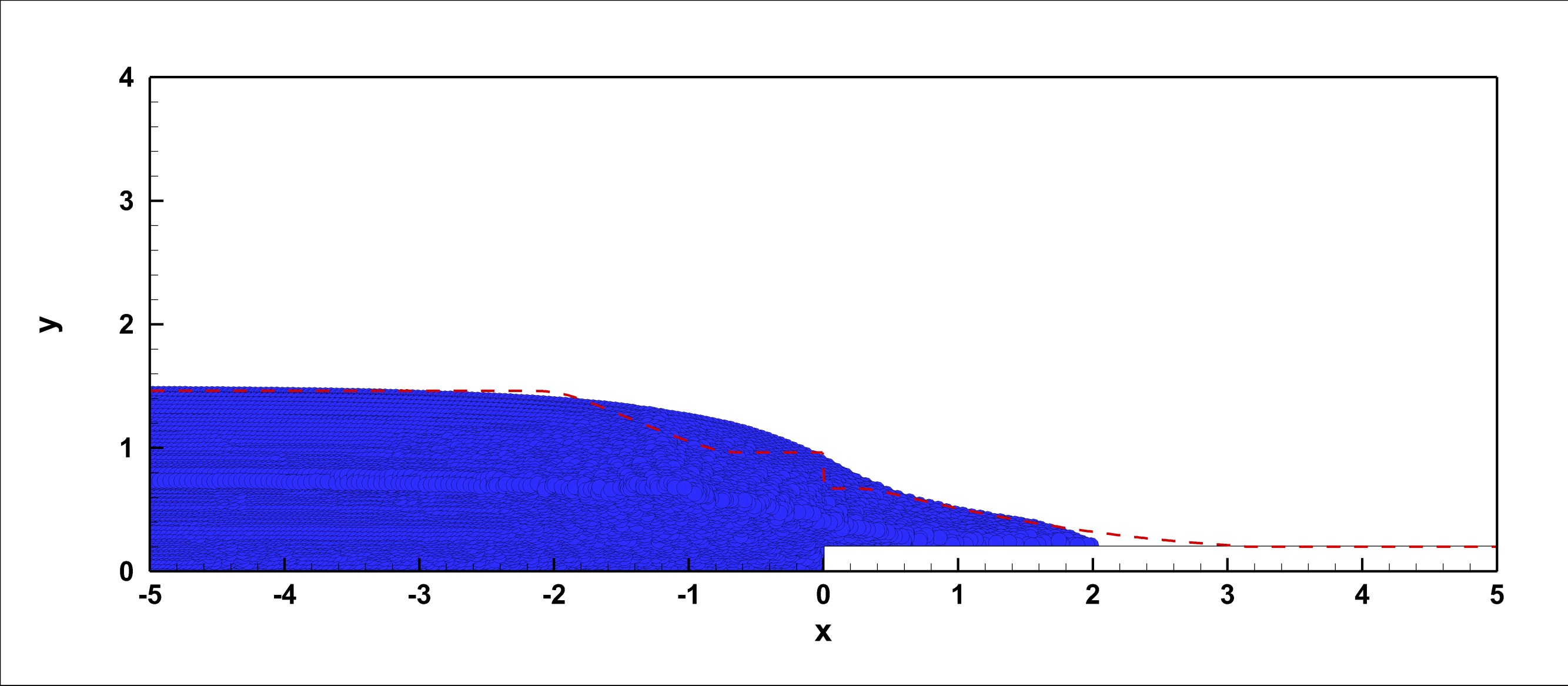}
	\includegraphics[width=0.49\linewidth]{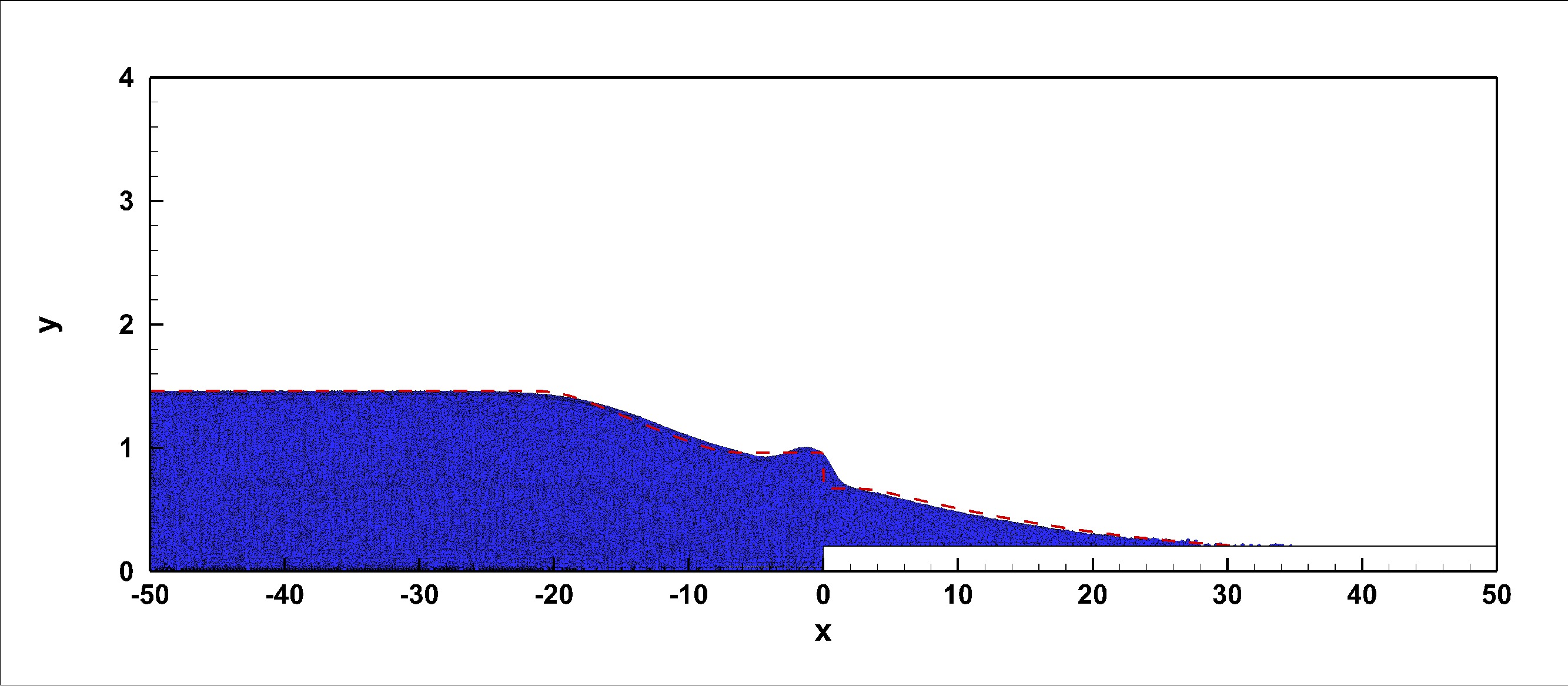} \\ \,
	\caption[Dambreak with step and dry bed.]{Dambreak with step and dry bed at time $t=0.5$ (left) and $t=5.0$ (right). Top row: we compare our numerical results (density contours) with the interface water-air taken from \cite{dumbser2011simple} (white dashed line), and the numerical solution obtained using the shallow water model (black dashed line).  Bottom row: we report the results obtained with the SPH method of \cite{SPH3D} (blue circles) and the shallow water model (red dashed line). }
	\label{fig.Dambreak_withStep_dry}
\end{figure*}

The next test problem is initialized as the one of Section \ref{ssec.drybedWithoutStep}, with the introduction of a \textit{step} of height $h=0.2$ between $x=0$ and $x=50$ modeled through reflective boundary conditions. Hence the computational domain becomes $\Omega= \left([-50,50] \times[0,4]\right) \backslash \left([0,50]\times[0,0.2]\right)$ covered by a Cartesian mesh of $\left(4000\times400\right) \backslash \left( 2000\times20\right)$, and the liquid occupies $\Omega_\ell = [-50,0] \times[0,0.4618]$.

The computational results obtained at the early time $t = 0.5$ and
at the later time $t = 5$ are depicted in Figure~\ref{fig.Dambreak_withStep_dry}, where they are compared
with a numerical solution computed with a third order path-conservative
WENO finite volume scheme \cite{ADERNC}\cite{USFORCE2}\cite{OsherNC} applied to the
$1$D shallow water equations with variable bottom. As in the
previous case, there are significant discrepancies between the shallow
water model and the two-phase model at early times,
but for large times the agreement is excellent, despite the initial
deviation between the two models.
Moreover our results are in perfect agreement with the ones obtained in \cite{dumbser2011simple} with a third order scheme and with the SPH method of \cite{SPH3D}.

\subsubsection{Wet bed with step}

\begin{figure*}
	\centering
	\includegraphics[width=0.49\linewidth]{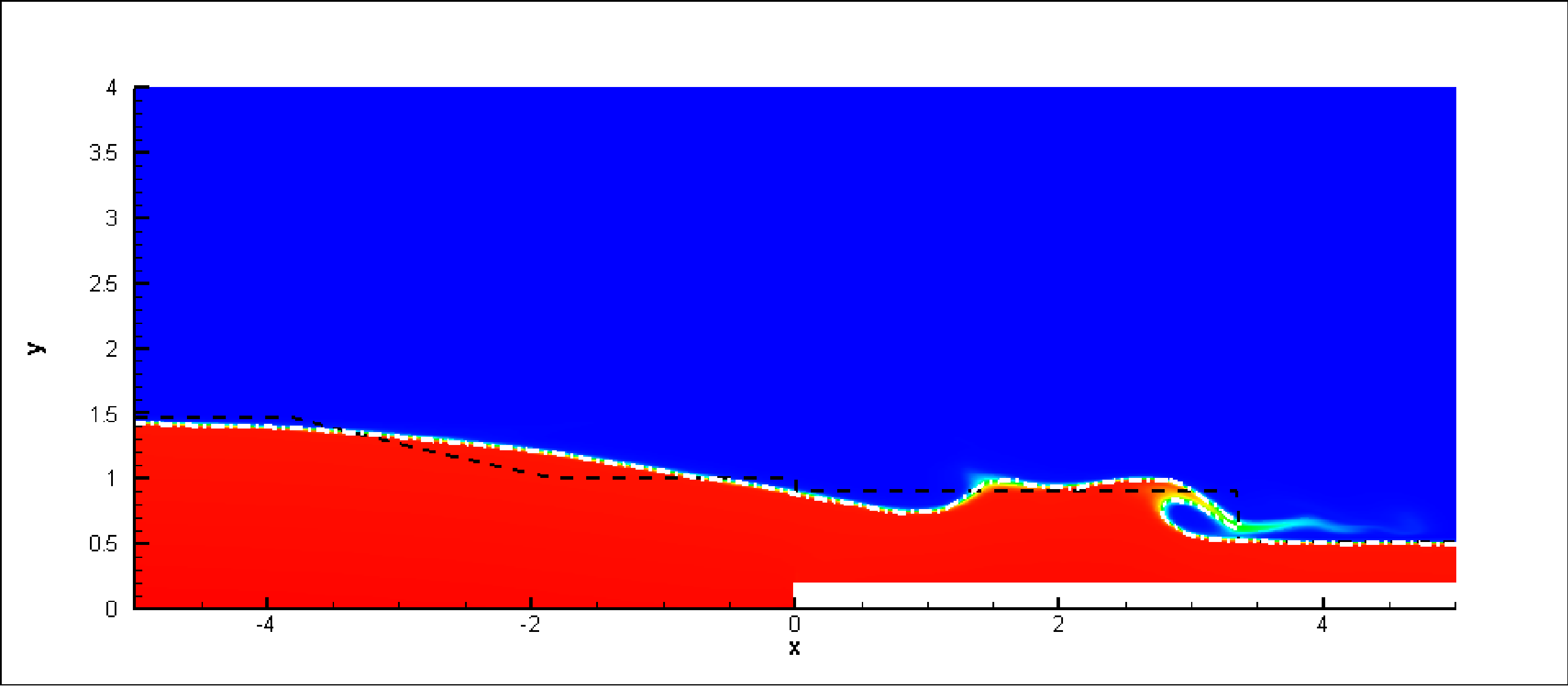}
	\includegraphics[width=0.49\linewidth]{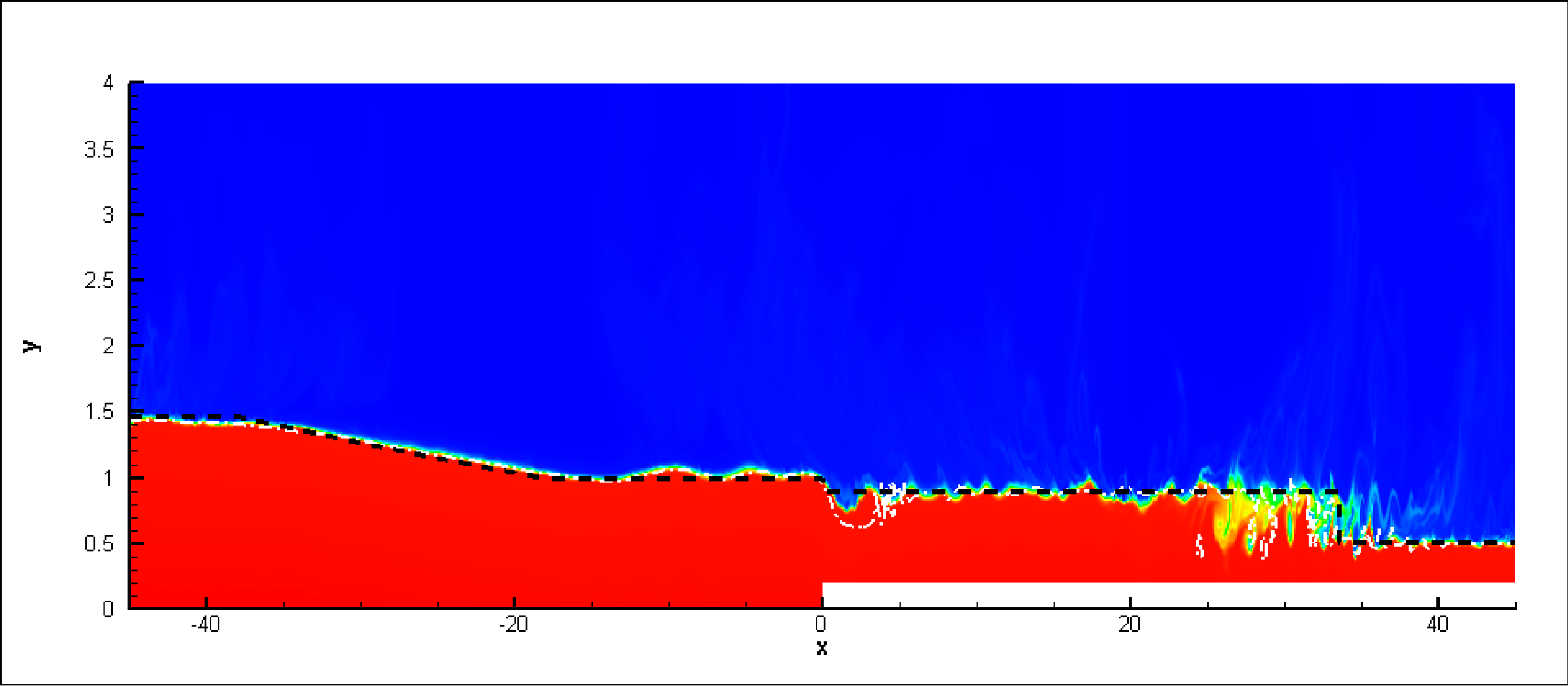} \\ \,
	\includegraphics[width=0.49\linewidth]{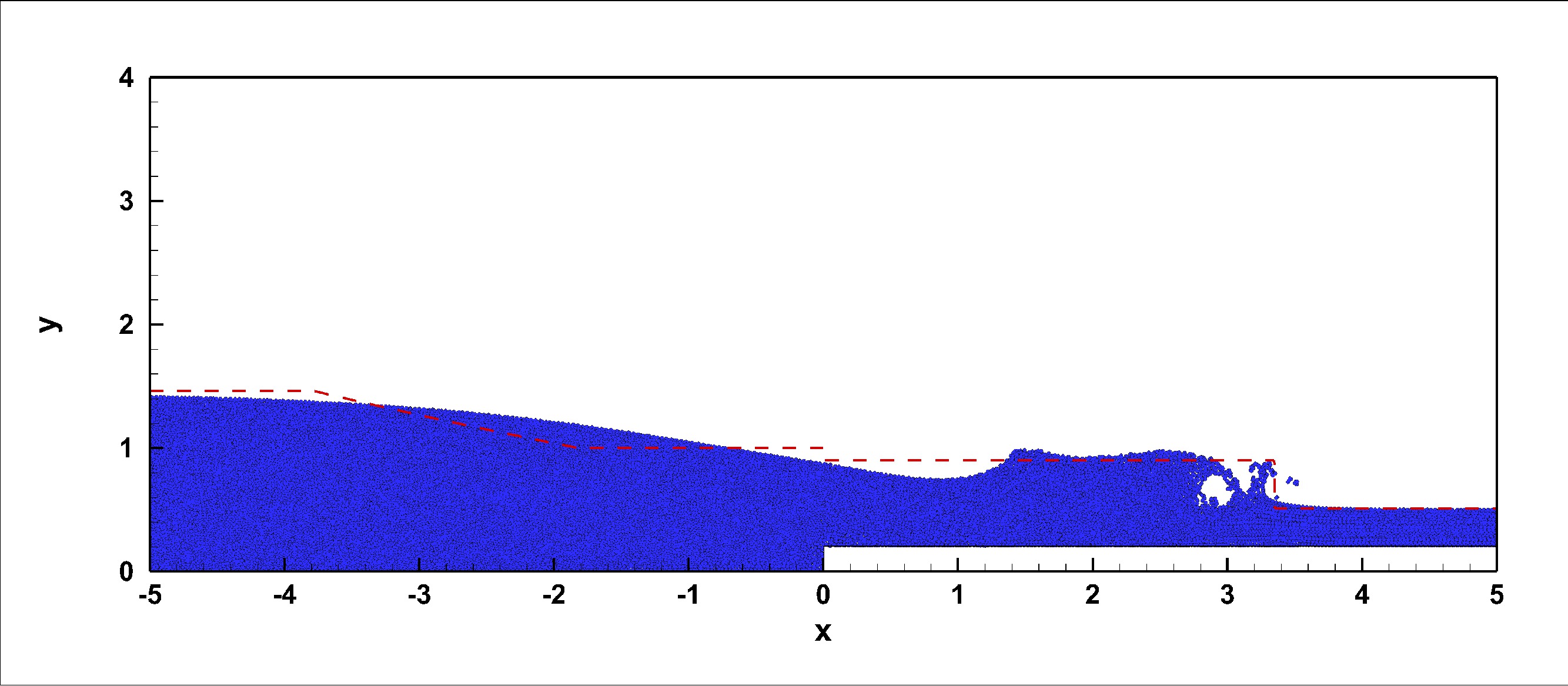}
	\includegraphics[width=0.49\linewidth]{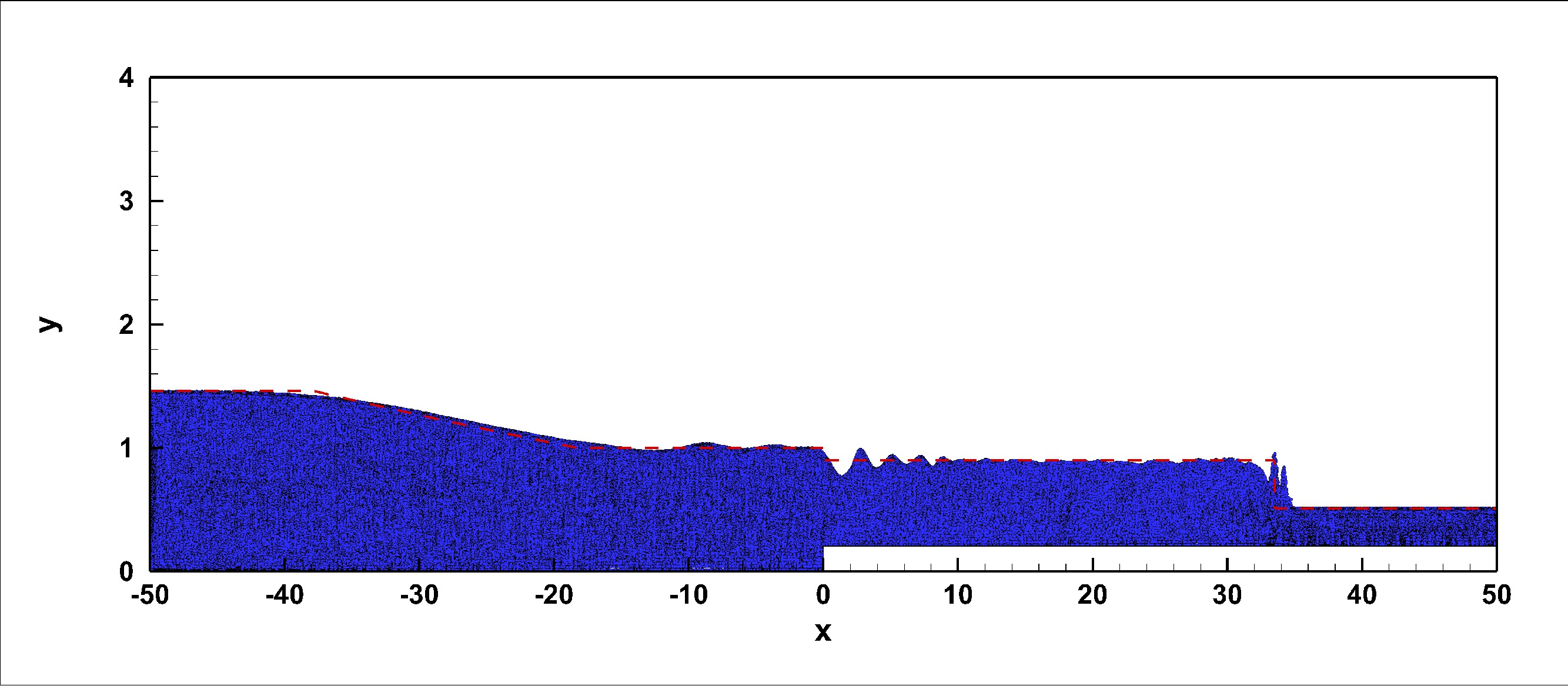}\\ \, 
	\caption[Dambreak with step and wet bed.]{Dambreak with step and wet bed at time $t=1.0$ (left) and $t=10.0$ (right). Top row: we compare our numerical results (density contours) with the interface water-air taken from \cite{dumbser2011simple} (white dashed line), and the solution obtained using the shallow water model (black dashed line).  Bottom row: we report the results obtained with the SPH method of \cite{SPH3D} (blue circles) and the shallow water model (red dashed line). }
	\label{fig.Dambreak_withStep_wet}
\end{figure*}

\begin{figure}
	\centering
	\includegraphics[width=0.6\linewidth]{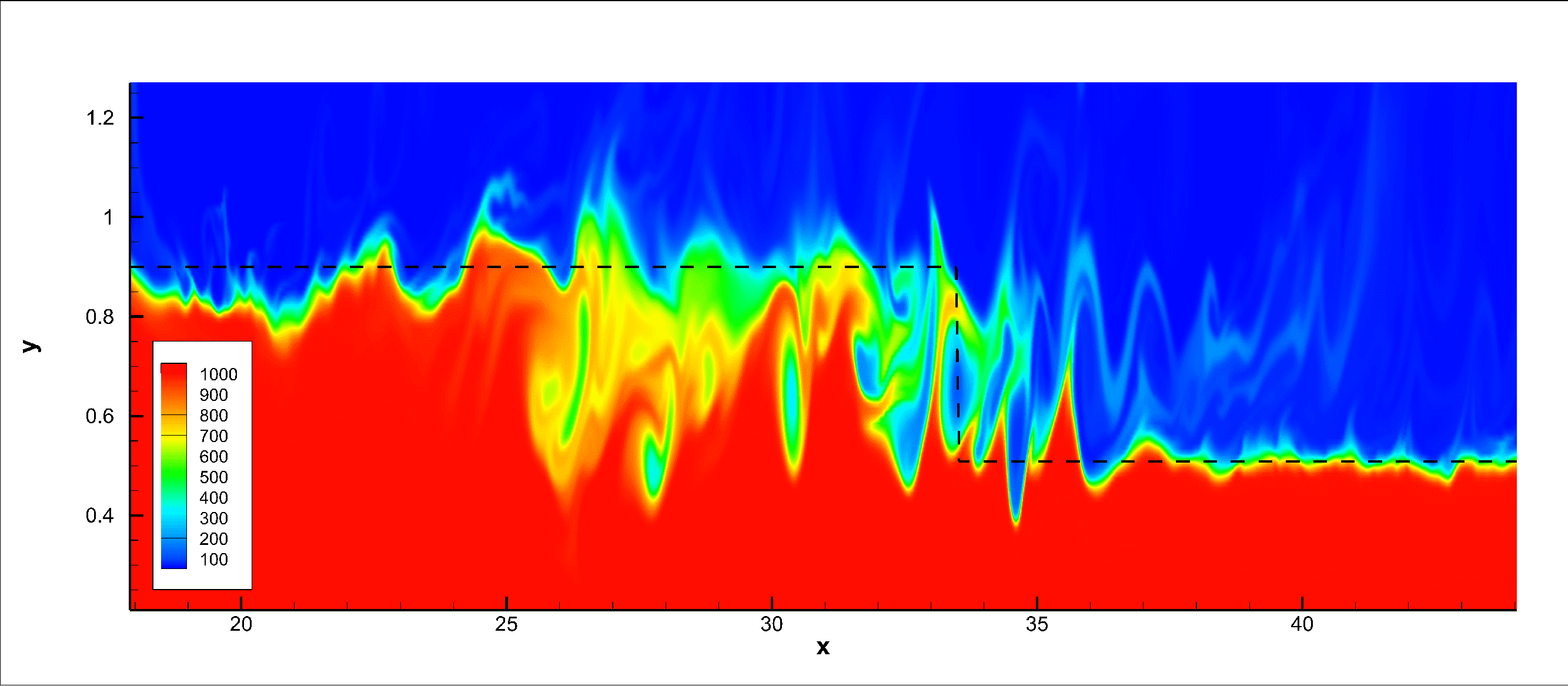}
	\caption{Dambreak with step and wet bed at time $t=10$, zoom into the complex free surface generated at time $t=10$.}
	\label{fig.Dambreak_withStep_wet_zoom_1}
\end{figure}

\begin{figure}
	\centering
	\includegraphics[width=0.6\linewidth]{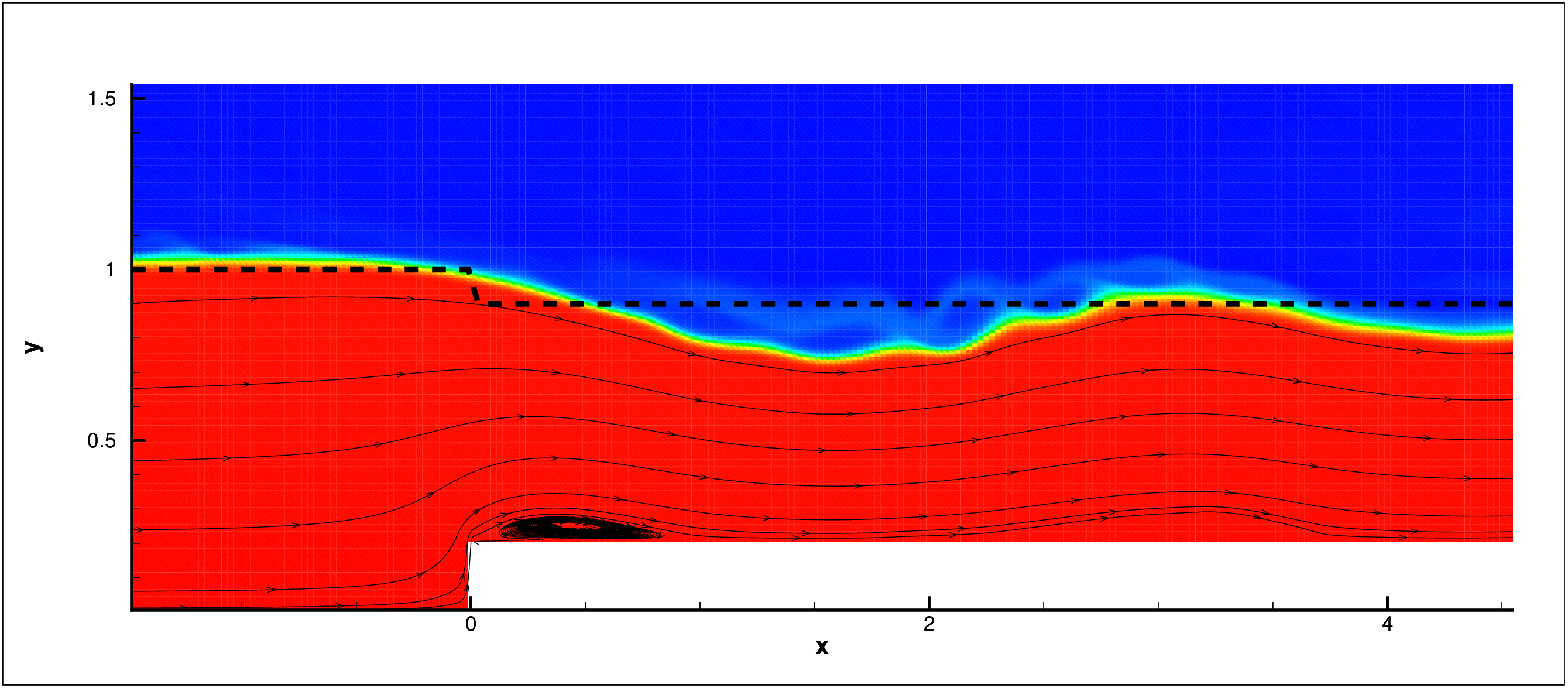}
	\caption[Dambreak with step and wet bed at time $t=10$, zoom into the recirculation zone in the vicinity of the bottom step.]{Dambreak with step and wet bed at time $t=10$, zoom into the recirculation zone in the vicinity of the bottom step. We depict our numerical results (density contour) and the exact solution of the shallow water equations with bottom step according to \cite{Bernetti}.}
	\label{fig.Dambreak_withStep_wet_zoom_2}
\end{figure}

For the last dambreak problem we consider the computational domain $\Omega= [-50,50] \times[0,4] - [0,50]\times[0,0.2]$ where the \textit{step} of height $h=0.2$ is modeled through reflective boundary conditions, and we cover it with a Cartesian  mesh of $\left(4000\times400\right) \backslash \left(2000\times20\right)$ elements.
The liquid is initially contained in $\Omega_\ell = [-50,0] \times[0,0.4618] \cup [0,50]\times[0.2,0.50873]$ and the other conditions are the same of Section \ref{ssec.wetBedWithoutStep}.
The results obtained with our second order well balanced scheme are depicted in Figure \ref{fig.Dambreak_withStep_wet} at the early time $t=1$ and at the later time $t=10$ and in both the cases we can observe the breaking waves phenomenon which, with respect to the previous test case with wet bed, appears due to the lower water level on the right of the dam. See also Figure \ref{fig.Dambreak_withStep_wet_zoom_1} for a zoom at $t=10$.
This considerable breaking of waves at the moving bore front is in agreement with the simulation obtained both in \cite{dumbser2011simple} and with the SPH method, while instead the shallow water model approximates this propagating wave front as a sharp moving discontinuity.
Moreover, these results are corroborated by experimental observations \cite{Janosi}\cite{IniDambreak} which clearly show that the dambreak into a wet region with rather low water level very quickly leads to breaking waves.

Whereas the shallow water model predicts only one rather small discontinuity
in the free surface directly over the bottom step, the two-phase
model predicts a much larger decrease of the free surface
profile and then, at about $x = 4$ (see Figure \ref{fig.Dambreak_withStep_wet_zoom_2}) there is again an increase in the free surface level to match the post-shock value. 
The physical explanation of this phenomenon can be found when tracing the streamlines in the vicinity of the bottom step where we observe a recirculation zone, which, for the rest of the flow, increases the apparent bottom step at $x = 0$.

\subsubsection{Comparison depth average velocity profiles}

\begin{figure*}
	\centering
	\includegraphics[width=0.45\linewidth]{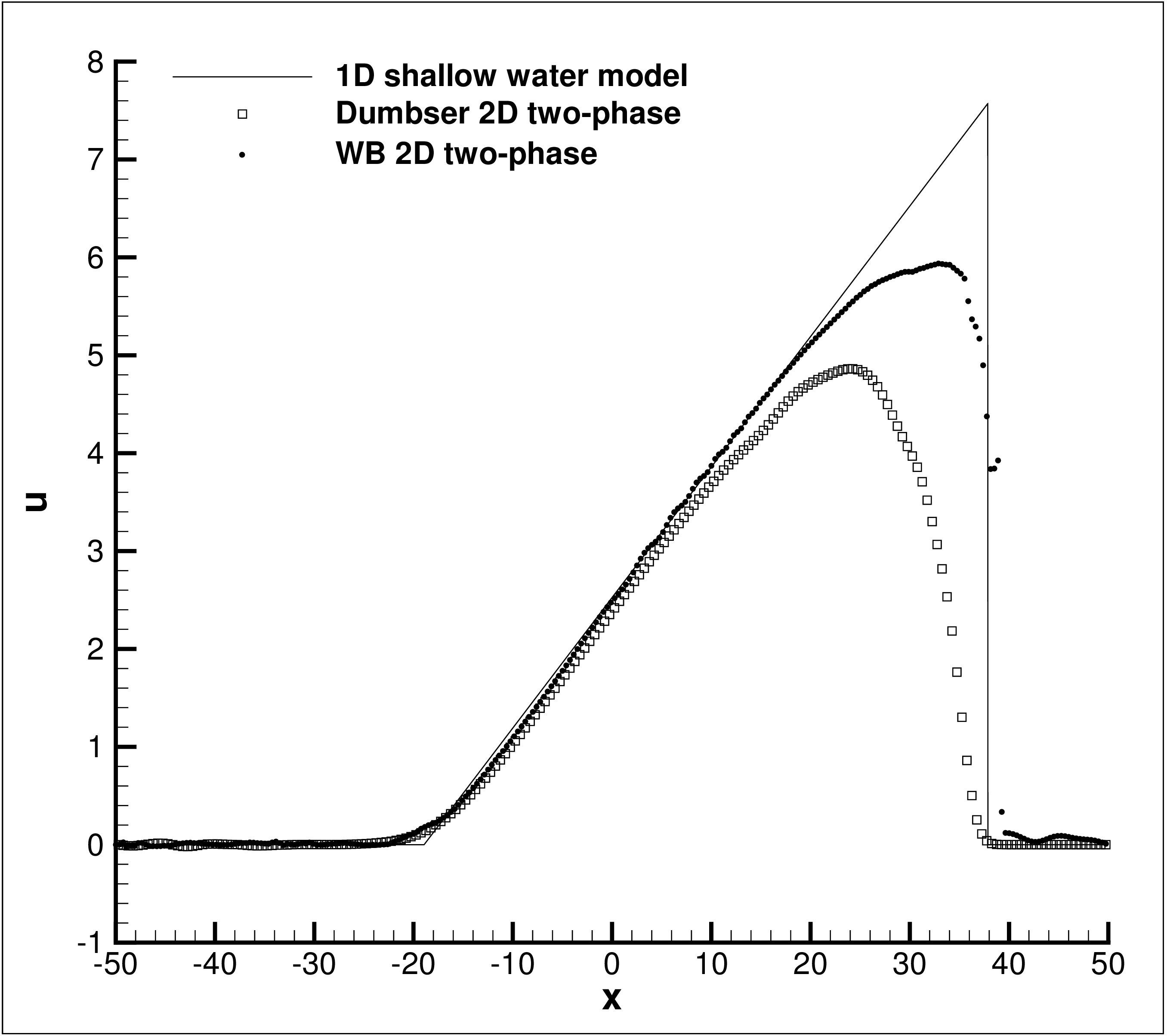}
	\ \includegraphics[width=0.45\linewidth]{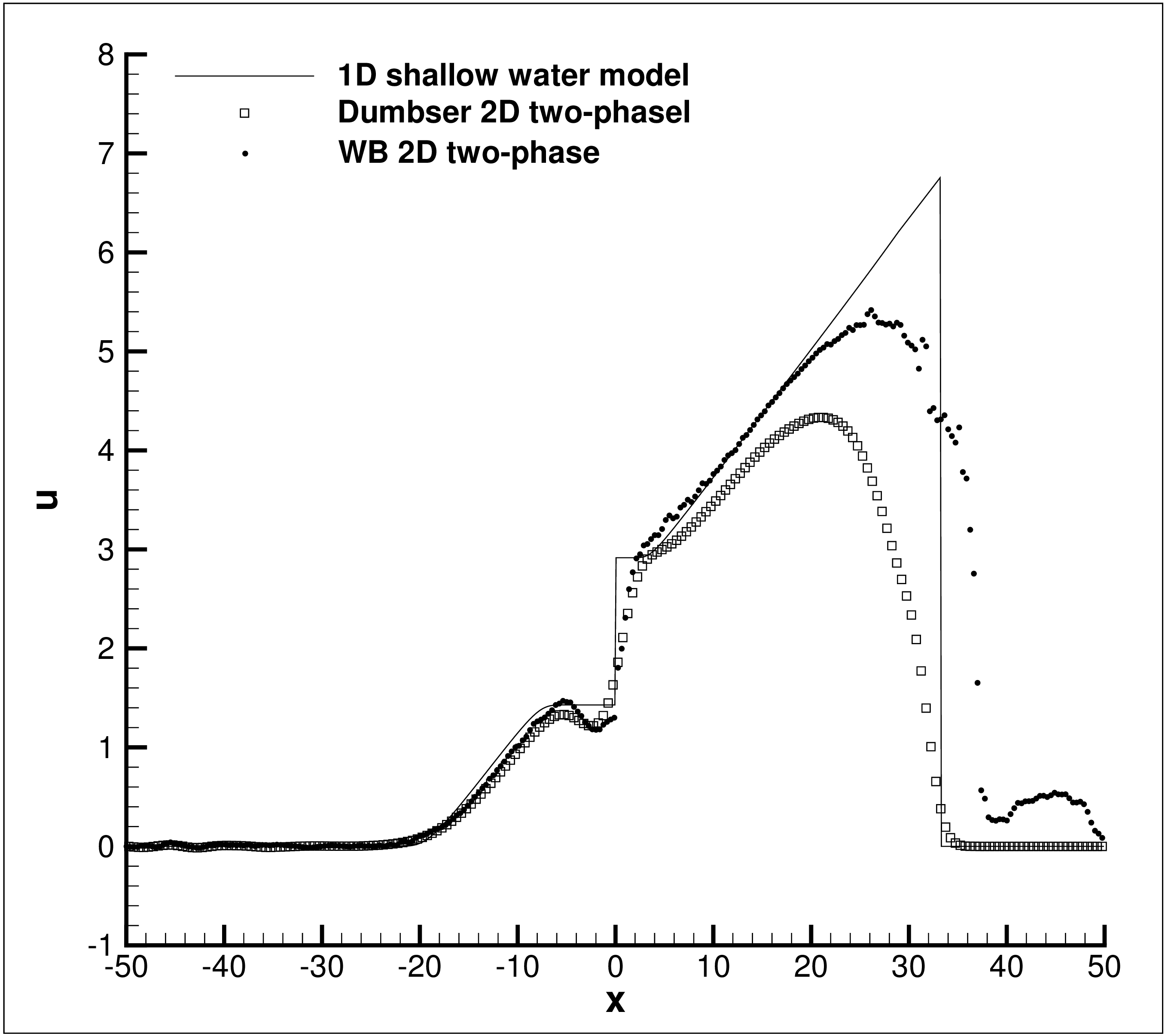}
	\caption{Comparison between the depth-averaged velocity profiles at the final times $t=5$ of the dambreak with dry bed without (left) with (right) bottom step.}
	\label{fig.depthAveregedVel}
\end{figure*}

In this section we compare the depth-averaged velocity profiles computed with our second order well balanced scheme and with the third order scheme proposed in \cite{dumbser2011simple} for the two-phase model against the results obtained with classical $1$D shallow water theory. In particular in $2$D the depth-averaged velocity profile is calculated as 
\be
\overline{u}(x) = \frac{\int_0^4 u(x,y)\alpha(x,y) dy}{\int_0^4 \alpha(x,y) dy}.
\ee 
The results are compared in Figure \ref{fig.depthAveregedVel} and our two-phase model, in agreement with laboratory experiments
 show a steep but not infinite gradient, in contrast to the shallow water model. We also
observe that the second order well-balanced method presented in this paper appears to be less dissipative than the 
third order non-wellbalanced WENO scheme used in \cite{dumbser2011simple}.

\subsection{Dambreak and impact against a vertical wall} 
\label{sec.Impacttest}

\begin{figure*}
	\centering
	\includegraphics[width=0.49\linewidth]{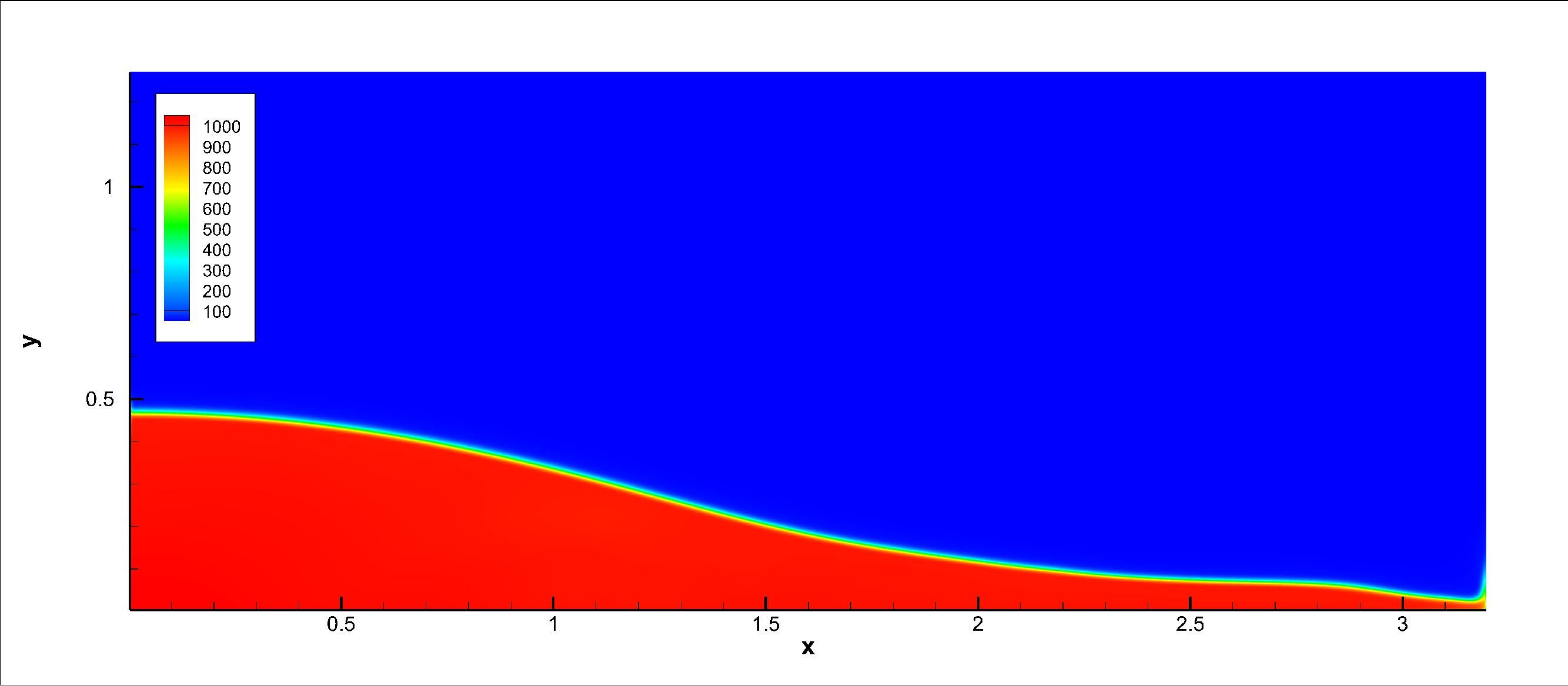}
	\includegraphics[width=0.49\linewidth]{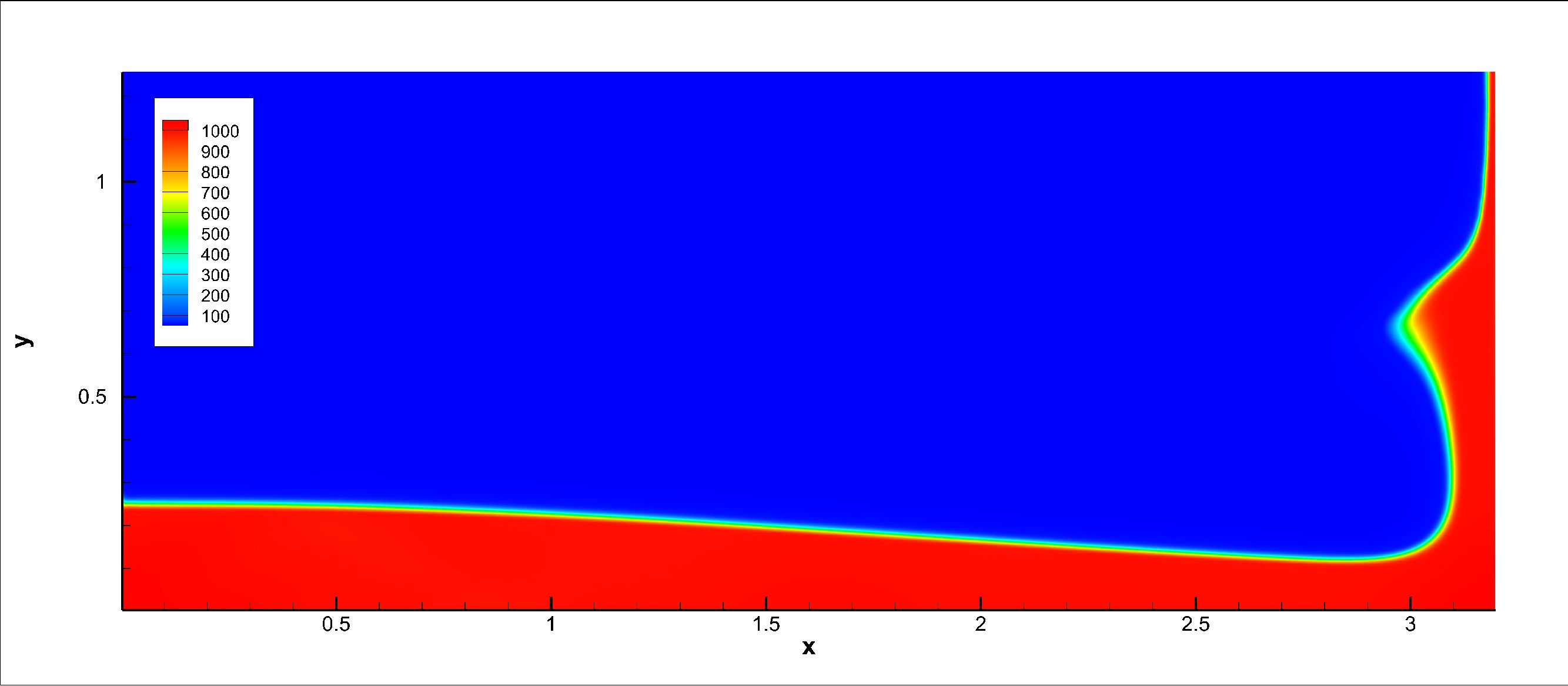} \\ \,
	\includegraphics[width=0.49\linewidth]{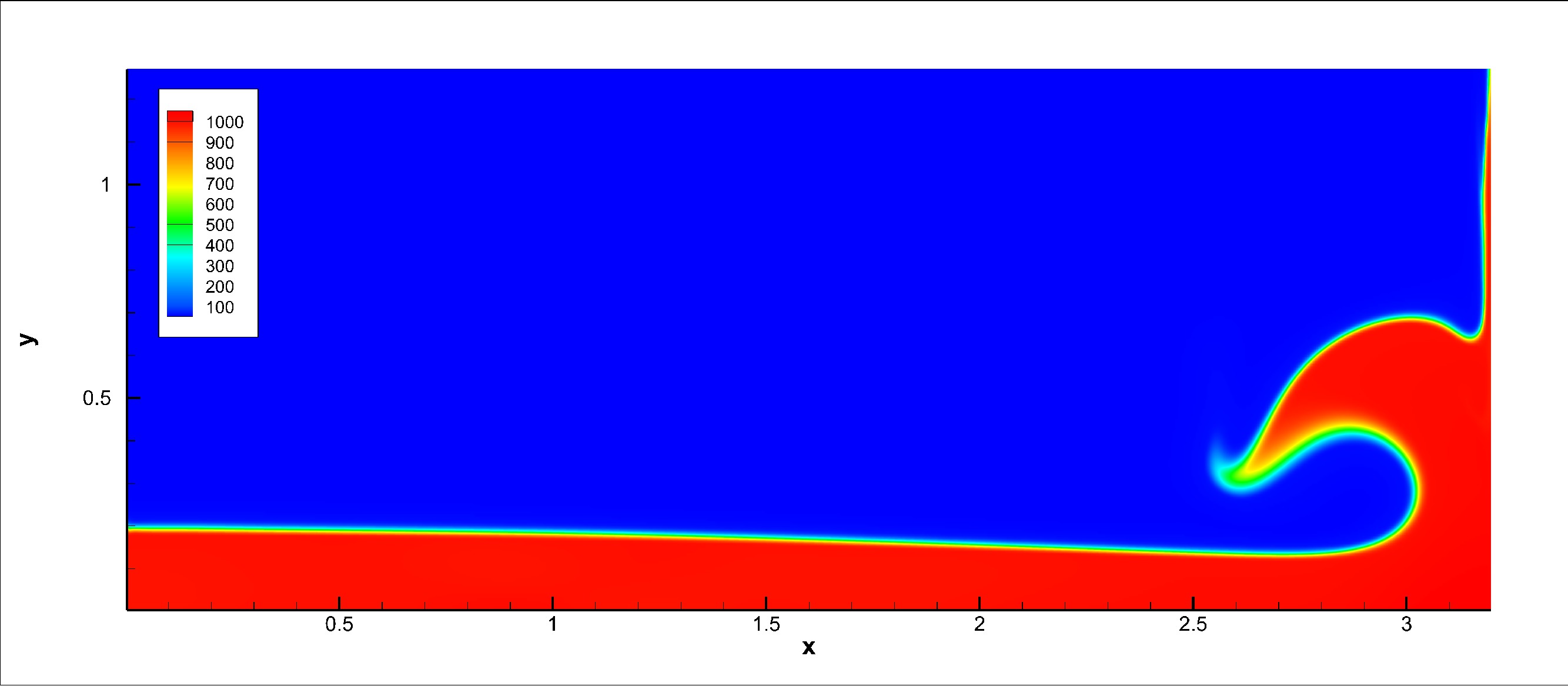}
	\includegraphics[width=0.49\linewidth]{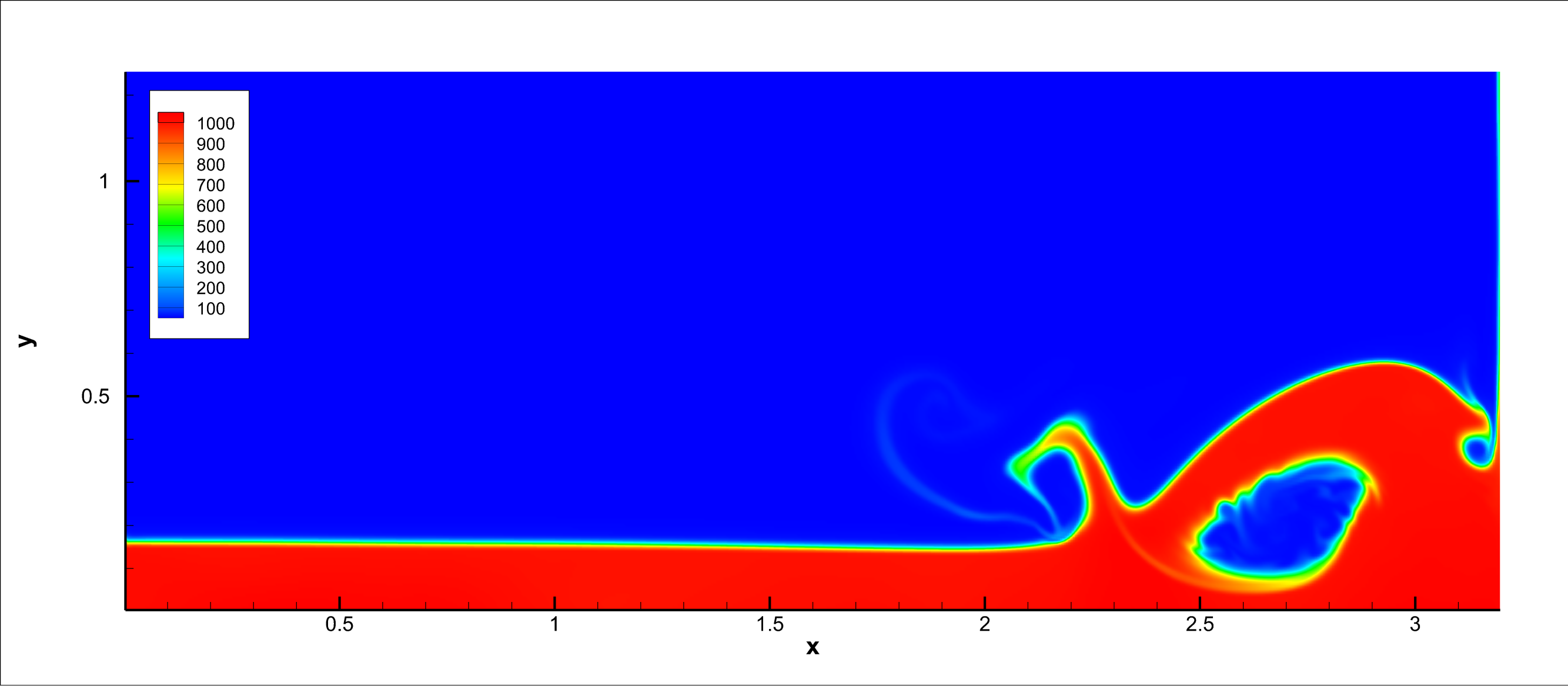} \\ \,
	\caption{Dambreak and impact against a vertical wall. In the figure we depict the results obtained with our second order well balanced scheme for the quantity $\alpha\rho$ at the times $t = 0.6, 1.2, 1.5, 1.7$.}
	\label{fig.My_impact}
\end{figure*}

The last test problem we want to show consists of a dambreak flow with successive impact against a vertical wall, which leads to the reflection of the incident wave and successive wave breaking. The setup of this test case is taken from \cite{SPH3D} and \cite{ColagrossiPhD}.
The computational domain is $\Omega = [0,3.2]\times[0,1.8]$ with reflective boundaries on all the sides apart from the top 
where we impose transmissive boundary conditions.
At $t=0$ the liquid occupies $\Omega_\ell = [0,1.2]\times[0.6]$, and the density distribution is given once again by \eqref{eq.rho_in_liquid}. We set $k_0=2.62\cdot 10^5$ and we cover the domain with a Cartesian mesh of $1000\times500$ elements.
The numerical results obtained with our well balanced second order method are depicted in Figure \ref{fig.My_impact} and can be compared with the previous cited references or with the results obtained in \cite{dumbser2011simple}.
This test shows the capability of our scheme to capture complex free surface flows as those produced by breaking waves. In particular, since the location of the free surface is implicitly given by the volume fraction $\alpha$, it is not constrained 
to be necessarily a a single-valued function, as it is in the standard shallow water context.


\paragraph*{Efficiency}

\begin{figure}
	\begin{center}	
		\includegraphics[width = 0.7\linewidth]{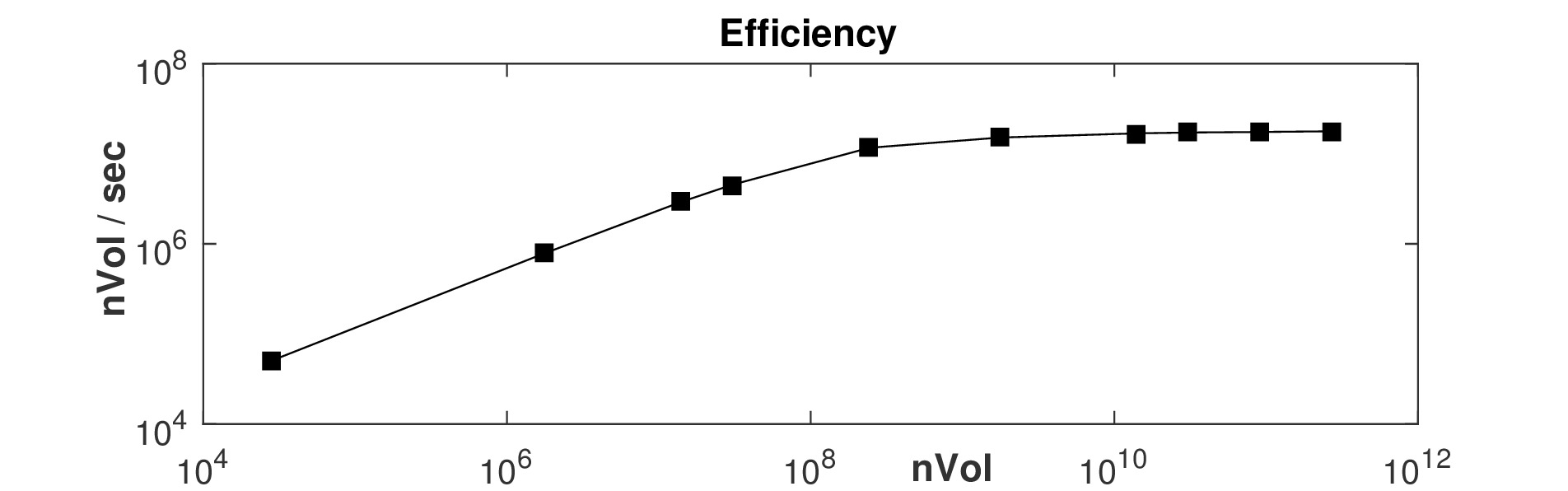}
		\caption{Impact against a vertical wall: efficiency. In the figure we show the number of volumes processed per seconds when employing finer and finer meshes to solve the dambreak problem of Section~\ref{sec.Impacttest} up to a final time $t_f=1$. (Note that with \texttt{nVol} on the x-axis we refer to the number of elements of the mesh multiplied by the number of time steps.)}
		\label{fig.BNsimp.efficiency}   
	\end{center} 
\end{figure}

In order to prove the efficiency of our implementation we have run this last test up to the final time $t_f=1$ on finer and finer meshes using the GeForce Titan Black with the double precision performance option enabled. 
The number of volumes processed per second with respect to the total number of processed space-time control volumes is reported 
in Figure \ref{fig.BNsimp.efficiency}: on sufficiently fine meshes, the new well-balanced scheme implemented in CUDA can update about ten million control volumes per second. 

\section{Conclusions}

In this work we have developed a new second order accurate well balanced path-conservative finite volume scheme for the simulation of complex non-hydrostatic free-surface flows. In particular, this is the first time that a well balanced scheme has been 
applied to the reduced two-phase flow model proposed in \cite{dumbser2011simple}. The mathematical model is derived from the  Baer-Nunziato model of  compressible multi-phase flows and allows to simulate even complex free-surface flows thanks to the use  
of the volume fraction function $\alpha$. The model does not make any of the classical assumptions of shallow water-type systems,  in particular the flow can be fully nonhydrostatic and the free surface is not constrained to be a single-valued function.  

Moreover, the studied system presents some non-conservative terms already in its original form, due to the advection equation for $\alpha$, which requires the use of a path-conservartive scheme. Thus, we have decided to rewrite also the gravity source term 
via a non-conservative product and to apply an Osher-Romberg Riemann solver at the element interfaces: indeed this kind of scheme  allows to treat the non-conservative terms in a well balanced way, so that equilibrium situations can be preserved up to machine  precision even for very long computational times. 
Furthermore, the use of an Osher-type scheme reduces the numerical dissipation at  the free surface.

The obtained results show an excellent agreement with analytical and experimental reference solutions both close and far away from  steady states, thanks also to the automatic detection of the water level of equilibrium, which extends the well balanced properties of the scheme also to situations where the equilibria are not exactly known a priori. 

Further research will concern the extension to more than second order of accuracy and to the three-dimensional case. In particular,  this last enhancement will require the design of an appropriate parallel data and code structure in CUDA in order to maintain the high computational efficiency of the parallel implementation.

\section*{Acknowledgments} 
The research presented in this paper has been partially financed by the European Research Council (ERC) under the 
European Union's Seventh Framework Programme (FP7/2007-2013) with the research project \textit{STiMulUs}, 
ERC Grant agreement no. 278267.
This research has been also supported by the Spanish Government and FEDER
through the research project MTM2015-70490-C2-1-R and the Andalusian Government research projects
P11-FQM-8179 and P11-RNM-7069.
Moreover this project has received funding from the European Union's Horizon 2020 research and innovation Programme under the 
Marie Sklodowska-Curie grant agreement no. 642768.

The authors also acknowledge funding from the Italian Ministry of Education, University and Research (MIUR) in the frame of the 
Departments of  Excellence Initiative 2018--2022 attributed to DICAM of the University of Trento. MD has also received support 
from the University of Trento in the frame of the Strategic Initiative \textit{Modeling and Simulation} (SIMaS).  

\bibliography{references}

\end{document}